\documentclass[12pt]{article}
 
\usepackage{agzt_macros}
\usepackage{amsrefs}

\newcommand{\ogc}{{free odd generation conjecture}\,\,}
\newcommand{\ogcdot}{{free odd generation conjecture.}\,\,}
\newcommand{\ogccomma}{{free odd generation conjecture,}\,\,}
 
\title{AGZT-Lectures on formal multiple zeta values}
\author{A. Burmester, N. Confurius, U. K\"uhn}

\begin{document}
\maketitle

\begin{abstract}
Formal multiple zeta values allow to study multiple zeta values by algebraic methods in a way that the open question about their transcendence is circumvented. In this note we show that Hoffman's basis conjecture for formal multiple zeta values is implied by the \ogc for the double shuffle Lie algebra. 
We use the concept of a post-Lie structure for a convenient approach to the multiplication on the double shuffle group. From this, we get a coaction on the algebra of formal multiple zeta values. 
This in turn allows us to follow the proof of Brown's  celebrated and unconditional theorem for the same result in the context of motivic multiple zeta values. We need the \ogc   twice: at first it gives a formula for the graded dimensions and secondly it is a key to derive a lift of the Zagier formula to the formal context.
\end{abstract}
 
\tableofcontents

\section{Introduction}

Multiple zeta values (MZVs) are real numbers defined as the convergent series
\begin{equation} \label{eq:ser_mzv}
\zeta(k_1, \dotsc, k_d) = \sum \limits_{n_1 > \cdots > n_d > 0} \dfrac{1}{n_1^{k_1} \cdots n_d^{k_d}},
\end{equation}
where $k_i$ are positive integers and the first component $k_1$ is strictly greater than $1$. These values were first considered by Euler in the 18th century and since then they have been studied in 
various contexts in number theory, knot theory and the theory of mixed Tate motives. There is  Hoffman's list of all related publications \cite{hoff_list}.  MZVs form a $\Q$-algebra $\Z$, which is contained in $\R$. 
One of the most challenging open question in the study of MZVs is the identification of all relations among them, even the question whether $\Z$ is graded by the weight is still open.

One of the important properties of MZVs is their representation in terms of iterated integral as follows:
\begin{equation} \label{eq:int_mzv}
\zeta(k_1, \dotsc, k_d) =  \int \limits_{1 > t_1 > \cdots > t_n > 0} \omega_1(t_1) \cdots \omega_n(t_n),
\end{equation}
where $n = k_1 + \cdots + k_d$ is the weight of the MZV, and $\omega_i(t_i) = dt_i/(1 - t_i)$ if $i \in \{k_1, k_1 + k_2, \dotsc, k_1 + \cdots + k_d\}$, and $\omega_i(t_i) = dt_i/t_i$ otherwise. The series representation \eqref{eq:ser_mzv} and the integral representation \eqref{eq:int_mzv} provide two different ways of expanding the product of two MZVs as linear combinations of MZVs, resulting in two distinct combinatorial interpretations. The equality of the products then allows us to generate a large family of relations among MZVs called  double shuffle relations. Nevertheless, these relations are not sufficient to capture all linear relations, for instance, the well-known identity $\zeta(2,1) = \zeta(3)$ due to Euler cannot be derived from them. In order to remedy this, Ihara, Kaneko, and Zagier extended the double shuffle relations by appropriate regularisations $\zeta_\ast(k_1, \dotsc, k_d)$ and $\zeta_\shuffle(k_1, \dotsc, k_d)$ for the divergent series and integrals respectivly. A comparision theorem for these two regularisations allowed them to introduce the so-called  extended double shuffle relations (EDS), which are widely believed to determine all linear relations among MZVs (see Conjecture 1 in \cite{ikz}).

There are two ways to study MZV's algebraically either by means of the formal multiple zeta values or by means of the motivic multiple zeta values. 

The algebra formal multiple zeta values $\Zf$ is the algebra spanned by symbols $\zf(k_1,...,k_d)$, which satisfy exactly the EDS  and no other relations. The work of
Racinet \cite{ra} allows to study  $\Zf$  in the context of Hopf algebras. By construction there is a surjective algebra morphism  $\Zf \to \Z$ given by 
\[
\zf(k_1,...,k_d) \to \zeta_\shuffle(k_1,...,k_d).
\]
The algebra of motivic multiple zeta values $\Z^{\mathfrak{m}}$ introduced and studied intensively by Goncharov, Deligne and Brown 
(see e.g.  \cite{de}, \cite{gon}, \cite{degon}, \cite{br}) 
is a Hopf algebra of functions on a certain group scheme associated to the fundamental group of $\mathbb{P}^1 \backslash \{0,1,\infty\}$. It is spanned by symbols $I^{\mathfrak{m}}(\ve_0;\ve_1,...,\ve_n;\ve_{n+1})$, where $\ve_i \in \{0,1\}$, modulo some relations in such a way that the period map $\Z^{\mathfrak{m}} \to \Z$ given by 
\[
I^{\mathfrak{m}}(\ve_0;\ve_1,...,\ve_n;\ve_{n+1}) \to   \int \limits_{\ve_{n+1} > t_1 > \cdots > t_n > \ve_0} \omega_{\ve_1}(t_1) \cdots \omega_{\ve_n}(t_n)
\]
is a surjective algebra morphism. For more details we refer to the book of Burgos-Gil and Fresan \cite{bgf}. 

Both approaches fit in the following abstract setting:

\begin{equation}  \label{diagram big picture introduction}
\begin{tikzcd}[baseline=(current  bounding  box.center)]
\big(\mathcal{A},\cdot,\Delta_G\big) \arrow[dddd,"\text{mod products}"'] &&&&& \big(\mathcal{U}(\mathfrak{g}),*,\co\big) \arrow[lllll,"\sim","\text{dual}" '] 
\\ \\ 
&& \big(G,*\big) \arrow[uull,leftrightarrow,"1:1"] &&&
\\ \\ 
\big(\operatorname{Indec}(\mathcal{A}) ,\delta\big) &&&&& \big(\mathfrak{g},\{-,-\}\big)
\arrow[lllll, "\sim", "\text{dual}"'] \arrow[uuuu, hookrightarrow] \arrow[uulll,leftrightarrow,"1:1"]
\end{tikzcd}. \end{equation} 

Here $G= \operatorname{Grp}(\widehat{\mathcal{U}}(\mathfrak{g}))$ is a graded, pro-unipotent group scheme, $\mathfrak{g}$ its Lie algebra, $\mathcal{U}(\mathfrak{g})$ the universal enveloping algebra of $\mathfrak{g}$ and the Hopf algebra $\mathcal{A}$ has two descriptions. 
It equals the graded dual of $\mathcal{U}(\mathfrak{g})$ as well as the Hopf algebra of functions on $G$. Finally $\operatorname{Indec}(\mathcal{A})$ is the Lie coalgebra of the indecomposable elements of $\mathcal{A}$.  
In the formal setup, we have $\mathcal{A} = \Zf \slash ( \zf(2) ) $  and Racinet denotes $G$  by $\operatorname{DM}_0$ and the Lie algebra $\mathfrak{g}$ by $\dm$. 
In the motivic setup $\mathcal{A} = \Z^\mathfrak{m} \slash ( \zeta^\mathfrak{m}(2) )$ and $G$ relates to the Galois group of the category of mixed Tate motives.

Explicit calculations show 
\[
\Delta_{\operatorname{DM}_0} = \DeltaGon,
\]
thus for both approaches we have the same formulae for the coproduct.
We like to emphasize the fact that our approach
to the coproduct $\Delta_{\operatorname{DM}_0}$ relies on the general theory of post-Lie algebras together with the work of Racinet, whereas
the original definition of the Goncharov coproduct $\DeltaGon$ in \cite{gon} 
was based on topological considerations for the path algebra. The latter is directly related to the representation of multiple zeta values by iterated integrals, which is a key to motivic multiple zeta values.
A small modification of this coproduct enables us to obtain the first important step in these lecture notes

\begin{Theorem} Set $\Af = \Zf \slash \big( \zf(2)\big)$. There is a well-defined  coaction 
\[ 
\DeltaGon:\Zf\to \Af\otimes\Zf,
\]
which is given by the same formulae as the Brown-Goncharov coaction for 
the motivic multiple zeta values. 
\end{Theorem}

Central for this notes is the following well-known conjecture for $\dm$, which is motivated by conjectures of Deligne (\cite{de}) and Y. Ihara (\cite[p. 300]{ih})
in the context of certain Galois actions and of Drinfeld \cite{dr} on his Grothendieck-Teichm\"uller Lie algebra. By work of Furusho \cite{fu}, we know that the Grothendieck-Teichm\"uller Lie algebra embedds into $\dm$.

\begin{Conjecture} \label{conj:odd_generators_intro} The double shuffle Lie algebra $\dm$ is a free Lie algebra with exactly one generator in each odd weight $w\geq 3$, i. e. 
\begin{equation*}
\dm\simeq\Lie(S),
\end{equation*}
where  the set $S$ is given by $S=\{ s_3, s_5,\ldots\}$.
We call this conjecture the \ogcdot
\end{Conjecture}

The main results we present in this lecture notes are the following.

\begin{Proposition}\label{prop:ogc_implications_intro} 
Assume the \ogc holds for $\dm$, then
\begin{enumerate}
\item Zagier's conjecture holds for  $\Zf$, i.e. 
\[
\sum_{w\ge 0} \dim (\Zfw{w})\, x^w = \frac{1}{1-x^2-x^3},
\]
where $\Zfw{w}$ is the subspace spanned by formal MZVs of weight $w$.
\item The formal zeta values $\zf(2)$ and $\zf(2r+1)$, $r \in \N$, are non-zero modulo products and  algebraically independent. 
\item The Kernel conjecture \ref{conj:kernel} holds for $\Zf$, i.e.  
\[\ker D_{<N}  \cap \Zfw{N}  = \Q \,\zf(N).\]
\end{enumerate}
\end{Proposition}

Using this proposition\footnote{In the motivic setup, the first claim is a theorem of Terasoma and Deligne-Goncharov.  The second and third claim 
are consequences of the construction of motivic multiple zeta values \cite{br},\cite{bgf}. In the long end they rely on Borel's theorem on the algebraic $K$-theory of $\Q$.} it is  not difficult to follow the lines of Brown's proof to derive the following results. 

\begin{Theorem} 
If the  \ogc holds, then all
the $\zf(k_1,\ldots,k_d)$ with $k_i \in \{2,3\}$ are linearly independent.  
\end{Theorem}

The dimension of the space spanned by the formal multiple zeta values from the above theorem in a fixed weight are the same as the ones we expect for all in Zagier's conjecture, therefore we get as corollary the verification of the Hoffman basis conjecture for the formal multiple zeta values.

\begin{Corollary}
If the \ogc holds, then  
the $\zf(k_1,\ldots,k_d)$ with $k_i \in \{2,3\}$ form a basis for $\Zf$ as a vector space.
\end{Corollary}

These notes are based on a series of talks we gave at the
\begin{center}
 \emph{{\bf A}rithmetische {\bf G}eometrie und {\bf Z}ahlen{\bf T}heorie Seminar}
\end{center}
 at the Universit\"at Hamburg  in the summer term of 2023. We like to thanks the audience for helpful remarks, which improved our understanding and this presentation. 
 Our motivation to study Brown's theorem in the context of formal multiple zeta values is that for multiple $q$-zeta values and  for multiple Eisenstein series  similar results either hold or conjecturally hold \cite{BaKu_conj}. Recent progress in that directions can be found in  \cite{bu}, \cite{AB_fqmzv}, 
 \cite{BaIt_fMES}, \cite{BMK_fdMES}.
 
Special thanks also go to Henrik Bachmann, Jose Burgos-Gil, Pierre Lochak, Dominique Manchon, Leila Schneps for various fruitful discussions related to these projects.

\section{Algebraic background}

We provide the general algebraic constructions, which we will use in all following sections.

\subsection{Hopf algebras}
\label{Algebraic background}

We start by a short presentation of Hopf algebras and their behaviour under duality. Detailed introductions into the theory of Hopf algebras can be found in \cite{ca}, \cite{foi}, and \cite{man}. In the following, let $R$ be any commutative ring. 
\begin{Definition}
A \emph{Hopf algebra} over $R$ is a tuple $(H,m,\eta,\Delta,\varepsilon,S)$, where $H$ is an $R$-algebra with the multiplication $m:H\otimes H\to H$ and the unit map $\eta:R\to H$,
\begin{align*} 
\Delta:H\to H\otimes H \qquad  (coproduct),\hspace{2,5cm}
\varepsilon: H\to R \qquad (counit)
\end{align*} 
are $R$-algebra morphisms, and
\[S:H\to H \qquad \qquad (antipode)\]
is a $R$-module morphism, such that the following compatibility conditions hold
\begin{itemize}
\item[(i)] coassociativity: 
\[(\id \otimes \Delta)\circ \Delta=(\Delta\otimes\id)\circ \Delta,\]
\item[(ii)] counitarity:
\[m\circ(\id\otimes\varepsilon)\circ \Delta=m\circ(\varepsilon\otimes\id)\circ\Delta=\id,\]
\item[(iii)] antipode property:
\[m\circ (\id\otimes S)\circ\Delta=m\circ(S\otimes\id)\circ\Delta=\eta\circ\varepsilon.\]
\end{itemize}
\end{Definition}
In the following, we will often omit the unit map, the counit or the antipode, if they are clear from the context or the explicit shape does not matter. 
\begin{Definition} \label{def grading}
We call a Hopf algebra $(H,m,\eta,\Delta,\varepsilon,S)$ \emph{graded} if there is a decomposition
\[H=\bigoplus_{i\in \mathbb{Z}_{\geq0}} H_i,\]
where each $H_i$ is a free $R$-submodule of finite rank, such that
\begin{itemize}
\item[(i)] $m(H_i\otimes H_j)\subset H_{i+j}$ for $i,j\geq0$,
\item[(ii)] $\Delta(H_n)\subset \bigoplus_{i+j=n} H_i\otimes H_j$ for $n\geq0$,
\item[(iii)] $S(H_i)\subset H_i$ for $i\geq0$.
\end{itemize}
\end{Definition} 
In this case, we have
\begin{align*}
\eta(R)\subset H_0, \qquad \varepsilon(H_i)=0\quad \text{for } i\geq1.
\end{align*}
Similarly, we call modules, algebras, and coalgebras graded if they satisfy the corresponding subsets of the above conditions.
 
\begin{Definition} \label{def completion} (i) Let $A$ be an $R$-module equipped with a \emph{descending filtration}, i.e., there is a chain of submodules \[A= \fil^{(0)} A\supset \fil^{(1)} A \supset \fil^{(2)} A \supset \fil^{(3)}A\supset \dots\ .\] 
The \emph{completion} $\widehat{A}$ of $A$ with respect to this filtration is defined by the inverse limit
\[\widehat{A}=\varprojlim_{j}\faktor{A}{\fil^{(j)} A}.\]
If $\widehat{A}=A$, then $A$ is called a \emph{complete} $R$-module.
\end{Definition} 
The completion $\widehat{A}$ of $A$ is also a filtered $R$-module via
\[\fil^{(j)}\widehat{A}=\varprojlim_{k>j} \faktor{\fil^{(j)} A}{\fil^{(k)} A}.\]

\begin{Proposition} \label{completion for grading} Assume that $A=\bigoplus_{i\geq0} A_i$ is a graded $R$-module. Then $A$ admits a descending filtration given by $\fil^{(j)} A= \bigoplus_{i\geq j} A_i$. Since $\faktor{A}{\fil^{(j)} A}=\bigoplus_{i=0}^{j-1} A_i$, the completion of $A$ is 
\[\widehat{A}=\varprojlim_{j} \faktor{A}{\fil^{(j)} A}=\prod_{i\geq 0} A_i.\] The completion $\widehat{A}$ is filtered by $\fil^{(j)}\widehat{A}=\prod_{i\geq j} A_i$. 
\end{Proposition} 

\begin{Definition} \label{def:completed Hofp algebra}
Let $(H,m,\eta,\Delta,\varepsilon,S)$ be a graded Hopf algebra. By extending the maps $m,\eta,\Delta,\varepsilon,S$ of $H$ to the completed module $\widehat{H}$, one obtains \emph{completed Hopf algebra} of $H$.
\end{Definition}

The completed Hopf algebra $(\widehat{H},m,\eta,\Delta,\varepsilon,S)$ is filtered, i.e., one has for all $i\geq0$
\begin{align*}
&m(\fil^{(i)} H\otimes \fil^{(j)} H)\subset \fil^{(i+j)} H,\qquad \Delta(\fil^{(i)}H)\subset \sum_{m+n=i} \fil^{(m)}H\otimes \fil^{(n)}H, \qquad \\ &S(\fil^{(i)}H) \subset \fil^{(i)}H.
\end{align*}
Evidently, we have the same construction for modules, algebras, and coalgebras.

\begin{Definition} \label{def associated graded} Let $A$ be  a filtered $R$-module. Then the \emph{associated graded module} $\gr A$ is defined by
\[\gr A=\bigoplus_{j\geq0} \faktor{\fil^{(j)} A}{\fil^{(j+1)} A}.\]
\end{Definition} 

One has $\gr A=\gr \widehat{A}$. In particular, if $A$ is a graded module, then $\gr \widehat{A}=A$.

If $A$ is a filtered $R$-module and all quotients $\faktor{A}{\fil^{(j)}A}$ are free modules of finite rank, then the module $\gr M$ is graded in the sense of Definition \ref{def grading}.

\begin{Definition} \label{associated graded Hopf algebra} Let $(H,m,\eta,\Delta,\varepsilon,S)$ be a filtered Hopf algebra over $R$. Then the \emph{associated graded Hopf algebra} is the $R$-module $\gr H$ equipped with the induced maps by $m,\eta,\Delta,\varepsilon$ and $S$.
\end{Definition} 

As before, we define the associated graded for modules, algebras, and coalgebras in the same way.

Hopf algebras behave nicely under duality pairings as introduced in \cite[Chapter 2, Section 2.1]{abe}.

\begin{Definition} \label{def dual modules}
Two $R$-modules $A$ and $B$ are \emph{dual}, if there is an $R$-linear map \[(\cdot\mid \cdot):A\otimes B\to R,\] such that 
\begin{itemize} \setlength\itemsep{-1mm}
\item[(i)] if $(a\mid b)=0$ for all $a\in A$, then $b=0$,
\item[(ii)] if $(a\mid b)=0$ for all $b\in B$, then $a=0$.
\end{itemize}
In this case, $(\cdot\mid\cdot)$ is called the \emph{duality pairing} of $A$ and $B$.
\vspace{0,3cm} \\
Let $A$ and $B$ be graded $R$-modules. If there is a duality pairing $(\cdot\mid \cdot):A\otimes B\to R$, such that
\[(A_i\mid B_j)=0 \quad \text{for all } i\neq j,\]
then $A$ and $B$ are \emph{graded dual}. In this case, we say that $(\cdot\mid \cdot)$ is a \emph{graded duality pairing}.
\end{Definition} 
\begin{Example} \label{usual dual spaces}
(i) Let $A$ be a free $R$-module of finite rank. Usually, the dual module is defined by \[A^*=\operatorname{Hom}_{R\operatorname{-lin}}(A,R).\] 
The modules $A$ and $A^*$ are also dual in the sense of Definition \ref{def dual modules}, the duality pairing is given by
\begin{align*}
(\cdot\mid\cdot):A^*\otimes A&\to R, \\
f\otimes a&\mapsto f(a).
\end{align*}
(ii) Let $A$ be a graded $R$-module. Then usually, its graded dual is defined by 	\[A^\vee=\bigoplus_{i\geq 0} A_i^*.\]
The modules $A$ and $A^\vee$ are also graded dual in the sense of Definition \ref{def dual modules}, the graded duality pairing is given by
\begin{align*}
A^\vee\otimes A&\to R,\\
f_i\otimes a_j&\mapsto\begin{cases} f_i(a_j), &\quad i=j, \\ 0 &\quad \text{else}\end{cases} \qquad (\text{where } f_i\in A_i^*,\ a_j\in A_j).
\end{align*}
In all following sections, we will use the notion $(-)^\vee$ exclusively for the graded dual.
\end{Example}

Let $A_1,B_1$ be dual $R$-modules for the pairing $(\cdot\mid\cdot)_1$, $A_2,B_2$ be dual $R$-modules for the pairing $(\cdot\mid\cdot)_2$, and $f:A_1\to A_2$ be an $R$-linear map. The \emph{dual map} to $f$ is the unique $R$-linear map $g:B_2\to B_1$ satisfying
\[(f(a),b)_2=(a,g(b))_1 \qquad \text{ for all } a\in A_1,\ b\in B_2.\]
\begin{Proposition} \label{graded dual Hofp algebras} Let $(H,m,\eta,\Delta,\varepsilon,S)$ be a (graded) Hopf algebra over $R$. If $H'$ is an $R$-module (graded) dual to $H$, then $H'$ equipped with the dual maps of $m,\ \eta,\ \Delta,\ \varepsilon$ and $S$ is also a (graded) Hopf algebra over $R$. \qed
\end{Proposition}

\subsection{Hoffman's quasi-shuffle Hopf algebras}
\label{subsec:quasi-shuffle}

We present a particular class of Hopf algebras, called quasi-shuffle Hopf algebras. Those we first introduced in \cite{h}, \cite{hi}, and all results are taken from there. Let $R$ be a commutative $\Q$-algebra with unit.

\begin{Notation} \label{not:alphabets} Let $\A$ be an alphabet, this means $\A$ is a countable set whose elements are called \emph{letters}. By $R \A$ denote the $R$-module spanned by the letters of $\A$ and let $\nca{R}{\A}$ be the free non-commutative algebra generated by the alphabet $\A$. The monic monomials in $\nca{R}{\A}$ are called \emph{words} with letters in $\A$, the set of all words is $\A^*$. Moreover, we denote by $\one$ the empty word. 

The \emph{length} of a word $w\in \A^*$ equals the number of its letters, i.e., the word $w=a_1\cdots a_n$ with $a_i\in\A$ has length $n$. 
We introduce the $j$-th letter function
\begin{align*}
\ve_j:\A^*&\to \A\cup \{\one\}, \\
a_1\cdots a_n&\mapsto\begin{cases}
a_j \quad & \text{ if } j\leq n, \\
\one & \text{ else}.
\end{cases}
\end{align*}
Instead of $w=\ve_1(w)\ve_2(w)\cdots \ve_n(w)\in \A^*$, we will often just write $w=\ve_1\cdots \ve_n$ where $n$ is the length of $w$. Furthermore, we call $u = \ve_{i_1}(w)\cdots \ve_{i_k}(w)$ a \emph{subword} of $w$ if $i_1<\dots<i_k\leq n$ for some integer $k\geq 0$. A subword $\ve_{i_1}\cdots\ve_{i_k}$ of $w$ is called \emph{strict} if $i_j+1 = i_{j+1}$ for all $j\in\{1,\dots,k-1\}$, i.e., if it consists of consecutive letters of $w$.
\end{Notation}

\begin{Definition} \label{def quasi-shuffle}
Let $\diamond:R \A\times R \A\to R \A$ be a commutative and associative product. Define the \emph{quasi-shuffle product} $\ast_\diamond$ on $\nca{R}{\A}$ recursively by $\one\ast_\diamond w=w\ast_\diamond\one=w$ and 
\begin{align*}
au\ast_\diamond bv=a(u\ast_\diamond bv)+b(au\ast_\diamond v)+(a\diamond b)(u\ast_\diamond v)
\end{align*}
for all $u,v,w\in \nca{R}{\A}$ and $a,b\in \A$.
\end{Definition} 

Note that the quasi-shuffle product $\quasish$ can be equally defined recursively from the left and from the right, since both product expressions agree \cite[Theorem 9]{zud}.

\begin{Example} \label{ex:shuffle}
Define 
\[a\diamond b=0 \qquad \text{for all } a,b\in \A,\] 
then we get the well-known \emph{shuffle product}, which is usually denoted by $\shuffle$. 
\end{Example} 

The \emph{deconcatenation coproduct} $\dec:\nca{R}{\A}\to \nca{R}{\A}\otimes\nca{R}{\A}$ is given for a word $w\in \A^*$ by
\begin{align} \label{def:dec}
\dec(w)=\sum_{w=uv}u\otimes v,
\end{align}
and the corresponding counit $\varepsilon:\nca{R}{\A}\to R$ is given for a word $w\in \A^*$ by
\begin{align*} 
\varepsilon(w)=\begin{cases} 1, \quad & w=\one, \\
0 & \text{else}. \end{cases}
\end{align*}

\begin{Theorem} (\cite[Theorem 3.1,3.2]{h}) The tuple $(\nca{R}{\A},\quasish,\one,\dec,\varepsilon)$ is an associative, commutative Hopf algebra. \qed 
\end{Theorem}

An explicit formula for the antipode of the Hopf algebra $(\nca{R}{\A},\quasish,\one,\dec,\varepsilon)$ is also given in \cite[Theorem 3.2]{h}.

For the shuffle algebra $(\nca{R}{\A},\shuffle)$, c.f. Example \ref{ex:shuffle}, there is an explicit generating set. Choose a total ordering on the alphabet $\A$, then the lexicographic ordering defines a total ordering on the set of all words $\A^*$. 

\begin{Definition} \label{def:Lyndon word} A word $w\in \A^*\backslash\{\one\}$ is called a \emph{Lyndon word} if we have for any non-trivial decomposition $w=uv$ that $w<v$.
\end{Definition}

\begin{Theorem} \label{shuffle algebra Lyndon words} (\cite[Theorem 4.9 (ii)]{Re}) The shuffle algebra $(\nca{R}{\A},\shuffle)$ is a free polynomial algebra generated by the Lyndon words of $\A$. \qed
\end{Theorem}

We will see that all quasi-shuffle algebras over the same alphabet $\A$ are isomorphic. In particular, the previous theorem holds for all quasi-shuffle algebras. \\
Let $(\nca{R}{\A},\quasish)$ be a quasi-shuffle algebra. By a \emph{composition} of a positive integer $n$ we mean an ordered sequence $I=(i_1,\ldots,i_r)$, such that $i_1+\dots+i_r=n$. Let $w=a_1\dots a_n\in \A^*$ be a word and $I=(i_1,\ldots,i_r)$ a composition of $n$, then define
\begin{align*} 
I[w]&=(a_1\diamond\dots\diamond a_{i_1})(a_{i_1+1}\diamond\dots \diamond a_{i_1+i_2})\dots (a_{i_1+\dots+i_{r-1}+1}\diamond\dots\diamond a_n)
\end{align*} 
and
\begin{align*}
\exp_{\quasish}(w)&=\sum_{I=(i_1,\ldots,i_r) \text{ composition of } n} \frac{1}{i_1!\cdots i_r!} I[w], \\
\log_{\quasish}(w)&=\sum_{I=(i_1,\ldots,i_r) \text{ composition of } n} \frac{(-1)^{n-r}}{i_1\cdots i_r} I[w].
\end{align*}
\begin{Theorem} \label{Hoffman iso} (\cite[Theorem 3.3]{h}) The map $\exp_{\quasish}$ is a Hopf algebra isomorphism
\[\exp_{\quasish}: (\nca{R}{\A},\shuffle,\dec)\xrightarrow{\sim} (\nca{R}{\A}, \quasish,\dec). \]
The inverse map is given by $\log_{\quasish}$. \qed
\end{Theorem} 

From Theorem \ref{shuffle algebra Lyndon words} and \ref{Hoffman iso}, one deduces the following.

\begin{Corollary} \label{qsa generated by Lnydon words}
Any quasi-shuffle algebra $(\nca{R}{\A},\quasish)$ is a free polynomial algebra generated by the Lyndon words of $\A$.
\end{Corollary}

We want to determine a dual of the quasi-shuffle Hopf algebra. Define a degree map on the letters in $\A$, such that $\deg(a)\geq1$ for all $a\in \A$. This induces a grading on $\nca{R}{\A}$ by
\[\deg(a_1\cdots a_n)=\deg(a_1)+\dots+\deg(a_n),\qquad a_1,\ldots,a_n\in\A.\] Denote by $\ncac{R}{\A}$ the completion with respect to this grading. There is a duality pairing
\begin{align} \label{pairing quasi-shuffle}
(\cdot\mid \cdot): \ncac{R}{\A}\otimes\nca{R}{\A}&\to R, \\
\phi\otimes w&\mapsto (\phi\mid w), \nonumber
\end{align}
where $(\phi\mid w)$ denotes the coefficient of $\phi\in \ncac{R}{\A}$ in $w\in \nca{R}{\A}$. We assume that $(\nca{R}{\A},\quasish,\dec)$ is a graded Hopf algebra with respect to the above degree map (Definition \ref{def grading}). Then the dual coproduct $\Delta_{\quasish}:\ncac{R}{\A}\to\ncac{R}{\A}\otimes\ncac{R}{\A}$ to $\quasish$ with respect to the above duality pairing is given by
\begin{align} \label{quasish dual coproduct}
\Delta_{\quasish}(\phi)=\sum_{u,v\in \A^*} (\phi\mid u\quasish v)\ u\otimes v.
\end{align}
Since the quasi-shuffle product $\quasish$ is graded and the homogeneous subspaces of $\nca{R}{\A}$ are finite dimensional, each coefficient in the coproduct is finite.
Moreover, denote the concatenation product by $\conc$.

\begin{Theorem} \label{dual to quasi-shuffle} The tuple $(\ncac{R}{\A},\conc,\one,\Delta_{\quasish},\varepsilon)$ is a complete cocommutative Hopf algebra. It is dual to the quasi-shuffle Hopf algebra $(\nca{R}{\A},\quasish,\one,\dec,\varepsilon)$ with respect to the pairing $(\cdot\mid\cdot)$ given in \eqref{pairing quasi-shuffle}. \qed
\end{Theorem} 

An explicit formula for the antipode of the Hopf algebra $(\ncac{R}{\A},\conc,\one,\Delta_{\quasish},\varepsilon)$ is given in \cite[p. 9]{h}.

\begin{Example} \label{duality shuffle product}
For the shuffle product $\shuffle$ given in Example \ref{ex:shuffle}, the dual coproduct is given by \[\co(a)=a\otimes\one+\one\otimes a \qquad \text{for all } a\in \A.\] So $(\ncac{R}{\A},\conc,\co)$ is a cocommutative Hopf algebra. 
\end{Example} 

\begin{Remark} \label{rem:dual pairing vs graded dual}
The associated graded Hopf algebra (cf Definition \ref{associated graded Hopf algebra}) to the completed Hopf algebra $(\ncac{R}{\A},\conc,\one,\Delta_{\quasish},\varepsilon)$ is just $(\nca{R}{\A},\conc,\one,\Delta_{\quasish},\varepsilon)$. By construction, $(\nca{R}{\A},\conc,\one,\Delta_{\quasish},\varepsilon)$ is exactly the graded dual to $(\nca{R}{\A},\quasish,\one,\dec,\varepsilon)$ in the sense of Example \ref{usual dual spaces}.
\end{Remark}

\subsection{The interaction of Hopf and Lie algebras}

We review some basic results on the interplay of Hopf algebras, group schemes, and Lie algebras, which will be applied in the following sections. 

We start with a basic example for a Hopf algebra, which occurs many times in the following. Let $R$ be any fixed commutative $\QQ$-algebra with unit.

\begin{Definition} Let $(\mathfrak{g},[-,-])$ be a Lie algebra over $R$. Then the \emph{universal enveloping algebra} of $\mathfrak{g}$ is
\[\mathcal{U}(\mathfrak{g})=\faktor{\mathcal{T}(\mathfrak{g})}{\langle x\otimes y-y\otimes x-[x,y]\mid x,y\in \mathfrak{g}\rangle}, \]
where $\mathcal{T}(\mathfrak{g})=\bigoplus_{j\geq0}\mathfrak{g}^{\otimes j}$ is the \emph{tensor algebra}. The product on $\mathcal{U}(\mathfrak{g})$ is induced by the concatenation product on the tensor algebra $\mathcal{T}(\mathfrak{g})$. 
\end{Definition}

The universal enveloping algebra of a Lie algebra is naturally equipped with a Hopf algebra structure. Let $(\mathfrak{g},[-,-])$ be any Lie algebra. Then the tensor algebra $\mathcal{T}(\mathfrak{g})$ is a Hopf algebra with the coproduct $\Delta:\mathcal{T}(\g)\to\mathcal{T}(\g)\otimes\mathcal{T}(\g)$ given by
\[\Delta(x)=x\otimes 1+1\otimes x,\qquad x\in \mathfrak{g},\]
the counit $\varepsilon:\mathcal{T}(\g)\to R$ given by 
\[\varepsilon(1)=1, \qquad \varepsilon(x)=0 \quad \text{ for } x\in \g,\] 
and the antipode $S:\mathcal{T}(\g)\to \mathcal{T}(\g)$ given by 
\[S(x_1\cdots x_n)=(-1)^nx_n\cdots x_1, \qquad  x_1,\ldots,x_n\in\g.\]
Since $\langle x\otimes y-y\otimes x-[x,y]\mid x,y\in \mathfrak{g}\rangle$ is a Hopf ideal in $\mathcal{T}(\mathfrak{g})$, also $\mathcal{U}(\mathfrak{g})$ becomes a Hopf algebra with the induced coproduct, counit, and antipode.

Let $(\g,[-,-])$ be a graded Lie algebra with $\operatorname{rank} L_0=0$. Then, the universal enveloping algebra $\mathcal{U}(\g)$ admits a grading with respect to the Hopf algebra structure. The graded dual of $\mathcal{U}(\g)$ is the symmetric algebra of the graded dual $\g^\vee$
\[\mathcal{S}(\g^\vee)=\faktor{\mathcal{T}(\g^\vee)}{\langle x\otimes y-y\otimes x\mid x,y\in \g^\vee\rangle}, \]
this means we have an algebra isomorphism 
\begin{align} \label{eq:iso_US}
\mathcal{U}(\g)^\vee\simeq \mathcal{S}(\g^\vee).
\end{align}
It is well-known that symmetric algebras are free polynomial algebras.

\begin{Example}
Consider the free Lie algebra $\Lie_R(\A)$ over some alphabet $\A$. Then, the universal enveloping algebra is given by
\[\mathcal{U}(\Lie_R(\A))=(\nca{R}{\A}, \conc,\Delta_{\shuffle}),\]
and thus the graded dual is (cf. Example \ref{duality shuffle product})
\[\mathcal{U}(\Lie(\A))^\vee=(\nca{R}{\A},\shuffle,\dec).\]
So by \eqref{eq:iso_US}, there is an algebra isomorphism
\[(\nca{R}{\A},\shuffle)\simeq (\mathcal{S}(\Lie_R(\A)),\cdot).\]
Moreover, $\mathcal{U}(\Lie(\A))$ is dual to $(\ncac{R}{\A},\shuffle,\dec)$ under the duality pairing given in \eqref{pairing quasi-shuffle} (cf Remark \ref{rem:dual pairing vs graded dual}).
\end{Example}

We review grouplike, primitive and indecomposable elements of Hopf algebras and their relationship. For this, we fix a commutative ring $R$ and a Hopf algebra $(H,m,\eta,\Delta,\varepsilon,S)$ over $R$.

\begin{Definition} 
An element $x\in H\backslash\{0\}$ is  called \emph{grouplike} if
\[\Delta(x)=x\otimes x.\]
The set of grouplike elements in $H$ is denoted by $\Grp(H)$. An element $x\in H$ is called \emph{primitive} if it satisfies
\[\Delta(x)=x\otimes 1+1 \otimes x.\]
By $\Prim(H)$ we denote the set of all primitive elements in $H$.
\end{Definition}

\begin{Theorem} \label{Thm grouplike primitive} 
The following holds.
\begin{itemize}
\item[(i)] The set $\Grp(H)$ equipped with the product and the unit of $H$ forms a group. For an element $x\in \Grp(H)$, the inverse element is given by $S(x)$.
Moreover, each grouplike element $x\in H$ satisfies $\varepsilon(x)=1$.
\item[(ii)] The set $\operatorname{Prim}(H)$ equipped with the commutator bracket $[x,y]=m(x\otimes y)-m(y\otimes x)$ is a Lie algebra. Furthermore, one has for each primitive element $x\in H$ that $\varepsilon(x)=0$ and $S(x)=-x$.
\end{itemize} 
\vspace{-0,5cm} \qed 
\end{Theorem}

\begin{Theorem} \label{CQMM theorem} (Cartier-Quillen-Milnor-Moore) Let $H$ be a graded cocommutative Hopf algebra, such that $\operatorname{rank} H_0=1$. Then there is a Hopf algebra isomorphism
\[H\simeq \mathcal{U}(\Prim(H)).\]
\vspace{-1cm} \\ \qed 
\end{Theorem} 

By passing to completions, we are able to relate the grouplike and primitive elements via an exponential map.
Let $H=\bigoplus_{j\geq0} H_j$ be a graded Hopf algebra and $\widehat{H}=\prod_{j\geq0} H_j$ its completion (Definition \ref{rem:dual pairing vs graded dual}). For any element $x\in \prod_{j\geq 1} H_j\subset \widehat{H}$, define
\[\exp_H(x)=\sum_{i\geq0} \frac{1}{i!}x^i, \]
where $x^i$ means applying the product map of $H$ exactly $(i-1)$-times to $x^{\otimes i}$.

\begin{Proposition} \label{log/exp for Prim/G} Let $H$ be a graded  Hopf algebra. Then there is a bijection
\begin{align*}
\Prim(\widehat{H}) &\xrightarrow{\sim} \Grp(\widehat{H}), \\
x &\mapsto \exp_H(x).
\end{align*}
\vspace{-1cm} \\ \qed 
\end{Proposition}

\begin{Definition} \label{def indecomposables}
The space of \emph{indecomposables} of $H$ is defined as
\[\Indec(H)=\faktor{\ker(\varepsilon)}{\ker(\varepsilon)^2}.\]
\end{Definition} 

Recall that if $H$ is a graded Hopf algebra with $H=\bigoplus_{j\geq0} H_j$, then we have
\[\ker(\varepsilon)=\bigoplus_{j\geq1} H_j.\]

Define the corresponding Lie cobracket $\delta$ to the coproduct $\Delta$ of $H$ as
\begin{align} \label{corresponding lie cobracket}
\delta=(\id-t)\circ \Delta:H\to H\otimes H,
\end{align}
where $t:H\otimes H\to H\otimes H$ simply permutes the tensor product factors.

\begin{Proposition} \label{indecomposables Lie cobracket}
The Lie cobracket $\delta$  defined in \eqref{corresponding lie cobracket} induces a Lie coalgebra structure on the space $\Indec(H)$ of indecomposables. \qed
\end{Proposition}

This Lie coalgebra structure is closely related to the Lie algebra structure on the primitive elements.

\begin{Theorem} \label{iso prim indecomposables} Let $H$ be a graded Hopf algebra, then also $\Indec(H)$ is a graded Lie coalgebra. There is an isomorphism of graded Lie algebras 
\[\Prim(H^\vee)\simeq \Indec(H)^\vee.\] \vspace{-1cm} \\ \qed
\end{Theorem} 

Here $H^\vee$ denotes the graded dual Hopf algebra of $H$ and $\Indec(H)^\vee$ denotes the graded dual Lie algebra of the Lie coalgebra $\Indec(H)$ (as in Example \ref{usual dual spaces}). 

Next, we shortly explain the relationship of affine group schemes and Hopf algebras. A detailed exposition of 
this interplay of algebraic geometry and abstract algebra is given in \cite{dg} and \cite{wa}. We also explain how grouplike, primitive, and indecomposable elements occur in this context.

Let $\QAlg$ be the category of commutative $\QQ$-algebras with unit, $\Sets$ be the category of sets, and $\Groups$ be the category of groups. 

\begin{Definition}
A functor $F: \QAlg\to\Sets$ is an \emph{affine scheme} if there is an object $A\in\QAlg$, such that $F$ is naturally isomorphic to the Hom-functor
\begin{align*}
\operatorname{Hom}_{\QAlg}(A,-): \QAlg&\to\Sets, \\
B&\mapsto\operatorname{Hom}_{\QAlg}(A,B).
\end{align*}
In this case, one says that $A$ represents the functor $F$.
\end{Definition} 

\begin{Theorem} \label{Yoneda} (Yoneda's Lemma) Let $E,F:\QAlg\to\Sets$ be two affine schemes represented by $A,B$. Then any natural transformation $\Phi:E\to F$ corresponds uniquely to an algebra morphism $\varphi:B\to A$. \qed
\end{Theorem} 

\begin{Definition} A functor $G:\QAlg\to \Groups$ is an \emph{affine group scheme} if there is $A\in\QAlg$, such that $G$ is naturally isomorphic to the Hom-functor $\operatorname{Hom}_{\QAlg}(A,-): \QAlg\to\Groups$.
\end{Definition} 

Any affine group scheme is also an affine scheme, we simply ignore the additional group structure.

\begin{Theorem} \label{affine group schemes Hopf algebras} (\cite[Subsection 1.4]{wa}) Let $F:\QAlg\to\Sets$ be an affine scheme represented by $A$. Then $F$ is an affine group scheme if and only if $A$ is a commutative Hopf algebra.
\qed \end{Theorem}

\begin{Example} \label{grou-like elements ags} For each $R\in\QAlg$, consider the dual quasi-shuffle Hopf algebra $(\ncac{R}{\A},\conc,\Delta_{\quasish})$ from Theorem \ref{dual to quasi-shuffle}. The functor
\begin{align*}
F:\QAlg&\to\Sets, \\
R &\mapsto\ncac{R}{\A}
\end{align*}
is an affine scheme represented by the commutative polynomial algebra $\QQ[(z_w)_{w\in \A^*}]$. The grouplike elements $\Grp(\ncac{R}{\A})$ for the coproduct $\Delta_{\quasish}$ form a group with the concatenation product (Theorem \ref{Thm grouplike primitive}). Hence restricting the images of the affine scheme $F$ to the grouplike elements $\Grp(\ncac{R}{\A})$, one obtains an affine group scheme
\begin{align*}
G:\QAlg&\to\Groups, \\
R &\mapsto \Grp(\ncac{R}{\A}).
\end{align*}
The affine group scheme $G$ is represented by the Hopf algebra $(\Qnca{\A},\quasish,\dec)$.
\end{Example}

Similar to the connection of Lie groups and Lie algebras, one can assign a Lie algebra functor to each affine group scheme.

\begin{Definition} \label{def Lie algebra for ags}
For $R\in\QQ\operatorname{-Alg}$, let $R[\varepsilon]=R[t] \slash (t^2)$ 
be the algebra of dual numbers, so $\varepsilon^2=0$. For an affine group scheme $G$, the corresponding \emph{Lie algebra functor} is
\begin{align*}
\g:\QAlg&\to \Liealg, \\
R&\mapsto \ker\Big(G\big(R[\varepsilon]\twoheadrightarrow R\big)\Big).
\end{align*}
\end{Definition}

Relating $\g(R)$ to the derivations on the representing Hopf algebra of $G$, which are left-invariant under the coproduct, gives the Lie algebra structure on $\g(R)$. 

The space $\g(R)$ consists of all elements in $G(R[\varepsilon])$ of the form $1+\varepsilon x$. Thus, one often identifies
\[\g(R)=\{x\mid 1+\varepsilon x\in G(R[\varepsilon]) \}.\]

Every affine group scheme $G$ is an inverse limit of algebraic affine group schemes, i.e., we have
\[G=\varprojlim G_n,\]
where each $G_n$ is represented by a finite dimensional Hopf algebra over $\QQ$.
Hence, we also have for the Lie algebra functor $\mathfrak{g}$ of $G$ that
\[\mathfrak{g}(R)=\varprojlim \mathfrak{g}_n(R), \qquad R\in \QAlg.\]
So, $\mathfrak{g}(R)$ is a completed filtered Lie algebra, where $\fil^{(j)}\mathfrak{g}(R)$ constists of all elements whose projection to $\mathfrak{g}_0(R),\ldots, \mathfrak{g}_{j-1}(R)$ is zero.

\begin{Proposition} \label{duality Lie algebra to ags and indecomposables}
Let $G$ be an affine group scheme represented by a graded Hopf algebra $H$ and denote by $\g$ the Lie algebra functor to $G$. Then one has
\[\gr\g(\QQ)\simeq\Indec(H)^\vee.\]
\vspace{-1cm} \\ \qed
\end{Proposition} 

By the above disscussion, $\mathfrak{g}(\QQ)$ is filtered, hence we can apply the construction from Definition \ref{def associated graded} to obtain the graded Lie algebra $\gr \mathfrak{g}(\QQ)$.

\begin{Proposition} \label{Lie algebra functor affine scheme}
Let $G$ be an affine group scheme. Then the Lie algebra functor $\g:\QAlg\to \Liealg$ is an affine scheme represented by $\mathcal{U}(\gr \g(\QQ))^\vee$. \qed
\end{Proposition} 

There is an important class of affine group schemes, for which there exists a natural isomorphism to their Lie algebra functors. These group schemes are called \emph{pro-unipotent}, a detailed discussion suitable for our context  can be for example found in \cite[Appendix A.6]{bu}, see also \cite{bgf}.

\begin{Theorem} \label{exp is iso g to G} (\cite[IV, Proposition 4.1]{dg})
Let $G$ be a pro-unipotent affine group scheme with Lie algebra functor $\g$. Then there is a natural isomorphism
\[\exp:\g\to G.\]
\vspace{-1cm} \\ \qed
\end{Theorem} 

The Baker-Campbell-Hausdorff series (\cite[p. 260]{mil}) gives the explicit relation between the Lie bracket on $\g$ and the group multiplication on $G$ under the isomorphism $\exp$. 

\begin{Example} \label{ex ags}
In Example \ref{grou-like elements ags} we considered the dual quasi-shuffle Hopf algebra $(\ncac{R}{\A},\conc,\Delta_{\quasish})$ and obtained the corresponding affine group scheme
\begin{align*}
G:\QAlg&\to\Groups, \\
R&\mapsto \Grp(\ncac{R}{\A}).
\end{align*}
The corresponding Lie algebra functor is given by
\begin{align*}
\g:\QAlg&\to \Liealg, \\
R&\mapsto \Prim(\ncac{R}{\A}),
\end{align*}
where we mean the primitive elements for the coproduct $\Delta_{\quasish}$. The Lie bracket is simply the commutator with respect to concatenation (cf. Theorem \ref{Thm grouplike primitive}).

The affine group scheme $G$ is pro-unipotent. So by Theorem \ref{exp is iso g to G}, there is a natural isomorphism \[\exp:\g\to G.\] 
Explicitly, this isomorphism is given by (cf. Theorem \ref{log/exp for Prim/G})
\begin{align*}
\exp(R):\Prim(\ncac{R}{\A})&\to\Grp(\ncac{R}{\A}), \\
f&\mapsto \exp(R)(f)=\sum_{i\geq 0}\frac{1}{i!}f^i.
\end{align*}
\end{Example} 

We summarize the results from this subsection in a diagram. Let $(G,\cdot)$ be a pro-unipotent affine group scheme, such that the representing Hopf algebra $(H,m_H,\Delta_H)$ is graded, commutative and satisfies $\operatorname{rank} H_0=1$. Moreover, let $\g$ be the Lie algebra functor associated to $G$. Then there is the following diagram 
\begin{equation} \label{diagram big picture general}
\begin{tikzcd}[baseline=(current  bounding  box.center)]
(H,m_H,\Delta_H) \arrow[dddd,"\text{Prop \ref{indecomposables Lie cobracket}}"'] &&&&&& (\mathcal{U}(\gr\g(\QQ)),\cdot,\Delta) \arrow[llllll,"\sim","\text{graded dual \eqref{duality H and U(g)}}" '] 
\\ \\ 
&&& (G,\cdot) \arrow[uulll,leftrightarrow,"1:1", "\text{Thm } \ref{affine group schemes Hopf algebras}"'] &&&
\\ \\ 
(\Indec(H),\delta) &&&&&& (\gr\g(\QQ),[-,-])
 \arrow[llllll, "\sim", "\text{graded dual (Prop \ref{duality Lie algebra to ags and indecomposables})}"'] \arrow[uuuu, hookrightarrow] \arrow[uulll,leftrightarrow,"1:1","\begin{matrix}\hspace{-1,5cm}\exp/\log\ +\ \gr \\ \text{ (Thm }\ref{exp is iso g to G})\end{matrix}"']
\end{tikzcd} \end{equation}
The upper duality is obtained from Theorem \ref{CQMM theorem}, \ref{iso prim indecomposables}, and Proposition \ref{duality Lie algebra to ags and indecomposables}
\begin{align} \label{duality H and U(g)}
H^\vee\simeq \mathcal{U}(\Prim(H^\vee))\simeq \mathcal{U}(\Indec(H)^\vee)\simeq\mathcal{U}(\gr\g(\QQ)).
\end{align}

\subsection{Derivations from coproducts} \label{subsec:hopf_deriva}

We fix a commutative $\QQ$-algebra $R$ with unit and let $(H,\ast,\one,\Delta,\varepsilon)$ be a graded Hopf algebra  over $R$ satisfying $\operatorname{rank} H_0=1$. Recall that we have \[\ker(\varepsilon)=\bigoplus_{w>0} H_w,\] so the space of indecomposables (Definition \ref{def indecomposables})
\[\Indec(H)=\ker(\varepsilon)/\ker(\varepsilon)^2\]  
inherits the grading. For each positive degree $w>0$, we get a canonical projection
\begin{align} \label{def:proj_pi} 
\pi_w:\ker(\varepsilon)\to \Indec(H)_w.
\end{align} 

\begin{Definition} \label{def:deriv_H}
For each $w \in\N$ define the map $D_{w}$ via the composition
\[
\xymatrix{
D_{w}: & H \ar[rr]^{\hspace{-1.5em}\Delta'} && \ker(\varepsilon)\otimes H \ar[rr]^{\hspace{-0.5em} \pi_{w}\otimes\id} && \Indec(H)_{w}\otimes H,
}\]
where we set $\Delta' = \Delta - \one\otimes\id$.  
\end{Definition}

Observe that we have by \cite[Proposition II.1.1]{man} \[\operatorname{im}(\Delta')\subset \ker(\varepsilon)\otimes H=\bigoplus_{w>0} H_w\otimes H.\]

\begin{Lemma} \label{lem:D_on_H_deriv}
	For each $w \in \NN$, the map $D_{w}\colon H \to \Indec(H)_{w} \otimes H$ from Definition \ref{def:deriv_H} is a derivation with respect to $\ast$, i.e., we have
	\[D_{w}(u\ast v)=D_{w}(u)\ast (\one\otimes v) + (\one\otimes u) \ast D_{w}(v).\]
\end{Lemma}

\begin{proof}
	For $u,v\in H$, we compute
	\begin{align*}
		D_{w}(u\ast v) &=(\pi_{w}\otimes\id)\Big(\Delta(u\ast v)-\one\otimes (u\ast v)\Big) \\
		&= (\pi_{w}\otimes\id)\Big(\Delta(u)\ast \Delta(v)-\one\otimes (u\ast v)\Big) \\
		&=(\pi_{w}\otimes\id)\Big(\Delta'(u)\ast \Delta'(v)+(\one\otimes u)\ast\Delta'(v)\\
		&\hspace{6cm}+\Delta'(u)\ast(\one\otimes v)\Big) 
	\end{align*} 
	All left tensor product factors in $\Delta'(u)$ and $\Delta'(v)$ are contained in $\ker(\varepsilon)$ by construction, hence the left tensor product factors of $\Delta'(u)\ast\Delta'(v)$ are in $\ker(\varepsilon)^2$. By definition of the projection $\pi_{w}$, we deduce 
	\begin{align*}
		D_{w}(u\ast v)&=(\pi_{w}\otimes\id)\Big((\one\otimes u)\ast\Delta'(v)+\Delta'(u)\ast(\one\otimes v)\Big) \\
		&=(\one\otimes u)\ast (\pi_{w}\otimes\id)\circ \Delta'(v)+(\pi_{w}\otimes\id)\circ \Delta'(u)\ast(\one\otimes v) \\
		&=D_{w}(u)\ast (\one\otimes v) + (\one\otimes u) \ast D_{w}(v).
	\end{align*} 
	
\end{proof}

\begin{Example}\label{ex:odd_free_lie} 
We endow the set $ S=\{ s_3,s_5 ,..., s_{2n+1},...\}$ with the weight function given by $\wt(s_{w}) =w$. This induces a grading on the Hopf algebra 
\[\mathcal{U}(\Lie(S))=\big(\Q\langle S\rangle, \conc,\co \big).\]
We write  
\begin{equation*}
\CUdual = (\Q\langle S\rangle, \shuffle, \dec).
\end{equation*}
for its graded dual (Example \ref{usual dual spaces}, \ref{duality shuffle product}). Since $\CUdual$ is a  graded Hopf algebra  we can study the derivations $D_w$ of Lemma \ref {lem:D_on_H_deriv}. For simplicity, we denote
\[L=\Indec(\CUdual).\]
Consider the element $s_3 s_5 s_7\in\CUdual$ of weight $15$. We have
	\begin{align*}
		\dec'(s_3 s_5 s_7) 
		&= \dec(s_3 s_5 s_7) - \one\otimes s_3 s_5 s_7\\
		&= s_3 s_5 s_7 \otimes \one + s_3 s_5 \otimes s_7 + s_3 \otimes s_5 s_7 \in \bigoplus_{w\geq1} \CUdual_{w} \otimes \CUdual.
	\end{align*}
	Observe, all the left factors are Lyndon words, thus they are non zero modulo shuffle products. Hence, we get the following non-zero images in $L_w\otimes \CUdual$
	\begin{align*}
		D_3(s_3s_5s_7) &= \pi_3(s_3 ) \otimes s_5s_7, \\
		D_{8}(s_3s_5s_7) &= \pi_{8}(s_3 s_5 )\otimes s_7,\\
		D_{15}(s_3s_5s_7) &= \pi_{15}(s_3s_5s_7)\otimes\one.
	\end{align*}
\end{Example}

\begin{Definition} \label{def:dec_Uf_general} Given the Hopf algebra
	$
	\CUdual = (\Q\langle S\rangle, \shuffle, \dec)
	$ from Example \ref{ex:odd_free_lie} 
	we define $\CUf$ by
	\begin{equation*}
		\CUf = \CUdual\otimes \Q[s_2].
	\end{equation*}
	and we extend the deconcatenation coproduct on $\CUdual$ to a coaction 
	\[\dec: \CUf\to \CUdual\otimes \CUf\]
	via $\dec(s_2) = 1\otimes s_2$. 
	We also set $s_{2n} = b_n s_2^n$, where the $b_n$'s  are non-zero rational numbers. 
\end{Definition}

We may assume that these $b_n$  are the Bernoulli numbers given by Proposition \ref{prop:zeta_222}.
The derivation $D$ extends naturally for this coaction. 

\begin{Example} For $r\in\N$, we compute $D_{2r+1}(s_{2r+1} s_2)$. First, we have
	\begin{align*}
		\dec'(s_{2r+1} s_2) &= \dec(s_{2r+1} s_2) - \one\otimes s_{2r+1} s_2 \\
		&= (s_{2r+1}\otimes \one + \one \otimes s_{2r+1}) (\one\otimes s_2) - \one\otimes s_{2r+1} s_2 \\
		&= s_{2r+1} \otimes s_2
		\in \bigoplus_{w\geq1}\CUdual_{w} \otimes \CUf. \end{align*}
	So finally we obtain
	\begin{equation*}
		D_{2r+1}(s_{2r+1} s_2) = \pi_{2r+1}(s_{2r+1})\otimes s_2 \in L_{2r+1}\otimes\CUf.
	\end{equation*}
\end{Example}

\begin{Definition}\label{def:sum_of_derivations} Keep the notation of Definition \ref{def:deriv_H}. 
	For each $N\in\NN$, we define
	\begin{equation*}
		D_{<N} = \bigoplus_{1\le2r+1 <N} D_{2r+1}.
	\end{equation*}
\end{Definition}

\begin{Example}\label{ex:s_even_in_ker}
	If $N=2n$ is even, then $\dec'(s_N) = 0$. Indeed using $b_n \neq 0$ we get  
	\[
	\dec(s_{2n}) = b_n \dec(s_2^n) = b_n (\one\otimes s_2)^n = \one\otimes s_{2n}.
	\]
	Thus $\dec'(s_N) = 0$ and also $D_{<N}(s_N)=0$ for all $N$.
\end{Example}

\begin{Proposition} 
	\label{prop:kernel_U_extended}
	Let $\CUf$ be as in Definition \ref{def:dec_Uf_general} 
	and let  $D_w$ be extended as above, then  
	we have for each $N\geq 2$ 
	\begin{equation*}
		\ker(D_{<N}) \cap \CUf_{N} = \Q \,s_N.
	\end{equation*}
\end{Proposition}

\begin{proof}  
	We first show that $s_N$ is contained in $\ker(D_{<N})$. 
	If $N$ is even, then we have seen in Example \ref{ex:s_even_in_ker} that $s_N \in \ker(D_{<N})$. If $N$ is odd, then $\dec'(s_N) = s_N \otimes \one$ and hence $\pi_{w}(s_N) = 0$ for each $w<N$. Hence $D_{<N}(s_N) = 0$.
	
	It remains to show that
	\begin{equation*}
		\ker(D_{<N}) \cap \CUf_{N} \subseteq \Q\, s_N.
	\end{equation*}
	Let $\xi\in\ker(D_{<N})\cap \CUf_{N}$ and write 
	\begin{equation}
		\label{eq:generic_element_UN}
		\xi = \alpha\, s_N  + \sum_{3\leq 2r+1 < N} s_{2r+1} u_r
	\end{equation}
	with $\alpha\in\Q$ and $u_r\in \big(\CUdual\otimes \Q[s_2]\big)_{N-2r-1}$. 
	For $2r+1<N$, we have by definition
 
	\begin{equation*}
		\dec'(\xi) = s_{2r+1} \otimes u_r + \begin{pmatrix} \text{tensor products whose left factor is a  } \\
			\text{product of the $s_i$  or of weight $\neq 2r+1$}\end{pmatrix},
	\end{equation*}
	and hence
	\begin{align} \label{eq:proof_ker} D_{2r+1}(\xi)&=(\pi_{2r+1}\otimes\id)\circ\dec'(\xi)\nonumber\\ 
		&=\pi_{2r+1}(s_{2r+1})\otimes u_r\\ \nonumber
		&\quad +\begin{pmatrix} \text{tensor products whose left factor is a concatenation   } \nonumber \\
			\text{ of the $s_i$ of weight $ 2r+1$ modulo shuffle products  }\end{pmatrix}. 
	\end{align}
	By assumption, we have $\xi\in \ker(D_{<N})$, in particular $D_{2r+1}(\xi)=0$. 
	It is $\pi_{2r+1}(s_{2r+1})$  non zero and  linear independent to all possible left factors of the second summand. 
	So we deduce from \eqref{eq:proof_ker} that $s_{2r+1}\otimes u_r=0$ and hence $u_r=0$. Since this holds for all $2r+1$, we get by \eqref{eq:generic_element_UN} that $\xi = \alpha\, s_N$.
\end{proof}

\section{The Ihara bracket and the Goncharov coproduct}
 
We want to apply the general constructions from Section \ref{Algebraic background} to a particular setup, which arises in the context of formal and motivic multiple zeta values. We will describe the operations purely algebraically in this section and explain their occurrence in the context of formal multiple zeta values later.

\subsection{Post-Lie algebras and the Grossman-Larson product}

We explain the rather new and abstract theory of post-Lie algebras and their universal enveloping algebras in general. In this context, we also introduce the Grossman-Larson product. All results in this subsection are taken from \cite{elm}.
Later we use this as a convenient algebraic approach to the Ihara bracket and the Goncharov coproduct.

In the following, we fix a commutative $\QQ$-algebra $R$ with unit. 

\begin{Definition}   \label{def:post-Lie}
A \emph{post-Lie algebra} $(\g, [-, -], \tr)$ is a Lie algebra $(\g, [-, -])$ over $R$ together with a $R$-bilinear product $\tr: \g\times\g\to\g$ such that for all $x,y,z \in \g$
\begin{align}
x \tr [y,z] &= [x \tr y,z] + [y,x \tr z], \label{def:post-lie1} \\
[x,y] \tr z &= x \tr (y \tr z) - (x \tr y) \tr z  -  y \tr (x \tr z) + (y \tr x) \tr z   
.\label{def:post-lie2}
\end{align}
\end{Definition}

\begin{Remark}
With the notation $L_x(y)= x \tr y$ we get that \eqref{def:post-lie1} is equivalent to $L_x$ being a Lie derivation. Equivalent to  \eqref{def:post-lie2} is the identity
\[L_{[x,y] +L_x(y) +L_y(x) }  =   L_x \circ L_y - L_y \circ L_x .\]
\end{Remark}

It is an easy exercise to check the following proposition.

\begin{Proposition} \label{prop:post-lie}
Let $ (\g, [-, -], \tr)$ be a post-Lie algebra. The bracket
\begin{equation}
\label{def:post-lie3}
\{ x,y \} = x \tr y - y \tr x + [x,y]
\end{equation}
satisfies the Jacobi identity for all $x, y \in \g$. Therefore, $\bar{\g}=(\g, \{-, -\})$ is also a Lie algebra. \qed
\end{Proposition}

\begin{Definition} \label{def:ext_tr} Let $(\g, [\cdot, \cdot], \tr)$ be a post-Lie algebra, and $\mathcal{U}(\g)$ the universal enveloping algebra of $(\g, [-,-])$. Let
\[\tr:\mathcal{U}(\g)\times \mathcal{U}(\g)\to\mathcal{U}(\g)\]
be the extension of the product $\tr:\g\times \g\to \g$ recursively given by
\begin{align}
&  x \tr 1 = 0 	\label{Ext0}\\
& 1 \tr A  = A						\label{Ext1}\\
& xA \tr y 	= x \tr (A \tr y)	 - (x \tr A) \tr y 	\label{Ext2}\\
& A  \tr BC	= (A_{(1)} \tr B)(A_{(2)} \tr C).	        \label{Ext3}	 
\end{align}
for all $A,B,C \in \mathcal{U}(\g)$ and $x,y \in \g$. Here, we use Sweedler's notation for the coproduct on $\mathcal{U}(\g)$:
\[\Delta(A)=A_{(1)}\otimes A_{(2)}.\]
\end{Definition}

A simple application of \eqref{Ext1} and \eqref{Ext3} leads to the following.

\begin{Lemma} \label{lem:Ext_derivation-formula} Let $(\g,[-,-],\tr)$ be a post-Lie algebra. We have for $x,t_1,\ldots,t_n\in\g$
\begin{equation*}
x \tr (t_1 \cdots t_n) = \sum_{i=1}^n t_1 \cdots t_{i-1} (x\tr t_i )t_{i+1}  \cdots t_n.
\end{equation*}
\end{Lemma}

\begin{Proposition} (\cite[Proposition 3.1]{elm})  \label{prop:post-lieUg} Let $(\g,[-,-],\tr)$ be a post-Lie algebra. Then the product $\tr:\mathcal{U}(\g)\times \mathcal{U}(\g)\to\mathcal{U}(\g)$ given in Definition \ref{def:ext_tr} is well-defined and unique. \qed
\end{Proposition} 

To prove this, use Lemma \ref{lem:Ext_derivation-formula} together with
\eqref{Ext3} and induction on the length of $B$ to extend $A \tr B$ to monomials $A$ and $B$. For more details we refer to \cite{elm} and in particular to \cite[Proposition 2.7]{og}. 

\begin{Definition}\label{def:grossman_larson}
Let $(\g,[-,-],\tr)$ be a post-Lie algebra. On the universal enveloping algebra $\mathcal{U}(\g)$ the \emph{Grossman-Larson product} is defined by
\begin{equation*}
A \glp B = A_{(1)}(A_{(2)} \tr B).
\end{equation*}
\end{Definition}

If $x\in \g$, then for all $A \in \mathcal{U}(\g)$ we get
$x \glp A = x \tr A + xA$. In particular, for $y \in \g$ we recover the Lie bracket of  Proposition \ref{prop:post-lie} via
\[x\glp y -y\glp x = x \tr y - y \tr x + [x,y] = \{x,y\}.\]
For any grouplike element $G \in \operatorname{Grp}(\widehat{\mathcal{U}}(\g))$ we get for any $A \in \mathcal{U}(\g)$ the formula 
\[G\glp A=G(G\tr A).\] 

Central for the application we have in mind is the following theorem.

\begin{Theorem} (\cite[Theorem 3.4]{elm}) \label{thm:new-universal} The Grossman-Larson product is associative and defines on $\mathcal{U}(\g)$  also the structure of
an  associative Hopf algebra $(\mathcal{U}(\g),\glp, \Delta)$.
Moreover, this Hopf algebra is isomorphic to the enveloping algebra $\mathcal{U}(\bar{\g})$. \qed
\end{Theorem}

For the proof, first observe that the Grossman-Larson product $\glp$ preserves the filtration on $\mathcal{U}(\g)$ given by the length of the monomials. The associated graded to the Grossman-Larson product with respect to this filtration is simply concatenation. So we obtain an isomorphism of graded Hopf algebras $\gr \mathcal{U}(\bar{\g})\to \mathcal{U}(\g)$, which results in a Hopf algebra isomorphism $\mathcal{U}(\bar{\g}) \to (\mathcal{U}(\g),\glp,\Delta)$. For a detailed proof in the analogue commutative setup we refer to \cite[Theorem 2.12]{og}.

The antipode $S$ of the Hopf algebra $(\mathcal{U}(\g),\glp,\Delta)$ differs from the standard antipode on the universal enveloping algebra $\mathcal{U}(\g)$. It can be computed recursively from
\begin{align} \label{eq:antipode_glp}
\glp\circ(S\otimes\id)\circ \Delta=\eta\circ\varepsilon,
\end{align}
where $\eta$ is the usual unit map und $\varepsilon$ the usual counit map of $\mathcal{U}(\g)$.

\begin{Example} \label{ex:antipode_glp} For $x,y\in\g$, we obtain
\begin{align*}
\glp\circ(S\otimes\id)\circ \Delta(xy)&=S(xy)\glp 1+S(x)\glp y+S(y) \glp x+S(1)\glp xy \\
&= S(xy)-x\glp y- y\glp x+ xy \\
&= S(xy)-x\tr y-xy-y\tr x-yx+xy \\
&=S(xy)-x \tr y-y\tr x-yx.
\end{align*}
Since we have $\eta(\varepsilon(xy))=0$, we deduce from \eqref{eq:antipode_glp}
\begin{align*}
S(xy)= x\tr y+y\tr x+yx.   
\end{align*}
The standard antipode on $\mathcal{U}(\g)$ is just given by $xy\mapsto yx$.

Similar computations show that we have for $x,y,z\in \g$
\begin{align*}
S(xyz)=\ & - zyx - x \tr zy - y \tr zx - z \tr yx + [x,y] \tr z + [x,z] \tr y + [y,z] \tr x \\
& - x \tr (z \tr y) - x \tr (y \tr z) - (y \tr x) \tr z - y \tr (z \tr x)  - (z \tr x) \tr y \\
& - (z \tr y) \tr x
\end{align*}
and the standard antipode on $\mathcal{U}(\g)$ is just $xyz\mapsto-zyx$.
\end{Example}

These examples show that a closed formula for the antipode $S$ is not at all obvious.

\subsection{The Ihara bracket}
  
We introduce the Ihara bracket and explain it in the context of post-Lie algebras. In the following, let $\fX$ be the free Lie algebra over $\QQ$ on the alphabet $\X=\{x_0,x_1\}$, and we denote its Lie bracket by $[-,-]$.

\begin{Definition}
For each $f \in \fX$ a special derivation is defined  by 
\[
d_f(x_i)=\begin{cases} 0 & \mbox{ if }  i=0, \\
[x_1,f]  &  \mbox{ if }  i=1. 
\end{cases}
\]
Let $\tr:\fX\times \fX\to \fX$ be the bilinear product given by
\[f\tr g=d_f(g).\]
\end{Definition}

\begin{Proposition}
The tuple $(\fX,[-,-],\tr)$ is a post-Lie algebra.
\end{Proposition}

\begin{proof}
By construction \eqref{def:post-lie1} is satisfied. The equation in \eqref{def:post-lie2} is equivalent to
\[d_{d_f(g)-d_g(f)+[f,g]}=d_f\circ d_g-d_g\circ d_f, \qquad f,g\in \fX.\]
On both sides are derivations, hence it suffices to check the equation on the generators $x_0,x_1$. For the first case, we directly obtain
\[d_{d_f(g)-d_g(f)+[f,g]}(x_0)=0=(d_f\circ d_g-d_g\circ d_f)(x_0).\]
For the second case, we compute
\begin{align*}
(\d_f\circ d_g-d_g\circ d_f)(x_1)&=d_f([x_1,g])-d_g([x_1,f]) \\
&=[x_1,d_f(g)]+[d_f(x_1),g]-[x_1,d_g(f)]-[d_g(x_1),f] \\
&=[x_1,d_f(g)]-[x_1,d_g(f)]+[[x_1,f],g]-[[x_1,g],f] \\
&=[x_1,d_f(g)]-[x_1,d_g(f)]+[x_1,[f,g]] \\
&=d_{d_f(g)-d_g(f)+[f,g]}(x_1),
\end{align*}
where the fourth step follows from the Jacobi identity for $[-,-]$.
\end{proof}

\begin{Definition} \label{def:Ihara}
The \emph{Ihara bracket} on $\fX$ is defined as
\[\{f,g\}=d_f(g)-d_g(f)+[f,g].\]
\end{Definition}

An immediate consequence of Proposition \ref{prop:post-lie} is the following.

\begin{Corollary}
The pair $(\fX,\{-,-\})$ is a Lie algebra.
\end{Corollary} 

\subsection{The Grossman-Larson product for the Ihara bracket}

We now study on the universal enveloping algebra \[\mathcal{U}(\fX)=(\QX,\conc,\co)\] 
the Grossman-Larson product (Definition \ref{def:grossman_larson})  determined by the Ihara bracket. We start with some formulas, which simplify the calculation of the extended product $\tr$.

\begin{Proposition}\label{prop:post_lie_on_letter}
Let $n\in \NN$. For all $a_1,\ldots,a_n\in \fX$, we have
\begin{align*}
(a_1\cdots a_n)\tr x_i=\begin{cases}
0 & \quad \text{if} \quad i=0, \\ 
[\cdots[[x_1,a_1],a_2],\ldots,a_n] & \quad \text{if} \quad i=1.
\end{cases}
\end{align*}
\end{Proposition}

\begin{proof}
For $n=1$, the claim follows from the definition of $\tr$. We prove the claim by induction on the number $n$ of Lie elements. By \eqref{Ext2} and Lemma \ref{lem:Ext_derivation-formula}, we have
\begin{align*} 
(a_1\cdots a_n)\tr x_0 &= a_1\tr \big((a_2\cdots a_n)\tr x_0\big)-\big(a_1\tr (a_2\cdots a_n)\big)\tr x_0 \\
&=a_1\tr 0 -\big(a_1\tr (a_2\cdots a_n)\big)\tr x_0 \\
&=-\sum_{i=2}^n \Big(a_2\cdots a_{i-1} (a_1 \tr a_i) a_{i+1}\cdots a_n\Big)\tr x_0 \\
&=0.
\end{align*}
The last step follows from the induction hypotheses, since $a_1 \tr a_i$ is also a Lie element. Similarly, if $a_2 \neq x_1$  we compute
\begin{align*} 
(a_1\cdots a_n)\tr x_1 &= a_1\tr \big((a_2\cdots a_n)\tr x_1\big)-\big(a_1\tr (a_2\cdots a_n)\big)\tr x_1 \\
&=a_1\tr [\cdots[[x_1,a_2],a_3],\ldots,a_n] -\big(a_1\tr (a_2\cdots a_n)\big)\tr x_1 \\
&=[\cdots[[x_1,a_1],a_2],\ldots,a_n]\\
&+\sum_{i=2}^n \Big([\cdots[[x_1,a_2],a_3],\ldots,a_1 \tr a_i],\ldots,a_n]\\
&\hspace{6cm}-\big(a_2a_3\cdots (a_1 \tr a_i)\cdots a_n\big)\tr x_1\Big) \\
&=[\cdots[[x_1,a_1],a_2],\ldots,a_n] \\
&+\sum_{i=2}^n \Big([\cdots[[x_1,a_2],a_3],\ldots,a_1 \tr a_i],\ldots,a_n] \\
&\hspace{5cm}-[\cdots[[x_1,a_2],a_3],\ldots,a_1 \tr a_i],\ldots,a_n]\Big) \\
&=[\cdots[[x_1,a_1],a_2],\ldots,a_n].
\end{align*}
For $a_2=x_1$ the claim follows by a small modification of the above calculations.
\end{proof}

\begin{Proposition}\label{prop:post_lie_factor_out_x0}
For $A,B\in \mathcal{U}(\fX)$, we have
\[
A\tr x_0 B = x_0(A \tr B).
\]
\end{Proposition}

\begin{proof}
By \eqref{Ext3}, we have
\[
A\tr x_0 B = (A_{(1)} \tr x_0) (A_{(2)} \tr B).
\]
Since $A_{(1)}$ consists of (possibly empty) products of elements in $\fX$, we deduce from Proposition \ref{prop:post_lie_on_letter} that $(A_{(1)} \tr x_0)$ vanishes unless $A_{(1)}$ corresponds to the empty word (cf. \eqref{Ext1}). So the summand $\one\otimes A$ is the only one in $\Delta(A)$ with a non-zero contribution. Hence 
\[
(A_{(1)} \tr x_0) (A_{(2)} \tr B) = (\one\tr x_0) (A\tr B) = x_0(A\tr B). \qedhere
\]
\end{proof}

\begin{Proposition}
For $A,B\in \mathcal{U}(\fX)$, we have
\begin{equation*}
x_1 A \tr B = 0.
\end{equation*}
\end{Proposition}

\begin{proof}
Write $B = b_1\cdots b_m$ where $b_i\in \X$.
We prove the claim via induction on $\dep(B)=\#\{i \mid b_i = x_1\}$. 
If $B = x_0^m$ then repeatedly applying Proposition \ref{prop:post_lie_factor_out_x0} yields 
\begin{equation*}
x_1 A \tr B = x_0^{m-1} (x_1 A\tr x_0) = 0,
\end{equation*}
where the last equality follows from Proposition \ref{prop:post_lie_on_letter}.
Now assume $\dep(B)=M\geq 1$, and let $j\in\{1,\dots,m\}$ denote the smallest integer such that $b_j = x_1$. By applying Proposition \ref{prop:post_lie_factor_out_x0} $j-1$ times we obtain that
\begin{equation*}
x_1A \tr B = x_0^{j-1} (x_1A\tr x_1 b_{j+1}\cdots b_m) 
\overset{\eqref{Ext3}}{=} x_0^{j-1} \big( (x_1A)_{(1)} \tr x_1\big) \big( (x_1A)_{(2)} \tr b_{j+1}\cdots b_m\big).
\end{equation*}
Since $\Delta(x_1 A) = (x_1\otimes \one + \one \otimes x_1)\Delta(A)$ it is clear that each tensor product in $\Delta(x_1 A)$ has at least one factor starting in $x_1$. If the left factor starts in $x_1$, then $\big( (x_1A)_{(1)} \tr x_1\big) = 0$ by Proposition \ref{prop:post_lie_on_letter} since $[x_1,x_1] = 0$. If the right factor starts in $x_1$ then $\big( (x_1A)_{(2)} \tr b_{j+1}\cdots b_m\big) = 0$ by induction hypothesis since $\dep(b_{j+1}\cdots b_m)=M-1$.
\end{proof}

\begin{Example} We calculate $x_0x_0\glp x_0x_1$. At first we observe
\begin{align*}
\co(x_0x_0)=x_0x_0\otimes\one+2x_0\otimes x_0+\one\otimes x_0x_0,
\end{align*}
and thus by Definition \ref{def:grossman_larson}  we need to determine 
\begin{align*}
x_0x_0 \glp x_0x_1 &= x_0x_0(\one\tr x_0x_1)+2x_0(x_0\tr x_0x_1)+\one(x_0x_0\tr x_0x_1) \\
&=x_0x_0x_0x_1+2x_0(x_0\tr x_0x_1)+x_0x_0\tr x_0x_1.
\end{align*}
For the third term, we obtain from Propositions \ref{prop:post_lie_factor_out_x0} and \ref{prop:post_lie_on_letter}
\begin{align*}
x_0x_0\tr x_0x_1=x_0(x_0x_0\tr x_1)=x_0[[x_1,x_0],x_0].
\end{align*}
Similarly, one computes for the second term
\begin{align*}
x_0\tr x_0x_1=x_0(x_0\tr x_1)=x_0[x_1,x_0].
\end{align*}
Altogether, we get
\begin{align*}
x_0x_0 \glp x_0x_1
&=x_0x_0x_0x_1+2x_0x_0[x_1,x_0]+x_0[[x_1,x_0],x_0] \\
&=x_0x_1x_0x_0.
\end{align*}
\end{Example}

We want to give a closed formula for the Grossman-Larson product corresponding to the Ihara bracket.

Recall the $j$-th letter function from Notation \ref{not:alphabets}, i.\,e.
\begin{align*}
	\ve_j:\X^*&\to \X\cup \{\one\}, \\
	x_{i_1}\cdots x_{i_n}&\mapsto\begin{cases}
		x_{i_j} \quad & \text{ if } j\leq n, \\
		\one & \text{ else}.
	\end{cases}
\end{align*}
Instead of $w=\ve_1(w)\ve_2(w)\cdots \ve_n(w)\in \X^*$, we will often just write $w=\ve_1\cdots \ve_n$ where $n$ is the \emph{weight}, i.e., the number of letters, of $w$.
The antipode $S$ of the Hopf algebra $(\mathcal{U}(\fX),\conc,\co)$ is given by
\begin{align} \label{eq:antipode_QX}
S: \mathcal{U}(\fX)&\to \mathcal{U}(\fX), \\
\ve_1\cdots \ve_n&\mapsto (-1)^n \ve_n\cdots \ve_1. \nonumber
\end{align}

\begin{Proposition} \label{thm:GL_for_Ihara} Let $A\in \mathcal{U}(\fX)$ and all words $w=x_0^{k_1}x_1\cdots x_0^{k_d}x_1x_0^{k_{d+1}}\in \mathcal{U}(\fX)$, we have
\begin{align*}
A\glp w= A_{(1)}x_0^{k_1}S(A_{(2)})x_1 A_{(3)}x_0^{k_2}\cdots x_0^{k_d}S(A_{(2d)})x_1 A_{(2d+1)}x_0^{k_{d+1}}.
\end{align*}
Here, we use the iterated Sweedler notation
\[\co^{n}(A)=A_{(1)}\otimes \cdots \otimes A_{(n+1)}, \qquad n\geq1.\]
\end{Proposition}

\begin{proof} For any $A\in \mathcal{U}(\fX)$ and any word $w=x_0^{k_1}x_1\cdots x_0^{k_d}x_1x_0^{k_{d+1}}\in \mathcal{U}(\fX)$, we obtain from Definition \ref{def:grossman_larson} and several applications of \eqref{Ext3} that
\begin{align*}
A\glp w&=A_{(1)}(A_{(2)}\tr w)=A_{(1)}(A_{(2)}\tr x_0^{k_1}x_1)\cdots (A_{(d+1)}\tr x_0^{k_d}x_1)(A_{(d+2)}\tr x_0^{k_{d+1}})
\end{align*}
From Proposition \ref{prop:post_lie_factor_out_x0}, we deduce that
\begin{align} \label{eq:proof_glp_ihara1}
A\glp w&=A_{(1)}x_0^{k_1}(A_{(2)}\tr x_1)\cdots x_0^{k_d}(A_{(d+1)}\tr x_1)x_0^{k_{d+1}}.
\end{align}
By Proposition \ref{prop:post_lie_on_letter} and the definition of the coproduct $\co$ (Example \ref{duality shuffle product}) and the antipode $S$ (see \eqref{eq:antipode_QX}), we have for $A=a_1\cdots a_n\in \mathcal{U}(\fX)$ 
\begin{align} \label{eq:proof_glp_ihara2}
A\tr x_1=[\cdots[[x_1,a_1],a_2],\ldots,a_n]=S(A_{(1)})x_1A_{(2)}.
\end{align}
Combining \eqref{eq:proof_glp_ihara1} and \eqref{eq:proof_glp_ihara2} gives the claimed formula
\begin{align*}
A\glp w= A_{(1)}x_0^{k_1}S(A_{(2)})x_1 A_{(3)}x_0^{k_2}\cdots x_0^{k_d}S(A_{(2d)})x_1 A_{(2d+1)}x_0^{k_{d+1}}.
\end{align*}
\end{proof}

From Theorem \ref{thm:new-universal}, we obtain the following.

\begin{Theorem} \label{cor:U_hopf}
    The tuple $(\mathcal{U}(\fX),\glp,\co)$ is a cocommutative Hopf algebra.
\end{Theorem}

This result was also stated in by Willwacher \cite[Proposition 7.1]{willwa}, but proved completely differently.

\subsection{The Goncharov coproduct}

The aim of this subsection is to complete the diagram \eqref{diagram big picture general} for the Ihara bracket. In the last subsection, we studied the universal enveloping algebra $(\mathcal{U}(\fX),\circledast,\co)$ of $(\fX,\{-,-\})$. So the next step is to consider the grouplike elements of $(\mathcal{U}(\fX),\circledast,\co)$. In order to get a non-empty set, we have to pass to the completed Hopf algebra $(\widehat{\mathcal{U}}(\fX),\circledast,\co)$ in the sense of Definition \ref{def:completed Hofp algebra}. So precisely, consider the set
\[\operatorname{Grp}(\widehat{\mathcal{U}}(\fX))=\{A\in \kXc{\QQ}\mid \co(A)=A\otimes A\}.\]
From Corollary \ref{cor:U_hopf} and Theorem \ref{Thm grouplike primitive} we deduce the following.

\begin{Corollary} \label{cor:Grp_glp}
The pair $(\operatorname{Grp}(\widehat{\mathcal{U}}(\fX)),\glp)$ is a group.
\end{Corollary}

\begin{Remark} \label{rem:GLP_ihara_grouplike}
By Theorem \ref{Thm grouplike primitive}, the grouplike elements $G\in\Grp(\widehat{\mathcal{U}}(\fX))$ satisfy \[S(G)=G^{-1},\] 
where the inverse is meant with respect to concatenation\footnote{In fact, the triple $(\Grp(\widehat{\mathcal{U}}(\fX), \conc, \glp)$ is a post-group. This is a general phenomenon for any post-Lie algebra: the corresponding set of group like elements in the completed enveloping algebra is endowed with two group structures: one for the "ordinary" product given by concatenation, the other for the Grossman Larson product. These product structures interact nicely. For an introduction to post-groups, see
\cite{Bai_2023}
and \cite{alkaabi2024free} }. 
Hence we have for all $G,H \in \Grp(\widehat{\mathcal{U}}(\fX))$
\[ G\glp H = G \kappa_G(H), \]
where $\kappa_G$ is the automorphism on $\widehat{\mathcal{U}}(\fX)$ with respect to concatenation given by 
\[
\kappa_G(x_i)=\begin{cases} x_0 & \mbox{ if  }\,  i=0, \\
 G^{-1} x_1 G  &  \mbox{ if  }\,  i=1. 
\end{cases}
\]
This formula equals the one given in \cite{ra}.
\end{Remark} 

\begin{Remark}
Viewing $(\fX,\{-,-\})$ as a post-Lie algebra an applying the construction of the Grossman-Larson product gives a very direct and explicit way to derive the group multiplication on $\Grp(\widehat{\mathcal{U}}(\fX))$ from the Lie bracket $\{-,-\}$ on $\fX$. So the Grossman-Larson product can be seen as the reverse construction of the usual linearization of the group multiplication to obtain the Lie bracket (cf \cite[Section II.2.2]{ra}). 
\end{Remark}

For the last step of completing the diagram \eqref{diagram big picture general}, we view the grouplike elements as an affine group scheme
\begin{align} \label{eq: ags from grouplike}
\QAlg &\to \Groups, \\
R&\mapsto (\Grp(\widehat{\mathcal{U}}(\fX\otimes R)),\circledast). \nonumber
\end{align}
Then following Example \ref{grou-like elements ags}, this affine scheme is represented by the algebra $(\QX,\shuffle)$. We want to describe the coproduct on $(\QX,\shuffle)$ induced by the multiplication $\circledast$ on the affine group scheme (cf Theorem \ref{affine group schemes Hopf algebras}).  We follow the calculations given in \cite[Subsection 3.10.6.]{bgf} in order to obtain the explicit formula for this coproduct.  

\begin{Notation} \label{not:I}
For $\ve_a,\ve_b\in \X$ and any $f\in \kXc{\QQ}$, we set
\begin{equation*}
\Iformal(\ve_a;f;\ve_b) = 
\begin{cases}
f, & (\ve_a,\ve_b) = (x_1,x_0), \\
S(f), & (\ve_a,\ve_b) = (x_0,x_1), \\
(f\mid \one)\cdot \one, & \ve_a = \ve_b.
\end{cases}
\end{equation*}
\end{Notation}

\begin{Definition} \label{def:gon} (\cite{gon})  
Let $\ve_1\cdots \ve_n \in \X^*$, then the \emph{Goncharov coproduct} is given by the  formula 
\begin{align} \label{eq:Gon_formula}
\DeltaGon(\ve_1\cdots \ve_n) =
\SumGon \Iformal(\ve_{i_p}; \ve_{i_p+1}\cdots\ve_{i_{p+1}-1};\ve_{i_{p+1}}) \otimes \ve_{i_1}\cdots\ve_{i_k},
\end{align}
where we set $i_0=0,\, i_{k+1}=n+1$   and $\ve_{0}=x_1,\, \ve_{n+1}=x_0$.
\end{Definition}

For any $f\in\kXc{\QQ}$ with constant term $0$ we have that $\Iformal(\ve;f;\ve)=0$ for any $\ve\in \X$ (Notation \ref{not:I}). Hence, there are usually a lot of vanishing terms in the Goncharov coproduct. We give examples and  graphical interpretations for this coproduct in Subsection \ref{subsec:visual}. 

The following lemma is a reformulation of \cite[Proposition 3.422]{bgf}. 

\begin{Lemma} \label{lem:glp_gon_dual} 
For all $G,H\in\Grp(\widehat{\mathcal{U}}(\fX))$, and $w\in \QX$, we have  the identity
\[(G\glp H\mid w)=(G\otimes H\mid \DeltaGon(w)),\]
where the duality pairing  $(\cdot\mid \cdot)$ is given in \eqref{pairing quasi-shuffle}.
\end{Lemma}
The equality holds also for an arbitrary element $H\in \widehat{\mathcal{U}}(\fX)$, $H$ being grouplike is not used in the proof. But in everything what follows, we only need the result as stated above.

\begin{proof}
We follow the calculations given in \cite[Subsection 3.10.6.]{bgf}. For $G\in\operatorname{Grp}(\widehat{\mathcal{U}}(\fX))$, we have by Theorem \ref{Thm grouplike primitive}
\begin{equation*}
\Iformal(\ve_0;G;\ve_{n+1}) = 
\begin{cases}
G, & (\ve_0,\ve_{n+1}) = (x_1,x_0), \\
G^{-1}, & (\ve_0,\ve_{n+1}) = (x_0,x_1), \\
\one, & \ve_0 = \ve_{n+1}.
\end{cases}
\end{equation*}
Let $G\in\operatorname{Grp}(\widehat{\mathcal{U}}(\fX))$ and $w\in \QX$. By Theorem \ref{thm:GL_for_Ihara}, the product $G\glp w$ can be computed as follows:
\begin{itemize}
\item if $w$ starts with $x_0$, put $G$ at the beginning,
\item between every $x_0$ and $x_1$ in $w$, insert $G^{-1}$,
\item between every $x_1$ and $x_0$ in $w$, insert $G$,
\item if $w$ ends with $x_1$, put $G$ at the end.
\end{itemize}
Hence, we obtain
\begin{align*}
G\glp(\ve_1\cdots \ve_n)=\Iformal(x_1;G;\ve_1)\ve_1\Iformal(\ve_1;G;\ve_2)\ve_2\ldots \ve_{n-1}\Iformal(\ve_{n-1};G;\ve_n)\ve_n I(\ve_n;G;x_0).
\end{align*}
For $G,H\in\operatorname{Grp}(\widehat{\mathcal{U}}(\fX))$ with
\[H=\sum_{w\in \X^*} h(w) \ve_1(w)\cdots \ve_{\wt(w)}(w),\qquad \ve_i(w)\in \X,\]
and $\xi_1\cdots \xi_n\in \X^*$, we compute 
\begin{align*}
(G\glp H&\mid \xi_1\cdots \xi_n) \\
&=\sum_{w\in \X^*}h(w)\Big(\Iformal\big(x_1;G;\ve_1(w)\big)\ve_1(w)\Iformal\big(\ve_1(w);G;\ve_2(w)\big)\\
&\hspace{6cm}\cdots \ve_{\wt(w)}(w)\Iformal\big(\ve_{\wt(w)}(w);G;x_0\big)\ \Big|\xi_1\cdots \xi_n \Big) \\
&=\Sum \Prod \Big(\Iformal\big(\xi_{i_p};G;\xi_{i_{p+1}}\big)\ \Big|\ \xi_{i_p+1}\cdots\xi_{i_{p+1}-1}\Big)\Big(H\ \Big|\ \xi_{i_1}\cdots \xi_{i_k}\Big). \\
\intertext{Since $\Iformal(\xi_a;-;\xi_b)$ is either an involution or constant, we deduce} 
(G\glp H&\mid \xi_1\cdots \xi_n) \\
&=\Sum \Prod \Big(G\ \Big|\ \Iformal\big(
\xi_{i_p};\xi_{i_p+1},\ldots,\xi_{i_{p+1}-1};\xi_{i_{p+1}}\big)\Big)\Big(H\ \Big|\ \xi_{i_1}\cdots \xi_{i_k}\Big). \\
\intertext{Any grouplike element $G$ for $\co$ satisfies by duality \[(G\mid w_1)(G\mid w_2)=(G\mid w_1\ \shuffle\ w_2)\] for all $w_1,w_2\in \QX$. Hence, we get}
(G\glp H&\mid \xi_1\cdots \xi_n) \\
&=\Sum \Big(G\ \Big|\ \Bigshuffle_{p=0}^k \Iformal\big(
\xi_{i_p};\xi_{i_p+1},\ldots,\xi_{i_{p+1}-1};\xi_{i_{p+1}}\big)\Big)\Big(H\ \Big|\ \xi_{i_1}\cdots \xi_{i_k}\Big) \\
&=\Sum \Big(G\otimes H\ \Big|\ \Bigshuffle_{p=0}^k \Iformal\big(\xi_{i_p};\xi_{i_p+1},\ldots,\xi_{i_{p+1}-1};\xi_{i_{p+1}}\big)\otimes \xi_{i_1}\cdots \xi_{i_k}\Big)\\
&=\big(G\otimes H\ \big|\ \DeltaGon(\xi_1\cdots \xi_n)\big),
\end{align*}
where $\xi_{i_0}=x_1$ and $\xi_{i_{k+1}}=x_0$.
\end{proof} 

Note that for any two words $u,v \in \QX$ we can find elements $G ,H \in\Grp(\widehat{\mathcal{U}}(\fX))$ such that $(G\otimes H\,| \, u \otimes v ) \neq 0$. Thus, Lemma \ref{lem:glp_gon_dual} uniquely determines the coproduct on $\QX$ dual to $\glp$ with respect to the pairing $(\cdot |\cdot)$.

Since $R\mapsto (\Grp(\widehat{\mathcal{U}}(\fX\otimes R)),\circledast)$ is an affine scheme represented by the algebra $(\QX,\shuffle)$, and by Lemma \ref{lem:glp_gon_dual} the multiplication $\circledast$ corresponds to the Goncharov coproduct on $(\QX,\shuffle)$, we get the following.

\begin{Proposition}\label{prop:QX_gon_hopf-alg}
The tuple $(\QX,\shuffle,\DeltaGon)$ is a weight-graded Hopf algebra.
\end{Proposition}

Summarizing the previous results gives the following diagram (cf \eqref{diagram big picture general})
\begin{equation}  \label{diagram big picture grouplike}
\begin{tikzcd}[baseline=(current  bounding  box.center), scale cd=0.965]
(\QX,\shuffle,\DeltaGon) \arrow[dddd,"\text{mod products}"'] &&& (\mathcal{U}(\fX),\glp,\co) \arrow[lll,"\sim","\text{graded dual}" '] 
\\ \\ 
& (\operatorname{Grp}(\widehat{\mathcal{U}}(\fX)),\glp) \arrow[uul,leftrightarrow,"1:1"] &&
\\ \\ 
(\Indec(\QX,\shuffle),\delta) &&& (\fX,\{-,-\})
\arrow[lll, "\sim", "\text{graded dual}"'] \arrow[uuuu, hookrightarrow] \arrow[uull,leftrightarrow,"1:1"]
\end{tikzcd} \end{equation} 

The bijective exponential map $\exp:\widehat{\fX}\to \operatorname{Grp}(\widehat{\mathcal{U}}(\fX))$ is given in Example \ref{ex ags}. By Theorem \ref{thm:new-universal}, this exponential map can also be used in this case.

\subsection{Combinatorics of the Goncharov coproduct} \label{subsec:visual}

Recall from Definition \ref{def:gon} that we have
\begin{align*}
\DeltaGon(\ve_1\cdots \ve_n) =
\SumGon \Iformal(\ve_{i_p}; \ve_{i_p+1}\cdots\ve_{i_{p+1}-1};\ve_{i_{p+1}}) \otimes \ve_{i_1}\cdots\ve_{i_k},
\end{align*}
where $\ve_1\dots\ve_n\in \X^*$, $i_0=0,\, i_{k+1}=n+1$ and $\ve_{0}=x_1,\, \ve_{n+1}=x_0$. 

We start by expressing the Goncharov coproduct in a more combinatorial way. 
\begin{Notation} \label{not:Gon_coprod} Let $w=\ve_1 \cdots \ve_n\in \X^*$ be a word of weight $n$. Then any subset $I=\{i_1,...,i_k\}$ of $I_n=\{1,...,n\}$ with $i_1<\cdots<i_k$ defines a word
\begin{align} \label{def:subword}
	v_I = \ve_{i_1} \cdots \ve_{i_k}.
\end{align}
So we express the Goncharov coproduct as
\begin{align} \label{eq:gon_comb}
	\DeltaGon(w) = \sum_{I \subset I_n } P_I(w) \otimes v_I, 
\end{align}
where
\begin{align} \label{eq:P_I}
	P_I(w)= \Bigshuffle_{p = 0}^k  \Iformal( \ve_{i_p}; u_p ;  \ve_{i_{p+1}} ).
\end{align}
Here we set $i_0=0$, $i_{k+1}=n+1$, $\ve_0=x_1$, $\ve_{n+1}=x_0$, and $u_p=\ve_{i_p+1}\cdots \ve_{i_{p+1}-1}$.
In the formula \eqref{eq:gon_comb} we refer to the words $v_I$ as \emph{subwords}
and to the $u_p$ in \eqref{eq:P_I} as \emph{quotient words}. Moreover, we call $\ve_{i_p}u_p\ve_{i_{p+1}}$ the \emph{enlarged subword} associated to $u_p$.
\end{Notation}

We explain now two ways to give an graphical interpretation of this formula.

On the one hand, there is the graphical interpretation introduced by Goncharov in \cite{gon}. Given a word $\ve_1\cdots\ve_n\in\X^*$, one locates the letters $\ve_1,\dots,\ve_n$ on the upper half of a semicircle and adds $x_1$ and $x_0$ at the left- and right-hand side, respectively. Then the term in the formula from Definition \ref{def:gon} that corresponds to the subword $v=\ve_{i_1}\cdots\ve_{i_k}$ is depicted by a polygon inscribed into the semicircle with vertices $x_1,\ve_{i_1},\dots,\ve_{i_k},x_0$.

For example the subword $\textcolor{blue}{\ve_2\ve_3\ve_4\ve_7\ve_n}$ corresponds to the summand
\begin{equation*}
\scalebox{0.9}{$
\textcolor{black}{\Iformal(\ve_0; \ve_1; \ve_2)
\Iformal(\ve_2; \one; \ve_3)
\Iformal(\ve_3; \one; \ve_4)
\Iformal(\ve_4;\ve_5\ve_6;\ve_7) \Iformal(\ve_7;\ve_8\dots\ve_{n-1};\ve_n)
\Iformal(\ve_n; \one; \ve_{n+1})
} \otimes \textcolor{blue}{\ve_2 \ve_3 \ve_4 \ve_7 \ve_n}
$}
\end{equation*}
and is depicted by
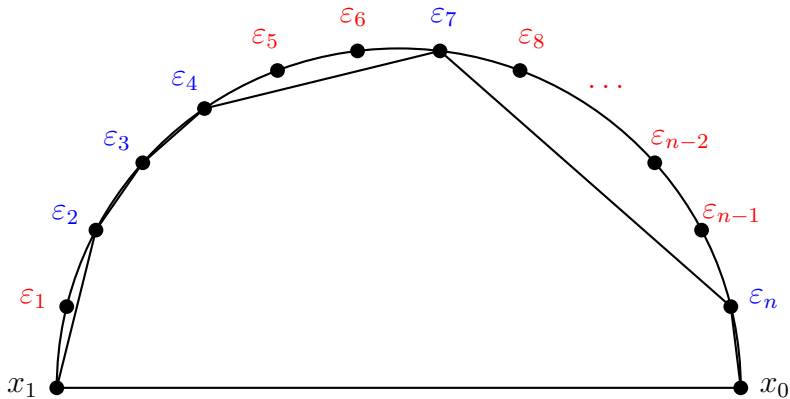
\begin{figure}
\begin{center}
\def\scale{4.5} 
\def\n{13} 
\begin{tikzpicture}[baseline = (current bounding box.north)]
	\coordinate (start) at (-\scale,0);
	\coordinate (end) at (\scale,0);
	
	\begin{scope}[nodes={fill=black, circle, scale=.5}]
		\clip (-\scale-.1,0) rectangle (\scale+.1,\scale+.1);
		\draw[thick] (start)+(end) circle (\scale);
	\end{scope}
		
	\draw[thick] (start) -- (end);
	
	\foreach \arc in {0,1,...,\n}{
		
		\coordinate (L\arc) at (180-\arc*180/\n : 1.1*\scale);
		
		\coordinate (N\arc) at (180-\arc*180/\n : \scale);
		
	}
	
	
	\draw[thick] (N0) to (N2)
	                  to (N3) 
	                  to (N4)
                      to (N7) 
                      to (N12)
                      to (N13);
	
	\node[fill=black, circle, scale=.5] at (N0) {};
	\node at (L0) {$x_1$};
	
	\foreach \num in {2,3,4,7}{
		\node[fill=black, circle, scale=.5] at (N\num) {};
		
		\node[blue] at (L\num) {$\ve_{\num}$};
	}
	
	\foreach \num in {1,5,6,8}{
	    \node[fill=black, circle, scale=.5] at (N\num) {};
		
		\node[red] at (L\num) {$\ve_{\num}$};
	}
	
	\node[red] at (180-9*180/\n : 1.08*\scale) {$\dots$};	
	
	\foreach \num in {10,...,\n}{
		\node[fill=black, circle, scale=.5] at (N\num) {};
	}

	\node[red] at (L10) {$\ve_{n-2}$};
	\node[red] at (L11) {$\ve_{n-1}$};
	\node[blue] at (L12) {$\ve_{n}$};
	\node at (L13) {$x_0$};
\end{tikzpicture}
\caption{Goncharov's semicircle}\label{fig:semicircle}
\end{center}
\end{figure}

On the other hand, the left hand factors in a summand of the Goncharov 
coproduct corresponds to a choice of distinct strict subwords $u_p$ of $\ve_1\cdots\ve_n$ (cf Notation \ref{not:alphabets}). For example, the choice of $\textcolor{red}{\ve_1}$, 
$\textcolor{red}{\ve_5\ve_6}$ and $\textcolor{red}{\ve_8\cdots \ve_{n-2}\ve_{n-1}}$ in the word $\ve_1\cdots\ve_n$  
corresponds to the above summand in its reduced  form, i.e., the trivial factors are omitted,

\begin{equation*}
\textcolor{red}{\Iformal(\ve_0; \ve_1; \ve_2) \Iformal(\ve_4;\ve_5\ve_6;\ve_7) \Iformal(\ve_7;\ve_8\dots\ve_{n-1};\ve_n)} \otimes \textcolor{black}{\ve_2 \ve_3 \ve_4 \ve_7 \ve_n}.
\end{equation*}
This can be visualized by \emph{eating worms} as follows:
\begin{figure}
\begin{center}
\includegraphics[width=.9\textwidth]{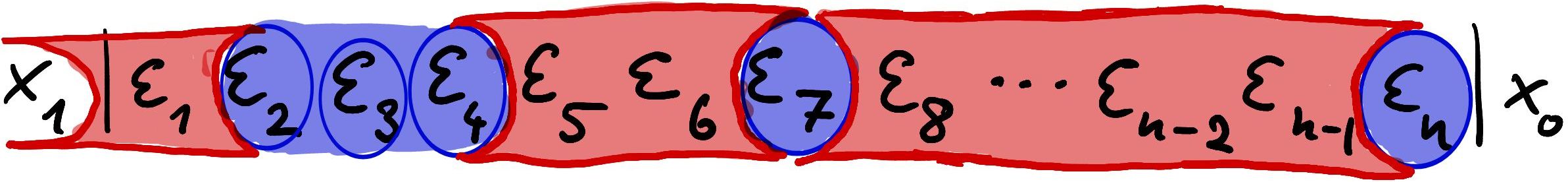}
\caption{''Eating worms''}\label{fig:eating_worms}
\end{center}
\end{figure}
We refer to the red and blue parts to ''red worms'' and ''blue worms'', respectively. The ''red worms'' have open mouths towards its boundaries which are either a ''blue worm'' or one of the outer boundaries $x_1$ and $x_0$.


\begin{Example}\label{exm:gon_coprod_visu}
We consider $x_0^{n-1}x_1\in\QX$ for $n\in\{1,2,3,4\}$.
For $n=1$ we compute directly that
\begin{align*}
\DeltaGon(x_1) = \Iformal(x_1;x_1;x_0) \otimes \one + \one \otimes x_1 = x_1 \otimes \one + \one \otimes x_1.
\end{align*}
The summands for $n=2$ are depicted via Goncharov's semicircle (cf. Figure \ref{fig:semicircle}) as follows
\begin{center}
\def\scale{1.12}
\def\n{3}
\begin{tikzpicture}[baseline = (current bounding box.north)]
	\coordinate (start) at (-\scale,0);
	\coordinate (end) at (\scale,0);
	
	\begin{scope}[nodes={fill=black, circle, scale=.5}]
		\clip (-\scale-.1,0) rectangle (\scale+.1,\scale+.1);
		\draw[thick] (start)+(end) circle (\scale);
	\end{scope}
	
	\draw[thick] (start) -- (end);
	
	\foreach \arc in {0,1,...,\n}{
		
		\coordinate (L\arc) at (180-\arc*180/\n : 1.3*\scale);
		
		\coordinate (N\arc) at (180-\arc*180/\n : \scale);
	}
	
	\draw[thick] (N0) to [bend left=10] (N3);
	
	\foreach \arc in {0,1,...,\n}{
	    \node[fill=black, circle, scale=.5] at (N\arc) {};
	}
	
	\node at (L0) {$x_1$};
	\node[red] at (L1) {$x_0$};
	\node[red] at (L2) {$x_1$};
	\node at (L3) {$x_0$};
\end{tikzpicture}
\begin{tikzpicture}[baseline = (current bounding box.north)]
	\coordinate (start) at (-\scale,0);
	\coordinate (end) at (\scale,0);
	
	\begin{scope}[nodes={fill=black, circle, scale=.5}]
		\clip (-\scale-.1,0) rectangle (\scale+.1,\scale+.1);
		\draw[thick] (start)+(end) circle (\scale);
	\end{scope}
	
	\draw[thick] (start) -- (end);
	
	\foreach \arc in {0,1,...,\n}{
		
		\coordinate (L\arc) at (180-\arc*180/\n : 1.3*\scale);
		
		\coordinate (N\arc) at (180-\arc*180/\n : \scale);
	}
	
	\draw[thick] (N0) to [bend right=0] (N1)
	                to [bend right=0] (N3);
	
	\foreach \arc in {0,1,...,\n}{
	    \node[fill=black, circle, scale=.5] at (N\arc) {};
	}
	
	\node at (L0) {$x_1$};
	\node[blue] at (L1) {$x_0$};
	\node[red] at (L2) {$x_1$};
	\node at (L3) {$x_0$};	
\end{tikzpicture}
\begin{tikzpicture}[baseline = (current bounding box.north)]
	\coordinate (start) at (-\scale,0);
	\coordinate (end) at (\scale,0);
	
	\begin{scope}[nodes={fill=black, circle, scale=.5}]
		\clip (-\scale-.1,0) rectangle (\scale+.1,\scale+.1);
		\draw[thick] (start)+(end) circle (\scale);
	\end{scope}
	
	\draw[thick] (start) -- (end);
	
	\foreach \arc in {0,1,...,\n}{
		
		\coordinate (L\arc) at (180-\arc*180/\n : 1.3*\scale);
		
		\coordinate (N\arc) at (180-\arc*180/\n : \scale);
	}
	
	\draw[thick] (N0) to [bend right=0] (N2)
	                to [bend right=0] (N3);
	
	\foreach \arc in {0,1,...,\n}{
	    \node[fill=black, circle, scale=.5] at (N\arc) {};
	}
	
	\node at (L0) {$x_1$};
	\node[red] at (L1) {$x_0$};
	\node[blue] at (L2) {$x_1$};
	\node at (L3) {$x_0$};
\end{tikzpicture}
\begin{tikzpicture}[baseline = (current bounding box.north)]
	\coordinate (start) at (-\scale,0);
	\coordinate (end) at (\scale,0);
	
	\begin{scope}[nodes={fill=black, circle, scale=.5}]
		\clip (-\scale-.1,0) rectangle (\scale+.1,\scale+.1);
		\draw[thick] (start)+(end) circle (\scale);
	\end{scope}
	
	\draw[thick] (start) -- (end);
	
	\foreach \arc in {0,1,...,\n}{
		
		\coordinate (L\arc) at (180-\arc*180/\n : 1.3*\scale);
		
		\coordinate (N\arc) at (180-\arc*180/\n : \scale);
	}
	
	\draw[thick] (N0) to [bend right=0] (N1)
                	to [bend right=0] (N2)
                	to [bend right=0] (N3);
	
	\foreach \arc in {0,1,...,\n}{
	    \node[fill=black, circle, scale=.5] at (N\arc) {};
	}
	
	\node at (L0) {$x_1$};
	\node[blue] at (L1) {$x_0$};
	\node[blue] at (L2) {$x_1$};
	\node at (L3) {$x_0$};
\end{tikzpicture}
\end{center}
Thus we get
\begin{align*}
\DeltaGon(x_0x_1)
&= \Iformal(x_1; x_0,x_1; x_0) \otimes \one + \Iformal(x_1; \one; x_0) \Iformal(x_0; x_1; x_0) \otimes x_0\\
&\hspace{0,4cm}+ \Iformal(x_1;x_0;x_1) I(x_1; \one ;x_0) \otimes x_1 + \Iformal(x_1;\one;x_0) \otimes x_0x_1\\
&= x_0x_1\otimes \one + \one \otimes x_0x_1.
\end{align*}
The summands for $n=3$ are visualized as ''eating worms'' as follows:
\begin{center}
\includegraphics[width=.8\textwidth]{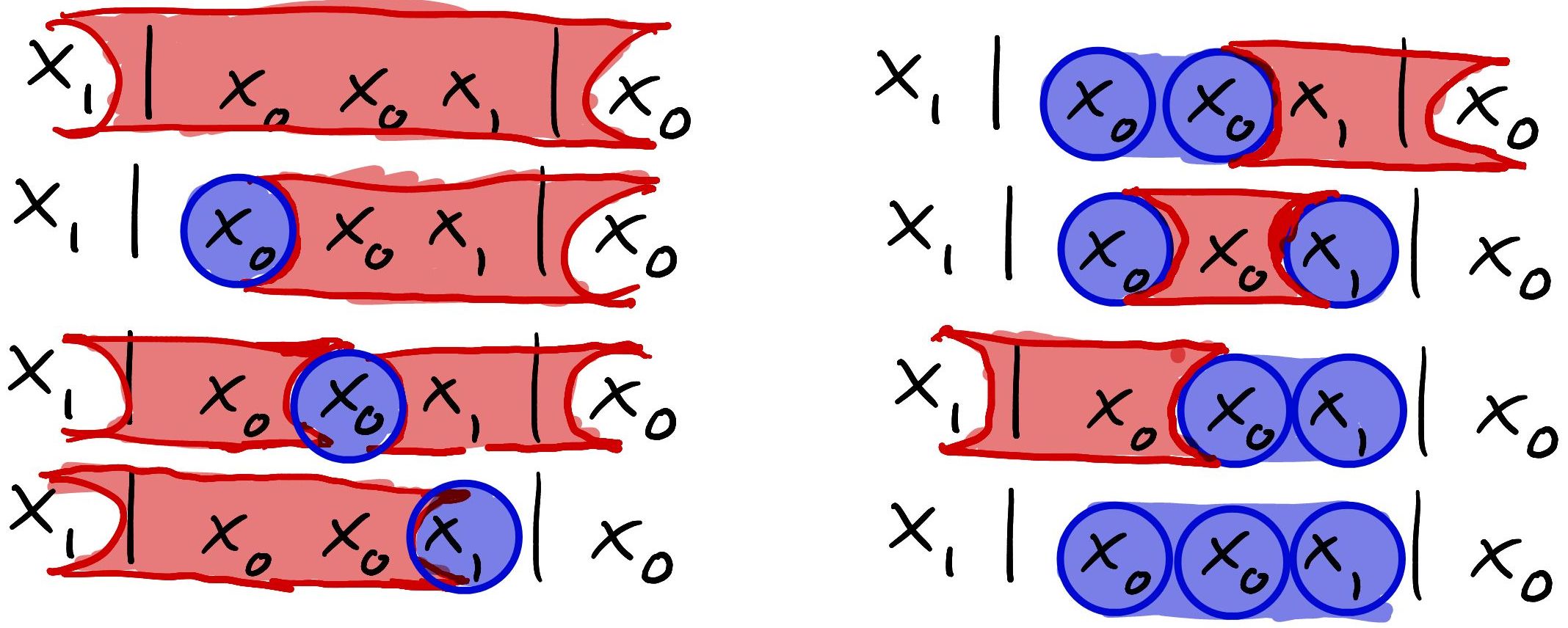}
\end{center}
Observe that half of the summands vanish by the convention given in Notation \ref{not:I} as at least one ''red worm'' starts and ends in the same letter. We thus have
\begin{align*}
\DeltaGon(x_0x_0x_1) &= x_0x_0x_1 \otimes \one - x_0 \otimes x_0x_1 + x_0 \otimes x_0x_1 + \one \otimes x_0x_0x_1 \\
&= x_0x_0x_1 \otimes \one + \one \otimes x_0x_0x_1.
\end{align*}
So $x_0^{n-1}x_1$ is primitive for $n=1,2,3$. 
However, for $n=4$ one similarly computes that
\begin{align*}
\DeltaGon(x_0x_0x_0x_1) &= x_0x_0x_0x_1 \otimes \one - 2x_0x_0 \otimes x_0x_1 + \one \otimes x_0x_0x_0x_1.
\end{align*}    
\end{Example}

\begin{Example}
We want to show the compatibility of $\DeltaGon$ and $\shuffle$ in a small example. It is  easy to see that
\begin{align*}
\DeltaGon(x_i) &= x_i \otimes \one + \one \otimes x_i \qquad (i=0,1),\\ 
\DeltaGon(x_0x_1) &= x_0x_1 \otimes \one + \one \otimes x_0x_1, \\
\DeltaGon(x_1x_0) &= x_1x_0\otimes \one + x_0 \otimes x_1 + x_1\otimes x_0 + \one \otimes x_1x_0.
\end{align*}
Now we compute
\begin{align*}
\DeltaGon(x_0 \shuffle x_1) &= \DeltaGon(x_0x_1 + x_1x_0)\\
&= (x_0x_1 + x_1x_0) \otimes\one + \one \otimes (x_0x_1 + x_1x_0) + x_0\otimes x_1 + x_1\otimes x_0\\
&= (x_0\otimes \one + \one\otimes x_0) \shuffle (x_1\otimes \one + \one \otimes x_1)\\
&= \DeltaGon(x_0) \shuffle \DeltaGon(x_1).
\end{align*}
\end{Example}

\subsection{The derivations \texorpdfstring{$D_w$}{Dw} for the Goncharov coproduct} \label{subsec:gon_deriva}

We make the derivations $D_w$, which we introduced in Subsection \ref{subsec:hopf_deriva} in a general context, explicit for the Goncharov coproduct.  
 
For simplicitly and to continue Brown's notation, we abbreviate
\[\LX=\Indec(\QX,\shuffle).\]
The weight, i.e., the number of letters, defines a grading on $(\QX,\shuffle)$ and hence also on the space $\LX$ of indecomposables. We denote the homogeneous subspaces of weight $w$ by $\QX_w$ resp. $\LXw{w}$. For each odd weight $2r+1$, $r\geq1$, the derivation $D_{2r+1}$ is given by
\[
\xymatrix{
D_{2r+1}: & \QX \ar[rr]^{\hspace{-3em}\DeltaGon'} && \bigoplus\limits_{w\geq1} \QX_w\otimes \QX \ar[rr]^{\hspace{0.5em} \pi_{2r+1}\otimes\id} && \LXw{2r+1}\otimes \QX.
}\]
Recall that $\DeltaGon'=\DeltaGon - \one\otimes\id$ and $\pi_{2r+1}:\bigoplus_{w\geq1} \QX_w\to \LXw{2r+1}$ is the canonical projection.

\begin{Proposition} \label{prop:shape_D}
Let $w = \ve_1\cdots \ve_n\in\X^*$, then each summand in $D_{2r+1} (w)$ corresponds to exactly one strict subword $u= \ve_{j+1}\dots\ve_{j+2r+1}$ with $0\le j \le n-2r-1$ of $w$.    
\end{Proposition}
\begin{proof} For each non-strict subword, the left tensor factor in the Goncharov coproduct \eqref{eq:Gon_formula} becomes a non-trivial product, and thus vanishes under the projection $\pi_{2r+1}$.
\end{proof}

\begin{Remark} \label{rem:D_quotientword}
Following the visualization of $\DeltaGon$ on a semicircle from Figure \ref{fig:semicircle}, we can depict the map $D_{2r+1}$ by omitting trivial factors of the following type
\def\scale{4}
\def\n{12}
\begin{center}
\begin{tikzpicture}[baseline = (current bounding box.north)]
	\coordinate (start) at (-\scale,0);
	\coordinate (end) at (\scale,0);
	
	\begin{scope}[nodes={fill=black, circle, scale=.5}]
		\clip (-\scale-.1,0) rectangle (\scale+.1,\scale+.1);
		\draw[thick] (start)+(end) circle (\scale);
	\end{scope}
	
	\draw[thick] (start) -- (end);
	
	\foreach \num in {0,1,...,\n}{
		\coordinate (L\num) at (180-\num*180/\n : 1.1*\scale);
		
		\coordinate (N\num) at (180-\num*180/\n : \scale);
	}
	
	\node at (L0) {$x_1$};
	\node at (L\n) {$x_0$};
	
	\node[blue] at (L1) {$\ve_{1}$};
	\node[blue] at (L3) {$\ve_{j-1}$};
	\node[blue] at (L4) {$\ve_j$};
	\node[red] at (L5) {$\ve_{j+1}$};
	\node[red] at (L7) {$\ve_{j+2r}$};
	\node[red] at (L8) {$\ve_{j+2r+1}$};
	\node[blue, anchor = south west] at (N9) {$\ve_{j+2r+2}$};
	\node[blue] at (L11) {$\ve_n$};
	
    \draw[thick] (N0) to[bend right = 15] (N1) 
        to (N3)
        to[bend right = 15] (N4)
        to (N9)
        to (N11)
        to[bend right = 15] (N12);

    \foreach \num in {1,3,4,9,11}{		
		\node[fill=black, circle, scale = .5] at (N\num) {};
	}
	
	\foreach \num in {0,5,7,8,12}{		
		\node[fill=black, circle, scale = .5] at (N\num) {};
	}		
	
	\node[red] at (L2) {$\dots$};
	\node[red] at (L6) {$\dots$};
	\node[red] at (L10) {$\dots$};
\end{tikzpicture}
\end{center}
\end{Remark}

\section{\texorpdfstring{The space $\B$ and  level lowering inspired by $\DeltaGon$}{Level lowering with d erivations} }

The subspace $\B \in \QX$ of words in $\two$ and $\three$ has a natural level filtration given by the number of occurences of the letter $\three$.  This subspace behaves nicely with respect to the Goncharov coproduct and 
we will show in this section that a variant of the derivations $D_{2r+1}$  made out of the Goncharov coproduct give rise to level lowering maps.

\subsection{Level filtration and level lowering maps}

\begin{Definition} \label{def:level}
We set 
\begin{equation*}
\B = \Q \langle x_0x_1,\, x_0x_0x_1 \rangle \subset \QX .
\end{equation*}	
The \emph{level} $\deg_3(w)$ of a word $w\in \B$ is given by the number of $x_0x_0x_1$ in $w$. This induces an ascending \emph{level filtration} on $\B$ given by 
	\begin{equation*}
		F_\ell \B = \operatorname{span}_\Q\Big\{  w \in \B \;\Big|\; \mbox{number of } x_0x_0x_1 \mbox{ in } w \leq \ell\Big\}
	\end{equation*}
	for all $\ell\in\N_0$.
\end{Definition}

\begin{Definition}\label{def:alphabet_23}
Consider the alphabet $\{\two,\three\}$ and let $\{\two,\three\}^*$ the set of all words with letters in $\{\two,\three\}$. Define the map
    \begin{align*}
        \bzd\colon \{\two,\three\}^* &\rightarrow \X^*
    \end{align*}
by $\bzd(\two) = x_0x_1$ and $\bzd(\three) = x_0x_0x_1$ and extend this with respect to concatenation. 

Obviously, the image $\bzd\big(\{\two,\three\}^*\big)\subset\X^*$ gives a $\Q$-basis for $\B$. 

The \emph{weight} of a word $u\in\{\two,\three\}^*$ is given by $\wt(u) = \wt\big(\bzd(u)\big)$, and similarly, the \emph{level} of $u\in \{\two,\three\}^*$ is given by $\deg_3(u)=\deg_3(\bzd(u))$.
\end{Definition}

\begin{Lemma} \label{lem:coaction_B23}
	We have
	\begin{equation*}
		\DeltaGon\colon \B \rightarrow \QX  \otimes \B.
	\end{equation*}
\end{Lemma}
 
\begin{proof}	
	Let $w= \ve_1\cdots\ve_n \in \B$ be a word. 
	It suffices to show, if $ v=\ve_{i_1}\cdots\ve_{i_k}$ is subword of $w$, which determines a non-trivial contribution in the Goncharov coproduct \eqref{eq:Gon_formula}, then $v\in\B$.
    It is $\ve_{i_1} = x_0$ and $\ve_{i_k}=x_1$, since otherwise the factors 
    $\Iformal(x_1; \ve_1\cdots \ve_{i_1-1}; \ve_{i_1})$ and $\Iformal(\ve_{i_k}; \ve_{i_k+1}\cdots \ve_{n} ; x_0)$ vanish. 
    If $x_0x_0x_0$ is a strict subword of $\ve_{i_1}\cdots\ve_{i_k}$,  then there is a vanishing factor $\Iformal(x_0;u;x_0)$.
    Indeed, there is a $j\in \{1,\dots,k-2\}$ such that
    $\ve_{i_j}\ve_{i_{j+1}}\ve_{i_{j+2}} = x_0x_0x_0$ and then at least one of factors $\Iformal(\ve_{i_j}; \ve_{i_j+1}\cdots\ve_{i_{j+1}-1};\ve_{i_{j+1}})$ and $\Iformal(\ve_{i_{j+1}}; \ve_{i_{j+1}+1}\cdots\ve_{i_{j+2}-1}; \ve_{i_{j+2}})$ vanishes since $w\in\B$ implies that $x_0x_0x_0$ is not a strict subword of $w$ and thus at least one of the quotient words is not empty (cf Notation \ref{not:Gon_coprod}). 
    One similarly shows that, if $x_1x_1$ is a strict subword of $\ve_{i_1}\cdots\ve_{i_k}$, then we obtain a vanishing factor $\Iformal(x_1;u;x_1)$. 
    Hence any  subword $ \ve_{i_1}\cdots\ve_{i_k} \notin\B$ contributes trivially and the claim follows.
\end{proof}

In particular, by Lemma \ref{lem:coaction_B23} the Goncharov coproduct restricts to
	\begin{equation*}
		\DeltaGon\colon F_\ell \B \rightarrow \QX \otimes F_\ell\B.
	\end{equation*}

\begin{Definition} \label{def:Dpart} For each $r\in\NN$, define the map $\partial_{2r+1}:\QX\to \QX_{2r+1}\otimes \QX$ by
\[
\partial_{2r+1}(w) =
\sum_{j=0}^{N-2r-1}  \Iformal(\ve_j;\ve_{j+1}\cdots \ve_{j+ 2r+1}; \ve_{j+2r+2}) \otimes  \ve_1\cdots\ve_j \ve_{j+2r+2}\cdots \ve_N
\]
for a word $w=\ve_1\cdots \ve_N\in \X^*$ and $\ve_0 = x_1$ and $\ve_{N+1} = x_0$.
\end{Definition}

\begin{Remark}\label{rem:diagr_partial_D}
Let $r\in\NN$. By Proposition \ref{prop:shape_D}, the linear map $\partial_{2r+1}$ fits into the following commutative diagram
\[
\xymatrix{
\QX \ar[rrr]^{\partial_{2r+1}} \ar[ddrrr]^{D_{2r+1}} &&&
\QX_{2r+1} \otimes \QX \ar[dd]^{\pi_{2r+1}\otimes \operatorname{id}} \\
& \\
&&& \LX_{2r+1}  \otimes \QX.
}
\]
\end{Remark}

The maps $\partial_{2r+1}$ introduced in 
Definition \ref{def:Dpart} reduce the level.
\begin{Lemma} \label{lem:D_reduces_level_B23} 
	For all $r,\ell\in\N$ we have
	\begin{equation*}
		\partial_{2r+1}(F_\ell \B) \subseteq \QX_{2r+1} \otimes F_{\ell-1} \B.
	\end{equation*}
\end{Lemma}

\begin{proof}
	Let $w= \ve_1 \cdots \ve_n \in F_\ell \B$. By Definition \ref{def:Dpart}, each summand in $\partial_{2r+1} (w)$ corresponds to a strict subword $u = \ve_{j+1}\cdots\ve_{j+2r+1}$ with $0 \le j  \le n-2r-1$ of $w$.
	If the enlarged word $\hat{u} = \ve_{j}\cdots\ve_{j+2r+2}$ does not contain $x_0x_0$ as a strict subword then it has to be $(x_0x_1)^{r+1} x_0$ or 
	$x_1 (x_0x_1)^{r+1}$ since it cannot contain $x_1x_1$ and we always have $\ve_0\ve_1= x_1 x_0 $ and $\ve_{n}\ve_{n+1}=x_1x_0$ by definition.
	In either case, we have $\ve_{j} = \ve_{j+2r+2}$ and therefore the factor $\Iformal(\ve_{j};u;\ve_{j+2r+2})$ vanishes. 
	We deduce that if the strict subword $u$ contributes to $\partial_{2r+1}(w)$, then $\hat{u}$ contains $x_0x_0$ as a strict subword and we have $\ve_j \neq \ve_{j+2r+2}$. Hence $\ve_1\cdots\ve_n\in F_\ell \B$ implies that $\ve_1\cdots\ve_j\ve_{j+2r+2}\cdots\ve_n\in F_{\ell-1}\B$.
\end{proof}

\begin{Definition} \label{def:ass_graded_B23}
	For all $\ell\in\N$ we set
	\begin{equation*}
		\gr_\ell^F(\B) = \faktor{F_\ell \B}{F_{\ell-1} \B}
	\end{equation*}
	and  $\gr_0^F= F_0\B$. For all $\ell\in\N$ we denote by $\pi_\ell^F\colon F_\ell \B \twoheadrightarrow \gr_\ell^F\B$
the natural projections. 
\end{Definition}	
	
By Lemma \ref{lem:D_reduces_level_B23}, the restricted map $(\operatorname{id}\otimes \pi_{\ell-1}^F)\circ \partial_{2r+1}\Big\rvert_{F_\ell\B}$ induces a map
	\begin{equation*}
		\partial_{2r+1}^{(\ell)} \colon \gr_\ell^F \B \rightarrow \QX_{2r+1} \otimes \gr_{\ell-1}^F \B
	\end{equation*}
for all $r,\ell\geq 1$.	

\begin{Definition} \label{def:B_1} Define the sub vector space
	\[
	\mathcal{B}^{1} = \gr_1^F \B  \oplus \,  (\gr_0^F \B)x_0  \subset \QX.
	\]
\end{Definition}

Explicitly, we have
\[\mathcal{B}^{1}=\operatorname{span}_\Q\{\pi_{1}^F((x_0x_1)^{r_1}x_0x_0x_1(x_0x_1)^{r_2})\mid r_1,r_2\geq0\}\oplus\operatorname{span}_\Q\{(x_0x_1)^r x_0 \mid r\geq0.\}\]
We have the following refinement of Lemma \ref{lem:D_reduces_level_B23}.

\begin{Theorem}\label{thm:graded_image_D_B23} 
Let $r,\ell\geq 1$. Then for each weight $N\geq 2r+1$, we have 
\begin{equation}\label{eq:filtered_partial_2r+1}
		\partial_{2r+1}^{(\ell)}\colon \gr_{\ell}^F (\B)_N \rightarrow \mathcal{B}^{1} \otimes \gr_{\ell-1}^F (\B)_{N-2r-1}.
	\end{equation}
\end{Theorem}
  
\begin{proof}
	Let $w= \ve_1 \cdots \ve_N \in\gr_\ell^F (\B)_N$ and $u = \ve_{j+1}\cdots\ve_{j+2r+1}$ with $0 \le j  \le N-2r-1$, be a strict subword of $w$.
	Assume the enlarged word $\hat{u}=\ve_{j}\cdots\ve_{j+2r+2}$ contains $x_0x_0$ at least two times as a strict subword. Then the summand in $\partial_{2r+1}(w)$ 
	correspoding to $u$ is contained in $\QX_{2r+1}\otimes F_{\ell-2}^F \B$ and thus vanishes under $\operatorname{id}\otimes \pi_{\ell-1}^F$.
    On the other hand, if the choice of $u$ contributes non-trivially to $\partial_{2r+1}(w)$ then it must contain $x_0x_0$ at least once as a strict subword 
    by Lemma \ref{lem:D_reduces_level_B23}.
	So we conclude, if $u$ contributes non-trivially, then it contains $x_0x_0$ exactly once as a strict subword. 
	The remaining cases for $u$ determine the factors 
 		\begin{enumerate}
		\item $\Iformal\big(x_1; (x_0x_1)^{r_1}  x_0 x_0 x_1   (x_0x_1)^{r_2};x_0\big) =  (x_0x_1)^{r_1}   x_0x_0x_1    (x_0x_1)^{r_2} $,
		\item $\Iformal\big(x_0; (x_1x_0)^{r_1} x_1 x_0 x_0   (x_1x_0)^{r_2} ;x_1\big) =  - (x_0x_1)^{r_2}   x_0x_0x_1   (x_0x_1)^{r_1} $,
		\item $\Iformal\big(x_1;(x_0x_1)^r  x_0; x_0\big) = (x_0x_1)^r x_0 $,
		\item $\Iformal\big(x_0;  (x_0x_1)^r  x_0;x_1\big)= -(x_0x_1)^r  x_0$.
	\end{enumerate}
	 By definition, all of these factors are mapped to elements of $\mathcal{B}^{1}$. Furthermore, we deduce from $\ve_j\neq \ve_{j+2r+2}$ that $\ve_1\cdots\ve_j\ve_{j+2r+2}\cdots\ve_N\in\gr_{\ell-1}^F (\B)_{N-2r-1}$ and the claim follows.
 \end{proof}

\begin{Remark} \label{rem:terms_in_partial_1}
The cases 1.-4. from the proof of Theorem \ref{thm:graded_image_D_B23} can be depicted as follows:
\def\scale{3}
	\def\n{23}
	\begin{center}
	\begin{tikzpicture}[baseline = (current bounding box.north)]
		\coordinate (start) at (-\scale,0);
		\coordinate (end) at (\scale,0);
		
		\begin{scope}[nodes={fill=black, circle, scale=.5}]
			\clip (-\scale-.1,0) rectangle (\scale+.1,\scale+.1);
			\draw[thick] (start)+(end) circle (\scale);
		\end{scope}
	
		\node[scale=.9] at (-\scale,\scale) {case 1.};
		
		\draw[thick] (start) -- (end);
		
		\foreach \num in {0,1,...,\n}{
			\coordinate (L\num) at (180-\num*180/\n : 1.1*\scale);
			
			\coordinate (N\num) at (180-\num*180/\n : \scale);
		}
  
		\draw[thick, bend right = 45] (N0) to (N1)
			to (N2)
			to (N3)
			to (N4)
			to (N5)
			to (N6)
			to (N7);
		
		\draw[thick] (N7) to node[below,pos=.95] {} (N15);
		
		\draw[thick, bend right = 45] (N15) to (N16)
			to (N17)
			to (N18)
			to (N19)
			to (N20)
			to (N21)
			to (N22)
			to (N23);
			
		\foreach \num in {0,1,3,4,5,...,18,19,21,22,\n}{
			\node[fill=black, circle, scale=.5] at (N\num) {};
		}
		
		\node at (L0) {$x_1$};
		\node at (L\n) {$x_0$};
		
		\foreach \num in {9,12,14}{
			\node[red] at (L\num) {$x_1$};
		}
		
		\foreach \num in {3,5,7,16,18,22}{
			\node[blue] at (L\num) {$x_1$};
		}
		
		\foreach \num in {8,10,11,13}{
			\node[red] at (L\num) {$x_0$};
		}
		
		\foreach \num in {1,4,6,15,17,19,21}{
			\node[blue] at (L\num) {$x_0$};
		}
		
		\node[blue] at (L2) {$\dots$};
		\node[blue] at (L20) {$\dots$};
	\end{tikzpicture}
	\begin{tikzpicture}[baseline = (current bounding box.north)]
		\coordinate (start) at (-\scale,0);
		\coordinate (end) at (\scale,0);
		
		\begin{scope}[nodes={fill=black, circle, scale=.5}]
			\clip (-\scale-.1,0) rectangle (\scale+.1,\scale+.1);
			\draw[thick] (start)+(end) circle (\scale);
		\end{scope}
		
		\node[scale=.9] at (\scale,\scale) {case 2.};
		
		\draw[thick] (start) -- (end);
		
		\foreach \num in {0,1,...,\n}{
			\coordinate (L\num) at (180-\num*180/\n : 1.1*\scale);
			
			\coordinate (N\num) at (180-\num*180/\n : \scale);
		}
		
		\draw[thick, bend right = 45, dotted] (N0) to (N1)
			to (N2)
			to (N3)
			to (N4)
			to (N5)
			to (N6);
		
		\draw[thick, dotted] (N6) to node[below,pos=.95] {} (N14);
		
		\draw[thick, bend right = 45, dotted] (N14) to (N15)
			to (N16)
			to (N17)
			to (N18)
			to (N19)
			to (N20)
			to (N21)
			to (N22)
			to (N23);
		
		\foreach \num in {0,1,3,4,5,...,18,19,21,22,\n}{
			\node[fill=black, circle, scale=.5] at (N\num) {};
		}
		
		\node at (L0) {$x_1$};
		\node at (L\n) {$x_0$};
		
		\foreach \num in {7,9,12}{
			\node[red] at (L\num) {$x_1$};
		}
		
		\foreach \num in {3,5,14,16,18,22}{
			\node[blue] at (L\num) {$x_1$};
		}
		
		\foreach \num in {8,10,11,13}{
			\node[red] at (L\num) {$x_0$};
		}
		
		\foreach \num in {1,4,6,15,17,19,21}{
			\node[blue] at (L\num) {$x_0$};
		}
		
		\node[blue] at (L2) {$\dots$};
		\node[blue] at (L20) {$\dots$};
	\end{tikzpicture}
	\begin{tikzpicture}[baseline = (current bounding box.north)]
		\coordinate (start) at (-\scale,0);
		\coordinate (end) at (\scale,0);
		
		\begin{scope}[nodes={fill=black, circle, scale=.5}]
			\clip (-\scale-.1,0) rectangle (\scale+.1,\scale+.1);
			\draw[thick] (start)+(end) circle (\scale);
		\end{scope}
		
		\node[scale=.9] at (-\scale,\scale) {case 3.};
		
		\draw[thick] (start) -- (end);
		
		\foreach \num in {0,1,...,\n}{
			\coordinate (L\num) at (180-\num*180/\n : 1.1*\scale);
			
			\coordinate (N\num) at (180-\num*180/\n : \scale);
		}
		
		\draw[thick, dashed, bend right = 45] (N0) to (N1)
			to (N2)
			to (N3);
		
		\draw[thick, dashed] (N3) -- node[below,pos=.2] {} (N11);
		
		\draw[thick, dashed, bend right = 45] (N11) to (N12)
			to (N13)
			to (N14)
			to (N15)
			to (N16)
			to (N17)
			to (N18)
			to (N19)
			to (N20)
			to (N21)
			to (N22)
			to (N23);
		
		\foreach \num in {0,1,3,4,5,...,18,19,21,22,\n}{
			\node[fill=black, circle, scale=.5] at (N\num) {};
		}
		
		\node at (L0) {$x_1$};
		\node at (L\n) {$x_0$};
		
		\foreach \num in {5,7,9}{
			\node[red] at (L\num) {$x_1$};
		}
		
		\foreach \num in {3,12,14,16,18,22}{
			\node[blue] at (L\num) {$x_1$};
		}
		
		\foreach \num in {4,6,8,10}{
			\node[red] at (L\num) {$x_0$};
		}
		
		\foreach \num in {1,11,13,15,17,19,21}{
			\node[blue] at (L\num) {$x_0$};
		}
		
		\node[blue] at (L2) {$\dots$};
		\node[blue] at (L20) {$\dots$};
	\end{tikzpicture}
	\begin{tikzpicture}[baseline = (current bounding box.north)]
		\coordinate (start) at (-\scale,0);
		\coordinate (end) at (\scale,0);
		
		\begin{scope}[nodes={fill=black, circle, scale=.5}]
			\clip (-\scale-.1,0) rectangle (\scale+.1,\scale+.1);
			\draw[thick] (start)+(end) circle (\scale);
		\end{scope}
		
		\node[scale=.9] at (\scale,\scale) {case 4.};
		
		\draw[thick] (start) -- (end);
		
		\foreach \num in {0,1,...,\n}{
			\coordinate (L\num) at (180-\num*180/\n : 1.1*\scale);
			
			\coordinate (N\num) at (180-\num*180/\n : \scale);
		}

		\draw[thick, dash dot, bend right = 45] (N0) to (N1)
			to (N2)
			to (N3)
			to (N4)
			to (N5)
			to (N6)
			to (N7)
			to (N8)
			to (N9)
			to (N10);		
		
		\draw[thick, dash dot] (N10) -- node[below right,pos=.7] {} (N18);
		
		\draw[thick, dash dot, bend right = 45] (N18) to (N19)
			to (N20)
			to (N21)
			to (N22)
			to (N23);
		
		\foreach \num in {0,1,3,4,5,...,18,19,21,22,\n}{
			\node[fill=black, circle, scale=.5] at (N\num) {};
		}
		
		\node at (L0) {$x_1$};
		\node at (L\n) {$x_0$};
		
		\foreach \num in {12,14,16}{
			\node[red] at (L\num) {$x_1$};
		}
		
		\foreach \num in {3,5,7,9,18,22}{
			\node[blue] at (L\num) {$x_1$};
		}
		
		\foreach \num in {11,13,15,17}{
			\node[red] at (L\num) {$x_0$};
		}
		
		\foreach \num in {1,4,6,8,10,19,21}{
			\node[blue] at (L\num) {$x_0$};
		}
		
		\node[blue] at (L2) {$\dots$};
		\node[blue] at (L20) {$\dots$};
		
	\end{tikzpicture}
	\end{center}
  
\end{Remark}
\vspace{0,7cm}
\begin{Remark} \label{rem:terms_in_partial_2}
For any word $w \in \B$ there is always  one $x_0$ left to a letter $x_1$. Therefore, if in case 1. from the  proof  of Theorem \ref{thm:graded_image_D_B23}
\[
    \hat{u} = \ve_j  \ve_{j+1} \cdots \ve_{j+2r+2} =  x_1 (x_0x_1)^{r_1} x_0 x_0 x_1 (x_0x_1)^{r_2}x_0
\]
is a strict subword of $w$ or equivalently if $j>1$,   
then we get also a case 2. contribution by
\[
    \hat{u}' = \ve_{j-1}  \ve_{j} \cdots \ve_{j+2r+1} =  x_0 (x_1 x_0)^{r_1} x_1x_0 x_0  (x_1x_0)^{r_2}  x_1.
\]
Conversely, each $\hat{u}'$ in case 2. determines a possibly enlarged strict subword $s$, as right to a letter $x_1$ must  always be a $x_0$. 

The dotted and dashed lines in the above picture depict such a pair of contributing subwords $\hat u$ and $\hat u'$.
\end{Remark}
  
Because of the commutative diagram in Remark \ref{rem:diagr_partial_D}, Lemma \ref{lem:D_reduces_level_B23} and Theorem \ref{thm:graded_image_D_B23} imply the following for the derivations $D_{2r+1}$.

 \begin{Proposition}
     \label{prop:level_red_derivations_X} The following holds.
    \begin{enumerate}
        \item For all $\ell,r\in\N$ we have
        \begin{equation*}
            D_{2r+1}(F_\ell \B) \subseteq \LX_{2r+1} \otimes F_{\ell-1} \B.
        \end{equation*}
        \item For all $\ell,r\in\N$ and each weight $N\geq 2r+1$ we have derivations
        \begin{equation*}
            D_{2r+1}^{(\ell)}\colon \gr_\ell^F(\B)_N \rightarrow \pi_{2r+1}\big(\mathcal{B}^1\big) \otimes \gr_{\ell-1}^F (\B)_{N-2r-1}.
        \end{equation*}
    \end{enumerate}
\end{Proposition}

\subsection{The linear map \texorpdfstring{$\partialphi$}{phi} is an isomorphism}

Let $N,\ell,r\in\N$ be positive integers throughout the rest of this subsection.

\begin{Definition}\label{def:B_1_to_c_w}
    We define a linear map
    \begin{align*}
        \phi\colon  \mathcal{B}^1 &\to \Q, \\
        \beta&\mapsto c_\beta = 
        \begin{cases}
        	c_{a,b} & \mbox{ if } \beta  \in \gr_1^F \big(\B\big) \\
        	2\cdot (-1)^n  & \mbox{ if }  \beta \in \gr_0^F(\B) x_0. 
        \end{cases}
    \end{align*}
   where
   \begin{equation*}
   	c_{a,b}^r = 2\cdot (-1)^r \left(\binom{2r}{2b+2} - (1-2^{-2r}) \binom{2r}{2a+1}\right).
   \end{equation*}
\end{Definition}

The behaviour of the numbers $c_{a,b}^r$ will be discussed in detail in Subsection \ref{subsec:numbers_c}.

Composing the map
\begin{equation*} 
		\partial_{2r+1}^{(\ell)}\colon \gr_{\ell}^F (\B)_N \rightarrow \mathcal{B}^{1} \otimes \gr_{\ell-1}^F (\B)_{N-2r-1}.
\end{equation*}
from Theorem \ref{thm:graded_image_D_B23} with $\phi$ on the first factor yields a map
\begin{align*}
    \partialphigr\colon \gr_\ell^F(\B)_N &\rightarrow \gr_{\ell-1}^F(\B)_{N-2r-1}\\
    w&\mapsto \phi \big(\partial_{2r+1}^{(\ell)} (w)\big).
\end{align*}
Here, we identified $\QQ\otimes \gr_{\ell-1}^F(\B)_{N-2r-1} \simeq \gr_{\ell-1}^F(\B)_{N-2r-1}$.

\begin{Definition}\label{def:filtered_partial_leN}
    For $N\geq 2r+1$ and $\ell\geq 1$ we define the linear map
     \begin{align*}
        \partialphi\colon \gr_\ell^F(\B)_N &\rightarrow \bigoplus_{3\leq 2r+1\leq N} \gr_{\ell-1}^F (\B)_{N-2r-1}\\
        w &\mapsto \sum_{3\leq 2r+1\leq N} \partialphigr(w).
     \end{align*}
\end{Definition}

Our goal is to show Theorem \ref{thm:matrix_invertible}, which says that $\partialphi$ is an isomorphism. 
We begin by fixing a basis for the domain and codomain of $\partialphi$, respectively.

\begin{Definition}\label{def:bases_partialphi}
    Recall the map $\bzd\colon \{\two,\three\}^*\to \X^*$ was given in Definition \ref{def:alphabet_23}. We set
    \begin{align*}
        \BNgr &= \bzd\big(\{u\in\{\two,\three\}^*\mid \wt(u) = N \text{ and }  \deg_3(u)= \ell \}\big), \\
        \BNfil &= \begin{cases}
        \bzd\big(\{u\in\{\two,\three\}^*\mid  \wt(u) < N-1, \,\,  \deg_3(u)=\ell-1  \}\big), & N\equiv \ell \bmod 2 \\
        \emptyset & N\not\equiv \ell \bmod 2
           \end{cases}.
    \end{align*}
\end{Definition}

\begin{Remark}\label{rem:parity_B23_level}
    It is easy to see that also $\BNgr$ is empty, if $N\not\equiv \ell \mod 2$.  \\
    If $\BNfil\neq\emptyset$, then it is the disjoint union of lower weights, i.\,e.
    \begin{equation*}
        \BNfil = \mathop{\dot\bigcup}_{n=0}^{N-3} \BNgr[n][\ell-1].
    \end{equation*}
\end{Remark}

\begin{Lemma}
\label{lem:B_N_bijection} 
    The map
    \begin{align*}
        \psi\colon \BNfil &\rightarrow \BNgr\\
        \bzd(u) &\mapsto \bzd\big(\two^{r-1}\three u\big)
    \end{align*}
    where $2r = N-1-\wt(u)$ is a bijection.
\end{Lemma}

\begin{proof}
    If $\BNfil=\emptyset$, we also have $\BNgr = \emptyset$ and vice versa by Remark \ref{rem:parity_B23_level}.
    So we can assume without loss of generality that $\BNfil \neq\emptyset$. Let $\bzd(u)\in \BNfil$. 
    Indeed, $\bzd\big(\two^{r-1}\three u\big)\in\BNgr$ by the choice of $r$. 
    Since $\ell\geq1$, all words in $\BNgr$ can be uniquely written as $\bzd(\two^{s-1}\three u')$ for some $s\in\N$ and a word $u'\in\X^*$ with $\bzd(u')\in\BNfil$. 
    The assignment $\bzd(\two^{s-1}\three u') \mapsto \bzd(u')$ clearly gives an inverse map to $\psi$. Hence $\psi$ is a bijection, since both sets are finite. 
\end{proof}

\begin{Corollary}\label{cor:bases_same_size}
    We have
    \begin{equation*}
        \# \BNgr  = \# \BNfil.
    \end{equation*}
    \qed
\end{Corollary}

\begin{Remark}
    The statement in Corollary \ref{cor:bases_same_size} can also be proven directly by combinatorial means. 
    We again disregard the case where either set is empty (cf. Remark \ref{rem:parity_B23_level}). 
    Then we have $N = 2m + 3\ell$ for some non-negative $m\in\mathbb{Z}$ and observe that
    \begin{equation*}
		\# \BNgr  = \binom{m+\ell}{\ell} = \sum_{0\leq m'\leq m} \binom{m'+\ell-1}{\ell-1} =  \#\BNfil ,
	\end{equation*}
	where the second equality is the so-called hockey-stick identity of binomial coefficients.
\end{Remark}

\begin{Definition}\label{def:B_N_ordering}
 We endow   $\BNgr$ with the  lexicographic order $\leq_{\text{lex}}$  with respect to the order $\two < \three$. With this we define an order on 
 $\BNfil$ by $\bzd(u) \leq \bzd(u')$, iff\,
$\psi(\bzd(u)) \leq_{\text{lex}} \psi(\bzd(u)')$.   
\end{Definition}  

In other words, in Definition \ref{def:B_N_ordering} we require that the bijection $\psi$ is an order preserving map. Observe we get 
$\bzd(u)\leq \bzd(u')$ if and only if $\wt(u)<\wt(u')$ or if $\wt(u) = \wt(u')$ and  $\bzd(u) \leq_{\text{lex}}  \bzd(u')$
The same order is used in \cite{bgf}. 

\begin{Remark}
 Brown \cite{br} uses a similar lexicographic order on $\BNgr$, but a different order on $\BNfil$. 
The order on $\BNfil$ can be recovered from the lexicographic order on $\BNgr$ via a variation   of the map $\psi\colon \BNgr \to \BNfil$ 
from Lemma \ref{lem:B_N_bijection} given by  
$u \mapsto u\three\two^{s-1}$. 
This difference in Brown's setup is due to the inverted MZV-notation.   
\end{Remark}

\begin{Proposition}\label{prop:partialphi_gr_even_coeffs}
    Let $w\in\BNgr$, $r\geq1$. If there is a decomposition $w=\beta v$ with $\beta \in\BNgr[2r+1][1]$, then 
    \begin{equation*}
        \partialphigr(w) = c_\beta\, v + (\text{terms with coefficients in $2\mathbb{Z}$}).
    \end{equation*}
    Otherwise, $\partialphigr(w)$ is given by a sum of terms with coefficients in $2\mathbb{Z}$.
\end{Proposition}

\begin{proof}
    Recall the four types of strict subwords $w$ that contribute non-trivially to $\partial_{2r+1}^{(\ell)}(w)$ from the proof of Theorem \ref{thm:graded_image_D_B23}. \\
    Unless $w = \beta v$, the first two cases come in pairs $\beta,\beta'$ by Remark \ref{rem:terms_in_partial_2}. They contribute with $c_\beta = c_{a,b}$ and $c_{\beta'} = - c_{b,a}$ for some non-negative $a,b\in\mathbb{Z}$. By Lemma \ref{lem:properties_c_beta} 2. we have $c_\beta + c_{\beta'}\in2\mathbb{Z}$.
    \\
    The last two cases contribute with the even coefficients $\pm 2$.
\end{proof}

We deduce for the linear maps $\partialphi$ from 
Definition \ref{def:filtered_partial_leN}
the following from Proposition \ref{prop:partialphi_gr_even_coeffs}.

\begin{Proposition}\label{prop:partialphi_even_coeffs}
    For $w\in \BNgr$ we have
    \begin{equation*}
        \partialphi(w) = \sum_{\substack{w = \beta v, \\ \beta\in\BNgr[2r+1][1],r\geq1}} c_\beta\, v + (\text{terms with coefficients in }2\mathbb{Z}).
    \end{equation*}
 \qed
\end{Proposition}

\begin{Definition}\label{def:M_N_matrix}
Let $N,\ell\in\N$ with $N = 2m + 3\ell$ for some non-negative $m\in\mathbb{Z}$. We define
   \begin{equation*}
        \MB = \big(m_{w,w'}\big)_{w\in \BNgr,\, w'\in \BNfil} 
    \end{equation*}
to be the representing matrix of $\partialphi$ with respect to the bases $\BNgr$ and $\BNfil$ from Definition \ref{def:B_N_ordering}.
\end{Definition}

Precisely, $m_{w,w'}\in\Q$ are for $w\in \BNgr$ defined by the equations
\[	\partialphi(w) = \sum_{w'\in \BNfil} m_{w,w'}\, w'.\]

We now deduce the following from Corollary \ref{cor:bases_same_size}.

\begin{Corollary}
    The matrix $\MB$ is quadratic. \qed
\end{Corollary}

Recall that our goal is to show that the level reducing map $\partialphi$ is an isomorphism. 
By Remark \ref{rem:parity_B23_level}, it suffices to study the non-trivial cases where $N = 2m+3\ell$ for some non-negative $m\in\mathbb{Z}$. 
The key idea of the proof is showing that the transformation matrices $\MB$ are invertible via the following lemma. Therefore, recall the $p$-adic valuation $\nu_p$ of a rational number also given in Definition \ref{def:p-adic_valuation}.

\begin{Lemma}\label{lem:p-adic_invertible}
    Let $p$ be a prime number and $n\in\N$. Let $A = (a_{ij})$ be an $n\times n$ matrix with entries in $\Q$ such that
    \begin{enumerate}
        \item[(i)] $\nu_p(a_{ij}) \geq 1$ for all $i>j$ and
        \item[(ii)] $\nu_p(a_{jj}) = \min\{\nu_p(a_{ij}) \mid i\in\{1,\dots,n\}\}\leq 0$ for all $j\in\{1,\dots,n\}$.
    \end{enumerate}
    Then $A$ is invertible.
\end{Lemma}

\begin{proof}
    To show that $\det(A)\neq 0$ it suffices to show that $\det(A')\neq 0$ where the matrix $A' = (a_{ij}')$ arises from $A$ by multiplying the $j$-th column with $p^{-\nu_p(a_{jj})}$ for all $j\in\{1,\dots,n\}$. 
    Since $\nu_p(a_{jj})\leq 0$, condition (i) still holds for $A'$. Condition (ii) implies that $\nu_p(a_{ij}')\geq 0$ for all entries of $i,j\in\{1,\dots,n\}$.
    By construction of $A'$, we have $\nu_p(a_{jj}') = 0$. So in particular, the entries on the diagonal are not zero. 
    It follows that $A'$ modulo $p$ is an upper triangular matrix with non-zero entries on the diagonal and is thus invertible.
\end{proof}

\begin{Theorem}[Brown \cite{br}]\label{thm:matrix_invertible}
 The map $\partialphi$ is an isomorphism of vector spaces.
\end{Theorem}
	
\begin{proof}
    We show that the matrix $\MB$ of the operator $\partialphi$, given in Definition \ref{def:M_N_matrix}, is invertible, since it satisfies the assumptions of Lemma \ref{lem:p-adic_invertible} for $p=2$.  \\
    Let $v\in\BNfil$. Recall the bijection $\psi\colon \BNfil\to \BNgr$ from Lemma \ref{lem:B_N_bijection} which is order-preserving by Definition \ref{def:B_N_ordering}. Let $2r = N-1-\wt(v)$ and let $w\in\BNgr$ with $w = \beta v$ for some $\beta\in \BNgr[2r+1][1]$.
    If $w\neq \psi(v)$, then $\beta$ does not end in $x_0x_0x_1$, hence $w > \psi(v)$. 
    We deduce from Proposition \ref{prop:partialphi_even_coeffs} that the entries of $\MB$ that are not in $2\mathbb{Z}$ are either on or above the main diagonal due to the orders on $\BNgr$ and $\BNfil$. This implies condition (i) from Lemma \ref{lem:p-adic_invertible}. \\
    The entries on the main diagonal of $\MB$ are given by the coefficients of $v$ in $\partialphi(\psi(v))$ with $r$ as above. By Lemma \ref{lem:properties_c_beta} 3. and Proposition \ref{prop:partialphi_even_coeffs}, these entries have a non-positive $2$-adic valuation and realize the minimum of this valuation within its column. Hence $\MB$ satisfies condition (ii) from Lemma \ref{lem:p-adic_invertible}.
\end{proof}

We illustrate the content of this section in two examples.

\begin{Example}
    In our first example, we compute the matrix $\MB[9][1]$ corresponding to the map $\partialphi[9][1]$.
    
    We have
    \begin{align*}
        \BNgr[9][1] &= \{\bzd(3,2,2,2), \bzd(2,3,2,2), \bzd(2,2,3,2), \bzd(2,2,2,3)\}
    \shortintertext{and}
        \BNfil[9][0] &= \{\bzd(2,2,2), \bzd(2,2), \bzd(2), \one\}.
    \end{align*}
    In order to compute, e.\,g., $\partialphi[9][1](\bzd(3,2,2,2))$ we first observe that
    \begin{align*}
        \partial_{3}^{(1)}\big(\bzd(3,2,2,2)\big) &= \Iformal(x_1;x_0x_0x_1;x_0) \otimes \bzd(2,2,2) + \Iformal(x_0;x_0x_1x_0;x_1) \otimes \bzd(2,2,2)\\
        &= (x_0x_0x_1 - x_0x_1x_0) \otimes \bzd(2,2,2).
    \end{align*}
    These terms can be depicted, respectively, as follows:
\begin{center}
	\def\scale{2.5} 
	\def\n{10} 
	\begin{tikzpicture}[baseline = (current bounding box.north)]
		\coordinate (start) at (-\scale,0);
		\coordinate (end) at (\scale,0);
		
		\begin{scope}[nodes={fill=black, circle, scale=.5}]
			\clip (-\scale-.1,0) rectangle (\scale+.1,\scale+.1);
			\draw[thick] (start)+(end) circle (\scale);
		\end{scope}
		
		\draw[thick] (start) -- (end);
		
		\foreach \arc in {0,1,...,\n}{			
			\coordinate (L\arc) at (180-\arc*180/\n : 1.2*\scale);
			
			\coordinate (N\arc) at (180-\arc*180/\n : \scale);
		}
		
		
		\foreach \num in {0,1,...,\n}{
			\node[fill=black, circle, scale=.5] at (N\num) {};
		}
		
		\foreach \num in {1,2}{
			\node[red] at (L\num) {$x_0$};
		}
		\foreach \num in {4,6,8}{
			\node[blue] at (L\num) {$x_0$};
		}
		
		\node[red] at (L3) {$x_1$};	
		
		\foreach \num in {5,7,9}{
			\node[blue] at (L\num) {$x_1$};
		}

		\node at (L0) {$x_1$};
		\node at (L\n) {$x_0$};

        \draw[thick] (N0) to (N4);
        \draw[thick, bend right = 15] (N4) to (N5)
            to (N6)
            to (N7)
            to (N8)
            to (N9)
            to (N10);
	\end{tikzpicture}
	\begin{tikzpicture}[baseline = (current bounding box.north)]
		\coordinate (start) at (-\scale,0);
		\coordinate (end) at (\scale,0);
		
		\begin{scope}[nodes={fill=black, circle, scale=.5}]
			\clip (-\scale-.1,0) rectangle (\scale+.1,\scale+.1);
			\draw[thick] (start)+(end) circle (\scale);
		\end{scope}
		
		\draw[thick] (start) -- (end);
		
		\foreach \arc in {0,1,...,\n}{			
			\coordinate (L\arc) at (180-\arc*180/\n : 1.2*\scale);
			
			\coordinate (N\arc) at (180-\arc*180/\n : \scale);
		}
		
		
		\foreach \num in {0,1,...,\n}{
			\node[fill=black, circle, scale=.5] at (N\num) {};
		}
		
		\foreach \num in {2,4}{
			\node[red] at (L\num) {$x_0$};
		}
		\foreach \num in {1,6,8}{
			\node[blue] at (L\num) {$x_0$};
		}
		
		\node[red] at (L3) {$x_1$};	
		
		\foreach \num in {5,7,9}{
			\node[blue] at (L\num) {$x_1$};
		}
		
		\node at (L0) {$x_1$};
		\node at (L\n) {$x_0$};

        \draw[thick, bend right = 15] (N0) to (N1);
		\draw[thick] (N1) to (N5);
        \draw[thick, bend right = 15] (N5) to (N6)
            to (N7)
            to (N8)
            to (N9)
            to (N10);
	\end{tikzpicture}
\end{center}
    Since $\phi(x_0x_0x_1) = 2$ and $\phi(x_0x_1x_0) = -1$ (cf. Definition \ref{def:B_1_to_c_w}) we obtain that
    \begin{equation*}
        (\phi\otimes \operatorname{id}) \circ \partial_{3}^{(1)}\big(\bzd(3,2,2,2)\big) = 3 \otimes \bzd(2,2,2)
    \end{equation*}
    which we identify with $3\,\bzd(2,2,2)$. Similarly, one computes that
    \begin{align*}
    	\partial_{5}^{(1)}\big(\bzd(3,2,2,2)\big) = &(x_0x_0x_1x_0x_1 - 2\, (x_0x_1)^2x_0) \otimes \bzd(2,2), \\
    	\partial_{7}^{(1)}\big(\bzd(3,2,2,2)\big) = &(x_0x_0x_1(x_0x_1)^2 +2\, (x_0x_1)^3x_0) \otimes \bzd(2), \\
    	\partial_{9}^{(1)}\big(\bzd(3,2,2,2)\big) = &x_0x_0x_1(x_0x_1)^3 \otimes \one. 
    \end{align*}
    Thus we obtain
    \begin{equation*}
        \partialphi[9][1](\bzd(3,2,2,2)) = 3\, \bzd(2,2,2) - \frac{15}{2}\, \bzd(2,2) + \frac{189}{16}\, \bzd(2) - \frac{223}{16}.
    \end{equation*}
    Similar computations for the remaining words in $\BNgr[9][1]$ give
    \begin{equation*}
        \MB[9][1] = 
        \begin{pmatrix}
            3 & -\frac{15}{2} & \frac{189}{16} & -\frac{223}{16} \\[1ex]
    		0 & -\frac{15}{2} & \frac{299}{8} & -\frac{889}{16} \\[1ex]
    		0 & 2 & -\frac{291}{16} & \frac{455}{16} \\[1ex]
    		-2 & 12 & -30 & \frac{641}{16} 
    	\end{pmatrix}.
    \end{equation*}
\end{Example}

\begin{Example}
    In our second example, we compute the matrix $\MB[10][2]$. We have
    \begin{align*}
        \BNgr[10][2] &= \{\bzd(3,3,2,2), \bzd(3,2,3,2), \bzd(3,2,2,3), \bzd(2,3,3,2), \bzd(2,3,2,3), \bzd(2,2,3,3)\}
    \shortintertext{and}
        \BNfil[10][1] &= \{\bzd(3,2,2), \bzd(2,3,2), \bzd(2,2,3), \bzd(3,2), \bzd(2,3), \bzd(3)\}.
    \end{align*}
    Similar to the previous example, one computes that
    \begin{align*}
		\partial_3^{(2)}\big(\bzd(3,2,2,3)\big) &= (x_0x_1x_0 - x_0x_0x_1 + x_0x_0x_1)\otimes \bzd(3,2,2) \\
		&+ (x_0x_0x_1 - x_0x_1x_0) \otimes \bzd(2,2,3), \\
		\partial_5^{(2)}\big(\bzd(3,2,2,3)\big) &= (x_0x_0x_1x_0x_1 - (x_0x_1)^2 x_0) \otimes \bzd(2,3) \\
		&+ ((x_0x_1)^2 x_0 - x_0x_0x_1x_0x_1 + x_0x_1x_0x_0x_1) \otimes \bzd(3,2), \\
		\partial_7^{(2)}\big(\bzd(3,2,2,3)\big) &= (x_0x_0x_1(x_0x_1)^2 - x_0x_0x_1(x_0x_1)^2 + (x_0x_1)^2 x_0x_0x_1) \otimes \bzd(3), \\
		\partial_9^{(2)}\big(\bzd(3,2,2,3)\big) &= 0.
	\end{align*}
	 By applying $\phi$ we obtain 
	 \begin{equation*}
	     \partialphi[10][2]\big(\bzd(3,2,2,3)\big) = -2\, \bzd(3,2,2) + 3\, \bzd(2,2,3) + 12\, \bzd(3,2) - \frac{15}{2}\, \bzd(2,3) -\frac{291}{16}\, \bzd(3).
	 \end{equation*}
	 Similar computations for the remaining elements of $\BNgr[10][2]$ yield
	 \begin{equation*}
	     \MB[10][2] = \begin{pmatrix}
	         3 & 0 & 0 & -12 & 0 & 28 \\[1ex]
		     0 & 3 & 0 & -\frac{11}{2} & 0 & 0 \\[1ex]
		     -2 & 0 & 3 & 12 & -\frac{15}{2} & -\frac{291}{16} \\[1ex]
		     0 & 0 & 0 & \frac{9}{2} & -10 & 0 \\[1ex]
		     0 & -2 & 0 & 0 & \frac92 & \frac{75}{8} \\[1ex]
		     0 & 0 & -2 & 0 & 12 & -\frac{291}{16} \\[1ex]
	     \end{pmatrix}.
	 \end{equation*}
\end{Example}

We conclude the examples with the observations that both $\MB[9][1]$ and $\MB[10][2]$ are invertible since 
\begin{align*}
    \det(\MB[9][1]) &= \frac{4865}{512}\\
    \det(\MB[10][2]) &= -\frac{435419}{64}
\end{align*}
and that all entries below the main diagonal are in $2\mathbb{Z}$, respectively.

\subsection{Some properties of the numbers \texorpdfstring{$c_{a,b}^r$}{C}}\label{subsec:numbers_c}

In this subsection we present some of the arithmetic properties of the numbers $c_{a,b}^r$ introduced in Definition \ref{def:B_1_to_c_w}. Later in Section \ref{subsec:lifting_zagier_to_fMZV}
these numbers occur in Zagier's Theorem and will play a crucial role in the main result of these notes. 

\begin{Definition} \label{def:coefficients_zagiers_formula}
	For integers $a,b,r\geq 0$ we set
	\begin{align*}
		A_{a,b}^r &= (1-2^{-2r}) \binom{2r}{2a+1}, \\
		B_{a,b}^r &= \binom{2r}{2b+2}.
	\end{align*}
Then, we get immediately
\begin{equation*}
	c_{a,b}^r = 2\cdot (-1)^r (B_{a,b}^r - A_{a,b}^r).
\end{equation*}
\end{Definition}

\begin{Definition}\label{def:p-adic_valuation}
	Let $p$ be a prime number and $q\in\QQ^\times$ non-zero. The \emph{$p$-adic valuation} $\nu_p(q)$ of $q$ is the integer $n\in\mathbb{Z}$ such that $q = p^n \frac{a}{b}$ with $a,b$ relatively prime to $p$. We further set $\nu_p(0) = \infty$.
\end{Definition}

It is not hard to verify that $p$-adic valuations satisfy the following basic properties for $q_1, q_2\in\Q$ 
\begin{align}
	\label{eq:p-adic_value_multiplicative}
	\nu_p(q_1 \cdot q_2) &= \nu_p(q_1) + \nu_p(q_2) \\
	\label{eq:p-adic_value_triangle-ineq}
	\nu_p(q_1 + q_2) &\geq \min\{\nu_p(q_1), \nu_p(q_2)\}
\end{align}
with equality in \eqref{eq:p-adic_value_triangle-ineq} if $\nu_p(q_1) \neq \nu_p(q_2)$.

\begin{Lemma}\label{lem:properties_c_beta}
	For all integers $a,b\geq 0$ we define 
	\begin{equation*}
		c_{a,b} = c_{a,b}^{a+b+1}.
	\end{equation*}
	These numbers satisfy
	\begin{enumerate}
		\item $c_{a,b}\in\mathbb{Z}[\frac12]$,
		\item $c_{a,b} - c_{b,a}\in 2\mathbb{Z}$ and
		\item $\nu_2(c_{a+b,0}) = \nu_2(c_{0,a+b}) \leq \nu_2(c_{a,b}) \leq 0$.
	\end{enumerate}
\end{Lemma}

\begin{proof}
	Write $n = a+b+1$.
	The first claim follows immediately from the definition
	\begin{equation*}
		c_{a,b} = 2\cdot (-1)^{n} \left(\binom{2n}{2b+2} - (1-2^{-2n}) \binom{2n}{2a+1}\right).
	\end{equation*}
	The second claim follows from \eqref{eq:final_lemma_claim_A}, i.\,e. $A_{a,b}^{n}-A_{b,a}^{n}=0$ hence
	\begin{equation*}
		c_{a,b} - c_{b,a} = 2\cdot (-1)^{n} (B_{a,b}^{n} - B_{b,a}^{n})\in 2\mathbb{Z}.
	\end{equation*}
	It remains to prove the third claim. We first show that $\nu_2((2n)!) < 2n$. By \eqref{eq:p-adic_value_multiplicative} we have 
	\begin{equation*}
		\nu_2((2n)!) = \sum_{i=1}^{2n} \nu_2(i).
	\end{equation*}
	We rearrange this sum and obtain
	\begin{align*}
		\sum_{i=1}^{2n} \nu_2(i) &= \sum_{k=1}^{\lfloor\log_2(2n)\rfloor} \#\big\{i\in\{1,\dots,2n\} \mid \text{$i$ is divisible by $2^k$}\big\}\\
		&= n + \left\lfloor\frac{n}{2}\right\rfloor + \left\lfloor\frac{n}{4}\right\rfloor + \dots + 1 < 2n.
	\end{align*}
	It follows that
	\begin{equation*}
		\nu_2\left(2^{-2n} \cdot \binom{2n}{2a+1}\right) < 0.
	\end{equation*}
	On the other hand, we clearly have $\nu_2\big(\binom{2n}{2b+2} - \binom{2n}{2a+1}\big) \geq 0$. 
	In particular, \eqref{eq:p-adic_value_triangle-ineq} becomes an equality for $q_1 = 2^{-2n} \cdot \binom{2n}{2a+1}$ and $q_2 = \binom{2n}{2b+2} - \binom{2n}{2a+1}$, hence
	\begin{align*}
		\nu_2(c_{a,b}) &= 1 + \nu_2\left( 2^{-2n} \cdot \binom{2n}{2a+1} + \binom{2n}{2b+2} - \binom{2n}{2a+1}\right) \\        
		&= 1 + \nu_2\left(2^{-2n} \cdot \binom{2n}{2a+1}\right) \leq 0.
	\end{align*}
	This proves the last inequality in the third claim.
	By writing $\binom{2n}{2a+1} = \frac{2n}{2a+1} \binom{2n-1}{2a}$ we further obtain from \eqref{eq:p-adic_value_multiplicative} that
	\begin{equation*}
		\nu_2(c_{a,b}) = 1 + \nu_2\left(2^{-2n} \cdot \frac{2n}{2a+1} \binom{2n-1}{2a}\right) = 2 - 2n + \nu_2(n) + \nu_2\left(\binom{2n-1}{2a}\right).
	\end{equation*}
	For $n\in\N$ fixed, this is minimal for $a = 0$ and $a = n-1$. The latter case is equivalent to $b = 0$.
\end{proof}

Later in Subsection \ref{subsec:lifting_zagier_to_fMZV} we need the following lemma.

\begin{Lemma} \label{lem:binomials}
	For $a,b\geq 0$ and $1\leq r\leq a+b+1$ we have
	\begin{equation*}
		c_{a,b}^r = \sum_{\substack{0\leq \alpha\leq a\\ 0\leq \beta\leq b\\ \alpha+\beta+1=r}} c_{\alpha,\beta}^r - \sum_{\substack{0\leq \alpha < a\\ 0\leq \beta\leq b\\ \alpha+\beta+1=r}}  c_{\beta,\alpha}^r + 2\cdot (-1)^r \big(\I(a\geq r) - \I(b\geq r)\big).
	\end{equation*}
	Here, $\I$ denotes the indicator function, i.e. $\I(A)=1$ if $A$ is a true statement and $\I(A) = 0$ else.
\end{Lemma}

\begin{proof}
	Using the quantities $A_{a,b}^r$ and $B_{a,b}^r$, it suffices to show that
	\begin{equation*}
		B_{a,b}^r - A_{a,b}^r = \sum_{\substack{0\leq \alpha\leq a\\ 0\leq \beta\leq b\\ \alpha+\beta+1=r}} (B_{\alpha,\beta}^r - A_{\alpha,\beta}^r) - \sum_{\substack{0\leq \alpha < a\\ 0\leq \beta\leq b\\ \alpha+\beta+1=r}}  (B_{\beta,\alpha}^r - A_{\beta,\alpha}^r) + (\I(a\geq r) - \I(b\geq r))
	\end{equation*}
	for $a,b\geq 0$ and $1\leq r\leq a+b+1$.

	We first show for all $\alpha,\beta\geq 0$ and $r = \alpha+\beta+1$ that
	\begin{align}
		\label{eq:final_lemma_claim_A}
		A_{\alpha,\beta}^r &= A_{\beta,\alpha}^r\\
		\label{eq:final_lemma_claim_B}
		B_{\alpha,\beta}^r &= B_{\beta+1,\alpha-1}^r.
	\end{align}
	Recall that we have for $0\leq k\leq n$
	\begin{equation*}
		\binom{n}{k} = \binom{n}{n-k}.
	\end{equation*}
	Identity \eqref{eq:final_lemma_claim_A} follows immediately for $n = 2r$ and $k = 2\alpha+1$. Identity \eqref{eq:final_lemma_claim_B} follows for $n = 2r$ and $k = 2\beta+2$ from 
	$$2(\alpha+\beta+1) - (2\beta +2) = 2\alpha = 2(\alpha-1) + 2.$$
	
	It suffices to show that
	\begin{align}
		\label{eq:final_lemma_claim_A2}
		A_{a,b}^r &= \sum_{\substack{0\leq \alpha\leq a\\ 0\leq \beta\leq b\\ \alpha+\beta+1=r}} A_{\alpha,\beta}^r - \sum_{\substack{0\leq \alpha < a\\ 0\leq \beta\leq b\\ \alpha+\beta+1=r}} A_{\beta,\alpha}^r \\
		\label{eq:final_lemma_claim_B2}
		B_{a,b}^r &= \sum_{\substack{0\leq \alpha\leq a\\ 0\leq \beta\leq b\\ \alpha+\beta+1=r}} B_{\alpha,\beta}^r - \sum_{\substack{0\leq \alpha < a\\ 0\leq \beta\leq b\\ \alpha+\beta+1=r}} B_{\beta,\alpha}^r + \I(a\geq r) - \I(b\geq r).
	\end{align}
	We first observe that $a\geq r$ is equivalent to 
	\begin{enumerate}
		\item $A_{a,\beta}^r$ vanishes for all $\beta\geq 0$,
		\item the summand for $\alpha = a$ does not appear in the first sum since $\alpha + \beta +1 = r$ would imply $\beta < 0$ and
		\item the summand for $\beta = 0$ appears in the second sum since then $\alpha = r-1 < a$.
	\end{enumerate}
	We deduce \eqref{eq:final_lemma_claim_A2} from the first two observations since
	\begin{equation*}
		\sum_{\substack{0\leq \alpha\leq a\\ 0\leq \beta\leq b\\ \alpha+\beta+1=r}} A_{\alpha,\beta}^r - \sum_{\substack{0\leq \alpha < a\\ 0\leq \beta\leq b\\ \alpha+\beta+1=r}} A_{\beta,\alpha}^r = A_{a,b}^r + \sum_{\substack{0\leq \alpha < a\\ 0\leq \beta\leq b\\ \alpha+\beta+1=r}} (A_{\alpha,\beta}^r - A_{\beta,\alpha}^r)
		\overset{\eqref{eq:final_lemma_claim_A}}{=} A_{a,b}^r.
	\end{equation*}
	
	We similarly observe for identity \eqref{eq:final_lemma_claim_B2} that $b\geq r$ is equivalent to
	\begin{enumerate}
		\item $B_{\alpha,b}^r$ vanishes for all $\alpha\geq 0$,
		\item the summand for $\beta = b$ does not appear in the first sum since $\alpha + \beta + 1 = r$ would imply $\alpha<0$ and
		\item the summand for $\alpha = 0$ appears in the first sum since then $\beta = r-1 < b$.
	\end{enumerate}
	So we have	
	\begin{align*}
		&\sum_{\substack{0\leq \alpha\leq a\\ 0\leq \beta\leq b\\ \alpha+\beta+1=r}} B_{\alpha,\beta}^r - \sum_{\substack{0\leq \alpha < a\\ 0\leq \beta\leq b\\ \alpha+\beta+1=r}} B_{\beta,\alpha}^r + \I(a\geq r) - \I(b\geq r) \\
		&\overset{\eqref{eq:final_lemma_claim_B}} {=} 
		B_{a,b}^r + \sum_{\substack{0\leq \alpha\leq a\\ 0\leq \beta < b\\ \alpha+\beta+1=r}} B_{\alpha,\beta}^r - \sum_{\substack{0\leq \alpha < a\\ 0\leq \beta\leq b\\ \alpha+\beta+1=r}} B_{\alpha+1,\beta-1}^r + \I(a\geq r) - \I(b\geq r)\\
		&= B_{a,b}^r + \sum_{\substack{0<\alpha\leq a\\ 0\leq \beta < b\\ \alpha+\beta+1=r}} B_{\alpha,\beta}^r + \I(b\geq r) - \sum_{\substack{0\leq \alpha < a\\ 0<\beta\leq b\\ \alpha+\beta+1=r}} B_{\alpha+1,\beta-1}^r - \I(a\geq r) \\
		&\quad + \I(a\geq r) - \I(b\geq r)\\
		\\
		&= B_{a,b}^r + \sum_{\substack{0<\alpha\leq a\\ 0\leq \beta < b\\ \alpha+\beta+1=r}} \Big( B_{\alpha,\beta}^r - B_{\alpha,\beta}^r \Big) \\
		&= B_{a,b}^r,
	\end{align*}
	where the second equality follows from the third observations above for $a\geq r$ and $b\geq r$, respectively, and the third equality follows from a shift of indices.
\end{proof}

\section{Formal multiple zeta values} 

\subsection{Multiple zeta values and their general structure}

We provide a short basic introduction in the theory of multiple zeta values, which will serve as the motivation for formal multiple zeta values. For a detailed exposition we refer to \cite{ba_notes}, \cite{bgf},  \cite{ikz}.

\begin{Definition}  \label{def:MZV}
	To integers $k_1\geq 2,\ k_2,\ldots,k_d\geq 1$, associate the \emph{multiple zeta value}
	\[\zeta(k_1,\ldots,k_d)=\sum_{n_1>\dots>n_d>0} \frac{1}{n_1^{k_1}\dots n_d^{k_d}} \in \mathbb{R}.\]
 	Denote the $\QQ$-vector space spanned by all multiple zeta values by
	\[\Z=\operatorname{span}_\QQ\{\zeta(k_1,\ldots,k_d)\mid d\geq 0,\ k_1\geq2,\ k_2,\ldots,k_d\geq 1\},\]
	where $\zeta(\emptyset)=1$.
	For an index $(k_1,\ldots,k_d)\in \NN^d$, define the \emph{weight} and \emph{depth} by
	\begin{align*}
		\wt(k_1,\ldots,k_d)=k_1+\dots+k_d, \qquad \dep(k_1,\ldots,k_d)=d.
	\end{align*}
	For simplicity, we will also refer to these numbers as the weight and depth of $\zeta(k_1,\ldots,k_d)$. 
\end{Definition} 

Numerical experiments have led to the following dimension conjectures for $\Z$.

\begin{Conjecture} \label{Zagier dimension conjecture} (\cite[p. 509]{zag}) \\
	1) The vector space $\Z$ is graded with respect to the weight, i.e.,  \[\Z=\bigoplus\limits_{w\geq0} \Zw{w},\] where $\Zw{w}$ is spanned by all multiple zeta values of weight $w$.
	\vspace{0,2 cm}\\ 
	2) The dimensions of the homogeneous subspaces $\Zw{w}$ are given by
	\[H_{\Z}(x)=\sum_{w\geq 0} \dim_\QQ(\Zw{w})x^w=\frac{1}{1-x^2-x^3}.
	\]
\end{Conjecture}

This conjecture implies that the numbers $d_w= \dim_\QQ(\Zw{w})$ should satisfy the recursion
\begin{align*}
	d_{w}=d_{w-2} +d_{w-3}    
\end{align*}
with initial values $d_1=0$ and $d_2=d_3=1$.

It is well-known that $\Z$ is not graded with respect to the depth, e.g., there is Euler's relation 
\begin{align} \label{Euler relation} 
	\zeta(2,1)=\zeta(3).
\end{align}
In these lectures we just want to mention the Broadhurst-Kreimer conjecture (\cite[(7)]{brkr}), which is a refinement of Zagier's dimension conjecture 
that in addition relies also on the depth filtration.

There exists also a suggestion for an explicit basis for $\Z$. We set
\begin{align*} \Ho=\operatorname{span}_{\QQ}\{\zeta(k_1,\ldots,k_d)\mid k_i\in \{2,3\} \text{ for } i=1,\ldots,d\}\subset \Z.
\end{align*} 

\begin{Conjecture} (\cite[Conjecture C]{h2}) A basis for $\Z$ is given by the Hoffman elements $\zeta(k_1,\ldots,k_d)$, $k_i\in \{2,3\}$. In particular, we have
	\[\Ho=\Z.\]
\end{Conjecture}

This conjecture would imply Zagier's dimension conjecture \ref{Zagier dimension conjecture}.

It is quite obvious, that there exists more structure on the space $\Z$.
\begin{Proposition} The space $\Z$ equipped with the usual multiplication of real numbers is an algebra. \qed
\end{Proposition}	
\label{stuffle shuffle} 

There are two ways of expressing the product of multiple zeta values, called the \emph{stuffle} and the \emph{shuffle product}. The stuffle product comes from the combinatorics of multiplying infinite nested sums. E.g., for $k_1,k_2\geq 2$, there is the simple calculation
\begin{align*} \zeta(k_1)\zeta(k_2)&=\left(\sum_{m>0}\frac{1}{m^{k_1}}\right)\left(\sum_{n>0}\frac{1}{n^{k_2}}\right)=\left(\sum_{m>n>0}+\sum_{n>m>0}+\sum_{m=n>0}\right) \frac{1}{m^{k_1}n^{k_2}}\\
	&= \zeta(k_1,k_2)+\zeta(k_2,k_1)+\zeta(k_1+k_2).
\end{align*}
The shuffle product is obtained from expressing multiple zeta values as iterated integrals (\cite[Theorem 1.108.]{bgf}). E.g., in depth $2$ the shuffle product reads for $k_1,k_2\geq 2$
\[\zeta(k_1)\zeta(k_2)=\sum_{j=2}^{k_1+k_2-1}\left(\binom{j-1}{k_1-1}+\binom{j-1}{k_2-1}\right)\zeta(j,k_1+k_2-j).\]

\subsection{Extended double shuffle relations}

To describe these two product expressions of multiple zeta values in general, we will use Hoffman's quasi-shuffle Hopf algebras (Subsection \ref{subsec:quasi-shuffle}). The comparison of these two product formulas lead us then to the (extended) double shuffle relations.

The shuffle product of the multiple zeta values can be described in terms of the previously studied shuffle Hopf algebra $(\QX,\shuffle,\dec)$. \\
Denote by $\h^0$ the subspace of $\QX$ generated by $\one$ and all words starting in $x_0$ and ending in $x_1$, so
\[\h^0=\QQ\one+x_0\QX x_1.\]

\begin{Theorem}  \label{shuffle hom} The map 
	\begin{align*}
		\zeta:(\h^0,\shuffle)&\to (\Z,\cdot), \\ 
		x_0^{k_1-1}x_1\cdots x_0^{k_d-1}x_1&\mapsto\zeta(k_1,\ldots,k_d)
	\end{align*}
	is a surjective algebra morphism compatible with notions of weight and depth for words and indices. 
	\qed
\end{Theorem} 

Recall that the \emph{weight} of a word $w\in \QX$ is the number of its letters, and by the \emph{depth} of a word $w\in \QX$ we mean the number of the letter $x_1$ in $w$.

To describe the stuffle product of multiple zeta values, we introduce a new alphabet.

\begin{Definition} \label{def stuffle} Consider the infinite alphabet $\Y=\{y_1,y_2,\ldots\}$. For a word in $\QY$, define the \emph{weight} and \emph{depth} by
	\begin{align*}
		\wt(y_{k_1}\cdots y_{k_d})=k_1+\dots+k_d,\qquad \dep(y_{k_1}\cdots y_{k_d})=d.
	\end{align*}
	Let the \emph{stuffle product} $\ast$ on $\QY$ be the quasi-shuffle product corresponding to 
	\[y_i\diamond y_j=y_{i+j} \qquad \text{ for } i,j\geq1.\] 
\end{Definition}

From Theorem \ref{dual to quasi-shuffle} we obtain the following.

\begin{Proposition}
	The tuple $(\QY,\ast,\dec)$ is a weight-graded commutative Hopf algebra. The complete dual Hopf algebra with respect to the pairing in \eqref{pairing quasi-shuffle} is given by $(\kYc{R}, \conc, \stco)$, where the coproduct $\stco$ is defined on the generators by
	\[\stco(y_i)=y_i\otimes \one+\one\otimes y_i+\sum_{j=1}^{i-1} y_j\otimes y_{i-j},\quad i=1,2,\ldots.\]
\end{Proposition} 

Denote by $\QY^0$ the subspace of $\QY$ generated by all words, which do not start in $y_1$.

\begin{Theorem} \label{stuffle hom} The map
	\begin{align*}
		\zeta:(\QY^0,\ast)&\to(\Z,\cdot), \\
		y_{k_1}\cdots y_{k_d}&\mapsto\zeta(k_1,\ldots,k_d)
	\end{align*}
	is a surjective algebra morphism compatible with the notions of weight and depth for words and indices. 
	\qed
\end{Theorem}

Comparing the shuffle and stuffle product formulas for multiple zeta values (Theorem \ref{shuffle hom}, \ref{stuffle hom}) gives the \emph{(finite) double shuffle relations} among multiple zeta values. Euler's relation given in \eqref{Euler relation} is not covered by the finite double shuffle relations, since there is no product decomposition in weight $3$. To get these kind of relations we will introduce regularizations.

\begin{Proposition} \label{stuffle reg} Let $T$ be a commutative variable and extend the stuffle product $\ast$ by $\QQ[T]$-linearity to $\QY^0[T]$. We have an algebra isomorphism
	\begin{align*}
		\operatorname{reg}_\ast:\QY^0[T]\to \QY, \\
		wT^n\mapsto w\ast y_1^{\ast n}.
	\end{align*}
	\vspace{-1cm} \\ \qed 
\end{Proposition} 

For any $w\in\QY$, set
\begin{align*}
	\zeta_\ast^T(w)=\zeta(\operatorname{reg}_\ast^{-1}(w))\in \Z[T],
\end{align*}
where we extend also the map $\zeta:\QY^0\to\Z$ by $\QQ[T]$-linearity to $\QY^0[T]$. We call
\begin{align*}
	\zeta_\ast(w)=\zeta_\ast^{T=0}(w)=\zeta(\operatorname{reg}_\ast^{-1}(w)|_{T=0})\in \Z
\end{align*}
the \emph{stuffle-regularized multiple zeta values}.
An immediate consequence of Proposition \ref{stuffle reg} is the following.

\begin{Theorem} \label{stuffle reg properties} The map $\zeta_\ast:\QY\to \Z$ given by $w\mapsto \zeta_\ast(w)$ is the unique map satisfying
	\begin{itemize}
		\item[(i)] $\zeta_\ast(w)=\zeta(w)$ for all $w\in \QY^0$,
		\item[(ii)] $\zeta_\ast(y_1)=0$,
		\item[(iii)] $\zeta_\ast(u)\zeta_\ast(v)=\zeta_\ast(u\ast v)$ for all $u,v\in \QY$.
	\end{itemize}
\end{Theorem} 

Similarly, there also exists a regularization with respect to the shuffle product for multiple zeta values.

\begin{Proposition} \label{shuffle reg} Let $T,U$ be a commutative variables and extend the shuffle product $\shuffle$ by $\QQ[T,U]$-linearity to $\h^0[T,U]$. There is an algebra isomorphism
	\begin{align*}
		\operatorname{reg}_{\shuffle}:\h^0[T,U]&\to \QX, \\
		wT^nU^m &\mapsto w\shuffle x_1^{\shuffle n}\shuffle x_0^{\shuffle m}.
	\end{align*}
	\vspace{-1cm} \\ \qed 
\end{Proposition} 

For any $w\in\QX$, set
\begin{align*}
	\zeta_{\shuffle}^{T,U}(w)=\zeta(\operatorname{reg}_{\shuffle}^{-1}(w))\in \Z[T,U],
\end{align*}
where the map $\zeta:\h^0\to\Z$ needs to be extended by $\QQ[T,U]$-linearity to $\h^0[T,U]$. We set
\begin{align*}
	\zeta_{\shuffle}^T(w)=\zeta(\operatorname{reg}_{\shuffle}^{-1}(w)|_{U=0})\in \Z[T],
\end{align*}
and, moreover, call
\begin{align*}
	\zeta_{\shuffle}(w)=\zeta_{\shuffle}^{T=U=0}(w)=\zeta(\operatorname{reg}_\ast^{-1}(w)|_{T=U=0})\in\Z
\end{align*}
the \emph{shuffle-regularized multiple zeta values}.
One derives the following from Proposition \ref{shuffle reg}.

\begin{Theorem} \label{shuffle reg properties} The map $\zeta_{\shuffle}:\QX\to \Z,\ w\mapsto \zeta_{\shuffle}(w)$ is the unique map satisfying
	\begin{itemize}
		\item[(i)] $\zeta_{\shuffle}(w)=\zeta(w)$ for all $w\in \h^0$,
		\item[(ii)] $\zeta_{\shuffle}(x_1)=\zeta_{\shuffle}(x_0)=0$,
		\item[(iii)] $\zeta_{\shuffle}(u)\zeta_{\shuffle}(v)=\zeta_{\shuffle}(u\shuffle v)$ for all $u,v\in \QX$.
	\end{itemize}
\end{Theorem} 

We want to relate these two regularizations of multiple zeta values. Define the $\Z$-linear map $\rho:\Z[T]\to\Z[T]$ by 
\begin{align} \label{map rho}
	\rho\left(\frac{T^m}{m!}\right)=\sum_{i=0}^m\gamma_i\frac{T^{m-i}}{(m-i)!},\qquad m=0,1,2,\ldots,
\end{align}
where the coefficients $\gamma_i\in \Z$ are defined by $\sum\limits_{i\geq0} \gamma_i u^i=\exp\left(\sum\limits_{n\geq 2} \frac{(-1)^n}{n}\zeta(n)u^n\right)$.

\begin{Theorem} \label{Relate stuffle- and shuffle-regularized} (\cite[Theorem 1]{ikz})
	For all $k_1,\ldots,k_d\geq 1$, one has \[\rho\big(\zeta_\ast^T(y_{k_1}\cdots y_{k_d})\big)=\zeta_\shuffle^T(x_0^{k_1-1}x_1\cdots x_0^{k_d-1}x_1).\] \vspace{-1cm} \\ \qed 
\end{Theorem} 

Combining the stuffle product formula for the stuffle-regularized multiple zeta values and the shuffle product formula for the shuffle-regularized multiple zeta values together with Theorem \ref{Relate stuffle- and shuffle-regularized} gives the \emph{extended double shuffle relations} among multiple zeta values.

\begin{Conjecture} \label{Conj all relations in Z} (\cite[Conjecture 1]{ikz}) All algebraic relations in the algebra $\Z$ of multiple zeta values are a consequence of the extended double shuffle relations.
\end{Conjecture} 

In particular, Conjecture \ref{Conj all relations in Z} would imply that the algebra $\Z$ is graded by weight, since the stuffle and the shuffle product are both homogeneous for the weight.

\begin{Example} \label{exm:stuffle_shuffle_prod}
	We calculate
	\begin{align*} 
		\zeta_\ast(y_1)\zeta_\ast(y_2)&=\zeta_\ast(y_1y_2)+\zeta_\ast(y_2y_1)+\zeta_\ast(y_3), \\
		\zeta_{\shuffle}(x_1)\zeta_{\shuffle}(x_0x_1)&=\zeta_{\shuffle}(x_1x_0x_1)+2\zeta_{\shuffle}(x_0x_1x_1)
	\end{align*}
	By Theorem \ref{Relate stuffle- and shuffle-regularized}, we have $\zeta_{\ast}(y_1y_2)=\zeta_{\shuffle}(x_1x_0x_1)$, and, moreover, $\zeta_\ast(y_1)=\zeta_{\shuffle}(x_1)=0$. Thus, the above equations reduce to
	\begin{align*}
		\zeta_\ast(y_2y_1)+\zeta_\ast(y_3)=2\zeta_{\shuffle}(x_0x_1x_1).
	\end{align*}
	By Theorem \ref{stuffle reg properties}, \ref{shuffle reg properties}, this is equivalent to 
	\begin{align*}
		\zeta(3)=\zeta(2,1).
	\end{align*}
	So we recover Euler's relations from the extended double shuffle relations.
\end{Example}

\begin{Remark} \label{rem:Ba_relations-poster} There are various kinds of relations obtained for multiple zeta values in the literature, and several of them are expected to give all algebraic relations in $\Z$. A graphical overview is given by H. Bachmann in \cite[Section 3.3]{ba_notes}.
\end{Remark}

\subsection{Formal multiple zeta values and some properties}

We have a canonical embedding
\begin{align*}
\iota:\QY &\hookrightarrow \QX, \\
y_{k_1}\cdots y_{k_d}&\mapsto x_0^{k_1-1}x_1\cdots x_0^{k_d-1}x_1,
\end{align*}  
which allows to transfer the stuffle product $\ast$ on $\QY$ (Definition \ref{def stuffle}) to the algebra $\h^1=\Q \one+\QX x_1$.
So motivated by Conjecture \ref{Conj all relations in Z}, we can reformulate the extended double shuffle relations purely algebraically in $\QX$.
\begin{Definition} \label{algebra Z^f} 
The algebra $\Zf$ of \emph{formal multiple zeta values} is given by
\[\Zf=\faktor{(\QX,\shuffle)}{\operatorname{R}_{\operatorname{EDS}}},\]
where $\operatorname{R}_{\operatorname{EDS}}$ is the ideal generated by the extended double shuffle relations. Following \cite{ikz}, $\operatorname{R}_{\operatorname{EDS}}$ is for example generated by $x_0$, $x_1$ and
\[ u\shuffle v - \iota\Big(\iota^{-1}(u)\ast \iota^{-1}(v)\Big),\qquad u \in \h^0, v\in \h^1.\]
\end{Definition} 

\begin{Remark}
Since the stuffle product $\ast$ as well as the shuffle product $\shuffle$ are graded for the weight, the ideal $\operatorname{R}_{\operatorname{EDS}}$ is weight-homogeneous. Therefore, the algebra $\Zf$ is weight-graded.
\end{Remark}

\begin{Remark}
We won't compare the formal multiple zeta values to the motivic multiple zeta values, which are formalisations that depend on a different set of relations. Their construction is requires a profound knowledge of modern arithmetic geometry, for detail we refer to \cite{bgf}. See also Remark  \ref{rem:Ba_relations-poster}.    
\end{Remark}   

Let us denote the canonical projection by
\begin{equation}\label{eq:def_zf}
 \zf: \QX \to \Zf,   
\end{equation} 
thus $\zf(w)$ is the class of $w\in \QX$ in the quotient algebra $\Zf$. The space $\Zf$ is a weight-graded algebra spanned by the elements $\zf(w)$, $w\in\X^*$, which satisfy exactly the extended double shuffle relations. In particular, we have for all $k \ge 1$ that
\begin{align}\label{eq:xxxvanishes}
 \zf(x_0^k)=\zf(x_1^k)=0.   
\end{align}

If we have a word $w=x_0^{k_1-1}x_1\cdots x_0^{k_d-1}x_1\in \h^0$, we often write instead of $\zf(w)$ also $\zf(k_1,\ldots,k_d)$. 
By Proposition \ref{shuffle reg} and \eqref{eq:xxxvanishes} these are also a spanning set for $\Zf$.
Additionally, the $\zf(k_1,\ldots,k_d)$ with $k_1\geq2$ also satisfy the usual stuffle product formulas.  

By construction,  the \emph{evaluation map} 
\begin{align} \label{evaluation map}
\Zf&\to\Z,\\
\zf(w)&\mapsto\zeta_{\shuffle}(w) \nonumber
\end{align}
is a surjective algebra morphism.
By Conjecture \ref{Conj all relations in Z}, this map should be an algebra isomorphism.

\begin{Proposition}\label{prop:zeta_222} 
For $k\geq0$, let $B_k$ be the $k$-th Bernoulli number and set 
\[b_{n} = (-1)^{n+1} B_{2n}\frac{24^n}{2 (2n)!}.\] Then we have for each $n\in\NN$
\begin{align*}
\zf( 2n )  & = b_{n} \zf(2)^n, \\
\zf(\{2\}^n)  &=   \frac{6^n}{(2n+1)!} \zf(2)^n.
\end{align*}
\end{Proposition}

Here we use the notation $\zf(\{2\}^n)=\zf(\underbrace{2,\ldots,2}_{n \text{ times}})$.

\begin{proof} 
Euler showed by analytic means  his famous  formula
\[
\zeta(2n) =  (-1)^{n+1} \frac{B_{2n}}{2(2n)!} (2\pi)^{2n}  = b_n \zeta(2)^n.
\]
For the formal multiple zeta values we deduce this formula from the extended double shuffle as follows: 
Considering the generating series for the formal multiple zeta values in depth two we deduce  as Gangl, Kaneko and Zagier in \cite{gkz} the identity  
\[
\sum_{\substack{ r,s \text{ even} \\ r+s =k } }   \zf( r, s ) = \frac{3}{4} \zf(k). 
\]
Now using  $
\zf( r) \zf(s) = \zf(r,s) +\zf(s,r) + \zf(r+s)
$
we get for the sum of products
\[
\sum_{\substack{ r,s \text{ even} \\ r+s =k } } \zf( r) \zf(s)  = \frac{k+1}{2} \zf(k).
\]
This yields a recursion for the numbers $b_n$ in the equation $\zf (2n)=b_n\zf(2)^n $, whose first terms are
\begin{align*}
\zf(4)& = \frac{2}{5} \zf(2)^2 ,\\
\zf(6) & = \frac{2}{7} \Big( \zf(2) \zf(4) + \zf(4) \zf(2) \Big) = \frac{8}{35} \zf(2)^3.
\end{align*}
 Now using a simple manipulation of the generating series for
 Bernoulli numbers\footnote{It is not difficult to check that $-\frac{1}{2} \left( \sqrt{6} x \operatorname{cot}(\sqrt{6} x) -1
  \right) $ is the generating series for the rational numbers $b_n$.
This approach to obtain this formula of Euler in the formal setting is elaborated in \cite[Corollary 2.12]{BaIt_fMES}. 
  }  or using the evaluation map and Euler's formula for the even zeta values 
 we deduce
 \[
 b_n= (-1)^{n+1} B_{2n}\frac{24^n}{2 (2n)!}.
 \]

In their paper \cite{hi} on quasi-shuffle algebras Hofmann and Ihara showed that for all $k\ge 2$ we have the identity of generating series
\begin{align*}
\exp\left( \sum_{i\ge 1}  \frac{(-1)^{i-1}}{i} \zf(i k ) t^i \right) = \sum_{n=0}^\infty \zf( \{k\}^n ) \,t^n .
\end{align*}     
If we specialize this to $k=2$ we see that $\zf( \{2\}^n )$ is a polynomial in even formal zeta values, however by the algebraic Euler formula this
equals a rational multiple of $\zf(2)^n$. As before we can determine this proportionality factor either by a direct calculation or by a comparison with the analytic 
identity 
\[
\zeta ( \{2\}^n ) = \frac{\pi^{2n}}{(2n+1)!}  = \frac{6^n}{(2n+1)!} \zeta(2)^n.
\]
due to Hoffman and Zagier.
\end{proof}

\begin{Remark} In the above proof we reduced identities for formal multiple zeta values to the determination of a rational number, which we computed by means of the evaluation map. Actually, we could even had calculated it without referring to the previous known analytic identities.
In Subsection \ref{subsec:lifting_zagier_to_fMZV} we see a much more elaborate lift of another set of identities satisfied by multiple zeta values to the formal setting. Those formulae due to Zagier can be seen as a refinement of the following proposition and so far no algebraic proof without referring to analytical identities is known.
\end{Remark}

\begin{Proposition}	\label{prop:sum_level_one}For all $n\in\N$ we have the identity
\begin{equation}\label{eq:palindromic_fmzv}
\begin{aligned}
    \zf( (x_0x_1)^n x_0) &= -2 \sum_{i=0}^{n-1} \zf(\{2\}^i,3,\{2\}^{n-1-i}) \\
    &= 2\sum_{i=1}^n (-1)^{i} \zf(2i+1)\, \zf(\{2\}^{n-i}).
\end{aligned}
\end{equation}
\end{Proposition}

\begin{proof}
For the first equality we consider the shuffle product
 identity in $\QX$
\begin{align*}\label{eq:B_one_in_H} 
    x_0 \shuffle (x_0x_1)^n = (x_0x_1)^n x_0 + \sum_{i=0}^{n-1} 2\, (x_0x_1)^i x_0x_0x_1 (x_0x_1)^{n-1-i} 
\end{align*}
and this amounts in the quotient algebra $\Zf$ to 
\begin{equation*}
    \zf(x_0) \zf(\{2\}^n) = \zf((x_0x_1)^n x_0) + \sum_{i=0}^{n-1} 2\, \zf(\{2\}^i,3,\{2\}^{n-1-i}).
\end{equation*}
The left-hand side vanishes since $\zf(x_0) = 0$ and the first identity follows.

For the second equality we use that the multiplication of formal zeta values with admissible indices satisfy the stuffle relation, we thus get
\begin{align*} 
    \zf(3) \zf(\{2\}^{n-1}) &=  \sum_{i=0}^{n-1}  \zf(\{2\}^i,3,\{2\}^{n-1-i}) 
    + \sum_{i=0}^{n-2}  \zf(\{2\}^{i},5,\{2\}^{n-2-i}) \\
    \zf(5) \zf(\{2\}^{n-2}) &=  \sum_{i=0}^{n-2}  \zf(\{2\}^{i},5,\{2\}^{n-2-i}) 
    + \sum_{i=0}^{n-3}  \zf(\{2\}^{i},7,\{2\}^{n-3-i}) \\
    &\;\vdots \\
    \zf(2n-1) \zf(2) &= \zf(2n-1,2) +\zf(2,2n-1) +\zf(2n+1).
\end{align*}
By taking the alternating sum we obtain 
\begin{equation*}
    -\sum_{i=0}^{n-1} \zf(\{2\}^i,3,\{2\}^{n-1-i}) = \sum_{i=1}^n (-1)^i \zf(2i+1) \zf(\{2\}^{n-i})
\end{equation*}
and the claim follows.
\end{proof}

\begin{Remark}\label{rem:zf_B1_level_1}
The identity \eqref{eq:palindromic_fmzv}  implies    
for the space $\mathcal{B}^1$ defined in Theorem \ref{thm:graded_image_D_B23} that 
$\zf\big(\mathcal{B}^1 \big)$ is contained in the vector space
generated by those $\zf(k_1,...,k_d)$, where one $k_i$ equals $3$ and all other $k_j$ equal $2$. 
\end{Remark}

\subsection{Racinet's approach and Ecalle's theorem} 

We explain Racinet's approach to formal multiple zeta values in terms of non-commutative power series (\cite{ra}). This means, we assign an affine group scheme and a Lie algebra to the formal multiple zeta values. This has two deep, structural consequences for the algebra $\Zf$. For more detailed expositions we refer to \cite{bu}, \cite{enr}.

\begin{Definition}
Let $\mathcal{A}$ be an alphabet, $R$ be some $\QQ$-algebra, and $\varphi:\nca{\QQ}{\A}\to R$ be a $\QQ$-linear map. Then the non-commutative generating series associated to $\varphi$ is
\[\operatorname{Gen}_{\varphi}(\A)=\sum_{w\in \A^*} \varphi(w)w\in \ncac{R}{\A}.\]
\end{Definition}

We assume that $(\nca{\QQ}{\A},\quasish)$ is any graded quasi-shuffle algebra with $\deg(a)\geq1$ for all $a\in \A$. Then we can compute the dual coproduct $\Delta_{\quasish}$ given in \eqref{quasish dual coproduct}.

\begin{Proposition} \label{gen series grouplike}
A $\QQ$-linear map $\varphi:\nca{\QQ}{\A}\to R$ is an algebra morphism with respect to $\quasish$ if and only if $\operatorname{Gen}_{\varphi}(\A)$ is grouplike for $\Delta_{\quasish}$,
\[\Delta_{\quasish}(\operatorname{Gen}_{\varphi}(\A))=\operatorname{Gen}_{\varphi}(\A) \otimes \operatorname{Gen}_{\varphi}(\A). \]
\vspace{-1cm} \\ \qed
\end{Proposition}

Motivated by this proposition, we consider the following sets.

\begin{Definition} \label{def DM} 
For any commutative $\QQ$-algebra $R$ with unit, let $\DM(R)$ be the set of all non-commutative power series $\phi\in \kXc{R}$ satisfying
\renewcommand{\arraystretch}{1,1}
\begin{center} \centering \begin{tabular}{cclc}
(i) & $(\phi|x_0)\ =\ (\phi|x_1)$ & $=$ & $0$, \\ 
(ii) & $\co(\phi)$ & $=$ & $\phi\hat{\otimes}\phi$, \\
(iii) & $\stco(\phi_\ast)$ & $=$ & $\phi_\ast\hat{\otimes}\phi_\ast$, 
\end{tabular} \end{center} \renewcommand{\arraystretch}{1}
where 
\[\phi_\ast=\exp\left(\sum_{n\geq 2}\frac{(-1)^{n-1}}{n}(\Pi_\Y(\phi)|y_n)y_1^n\right)\Pi_\Y(\phi)\in \kYc{R}\] and $\Pi_\Y$ is the $R$-linear extension of the projection 
\begin{align*} 
\QX&\to\QY, \\ 
x_0^{k_1-1}x_1\cdots x_0^{k_d-1}x_1x_0^{k_{d+1}-1}&\mapsto\begin{cases} y_{k_1}\cdots y_{k_d}, & \quad k_{d+1}=1 \\ 0 & \quad \textsf{else}\end{cases}.
\end{align*}
\vspace{0,2cm} \\
For each $\lambda\in R$, denote by $\DM_\lambda(R)$ the set of all $\phi\in \DM(R)$, which additionally satisfy
\[\hspace{-2,2cm}\text{(iv)} \qquad (\phi|x_0x_1) = \lambda. \]
\end{Definition} 

By Theorems \ref{stuffle reg properties}, \ref{shuffle reg properties}, \ref{Relate stuffle- and shuffle-regularized} and Proposition \ref{gen series grouplike}, the non-commutative generating series of the shuffle regularized multiple zeta values
\begin{align*}
\operatorname{Gen}_{\zeta_{\shuffle}}(\X)=\sum_{w\in \X^*} \zeta_{\shuffle}(w)w
\end{align*}
is an element in $\DM_{\pi^2/6}(\Z)$.

\begin{Theorem} \label{DM non-empty} (\cite[Theorem I]{ra2}) For each commutative $\QQ$-algebra $R$ and $\lambda\in R$, the set $\DM_\lambda(R)$ is non-empty. \qed
\end{Theorem} 

From \cite{dr} and \cite{fu} one deduces that there also exist elements $\phi$ in $\DM_\lambda(R)$ additionally satisfying
\begin{align*} \label{sol to DS equation even}
(\phi\mid x_0^{k}x_1)=0 \quad \text{ for } k\geq 1 \text{ even}.
\end{align*} 
The sets $\DM(R)$ and $\DM_\lambda(R)$ give rise to affine schemes represented by (quotient algebras of) $\Zf$.

\begin{Proposition} \label{DM represented by} (\cite[p. 107]{ra}) \\
(i) The functor $\DM_\lambda:\QAlg\to \Sets$ is an affine scheme represented by the algebra $\Zf$ of formal multiple zeta values. In particular, for $R\in\QAlg$ there is a bijection 
\begin{align*} 
\operatorname{Hom}_{\QAlg}(\Zf,R)&\to\DM(R), \\
\varphi&\mapsto\sum_{w\in \X^*} \varphi(\zf(w)) w.
\end{align*} 
(ii) The functors $\DM_\lambda:\QAlg\to \Sets$ are an affine schemes represented by the quotient algebras $\faktor{\Zf}{(\zf(2)-\lambda)}$. \qed \end{Proposition} 

$\DM_0$ is often called the \emph{double shuffle group}. To figure out the group structure for the affine scheme $\DM_0$, one needs to consider first the corresponding linearized space.

\begin{Definition} \label{def dm0} For any $R\in \QAlg$, let $\mathfrak{dm}(R)$ be the $R$-vector space of all non-commutative polynomials $\psi\in\kX{R}$ satisfying
\renewcommand{\arraystretch}{1,1}
\begin{center} \centering \begin{tabular}{cclc}
(i) & $(\psi\mid x_0)\ =\ (\psi\mid x_1)$ & $=$ & $0$, \\ 
(ii) & $\co(\psi)$ & $=$ & $\psi\otimes\one+\one\otimes\psi$, \\
(iii) & $\stco(\psi_\ast)$ & $=$ & $\psi_\ast\otimes\one+\one\otimes\psi_\ast$, \\
\end{tabular} \end{center} \renewcommand{\arraystretch}{1}
where 
\[\psi_\ast=\Pi_\Y(\psi)+\sum\limits_{n\geq 2}\frac{(-1)^{n-1}}{n}(\Pi_\Y(\psi)\mid y_n)y_1^n\in \kY{R}\] 
and $\Pi_\Y$ is the $R$-linear extension of the canonical projection $\QX\to\QY$.
\vspace{0,2cm} \\
By $\dm(R)$ denote the subspace of all $\psi\in\mathfrak{dm}(R)$ additionally satisfying
\[\hspace{-3,5cm}\text{(iv)} \qquad (\psi\mid x_0x_1)=0.\]
\end{Definition} 

Denote $\dm=\dm(\QQ)$. Then one has $\dm(R)=\dm\otimes R$. The space $\dm$ is often called the \emph{double shuffle Lie algebra}, the name will be justified in Theorem \ref{Racinet theorem}.

\begin{Example} \label{elements dm} The space $\dm$ can be computed algorithmically in small weights, see  e.g. \cite{enr},  \cite{bu_master}. One obtains the following elements in $\dm$ up to weight $5$
\begin{align*}
\xi(3)&=[x_0,[x_0,x_1]]+[[x_0,x_1],x_1],\\
\xi(5)&=[x_0,[x_0,[x_0,[x_0,x_1]]]]+2[[x_0,[x_0,[x_0,x_1]]],x_1]+\frac{1}{2}[[x_0,[x_0,x_1]],[x_0,x_1]]\\
&\hspace{0,5cm}+2[x_1,[x_1,[x_0,[x_0,x_1]]]]-\frac{3}{2}[[x_0,x_1],[[x_0,x_1],x_1]]+[[[[x_0,x_1],x_1],x_1],x_1].
\end{align*}
\end{Example}

\begin{Proposition} \label{dm0 even weigth depth 1} (\cite[IV, Proposition 2.2]{ra}) For each $k\geq 2$ even and $\psi\in \dm(R)$, one has
\[(\psi\mid x_0^{k-1}x_1)=0.\]
\vspace{-1cm} \\ \qed 
\end{Proposition} 

The Lie bracket on the spaces $\dm(R)$ is given by the Ihara bracket $\{-,-\}$ from Definition \ref{def:Ihara}.

\begin{Theorem} \label{Racinet theorem} (\cite[IV, Proposition 2.28., Corollary 3.13.]{ra})  Let $R$ be a commutative $\QQ$-algebra with unit. 
\vspace{0,3cm}\\
1) The space $\dm(R)$ is a Lie subalgebra of $(\fX\otimes R,\{-,-\})$.
\vspace{0,3cm} \\
2) For all $\phi_1,\phi_2\in \DM_\lambda(R)$, there exists a unique element $\psi$ in the completed Lie algebra $\widehat{\dm}(R)$ such that
\[\exp(s_\psi)(\phi_1)=\phi_2.\] \qed
\end{Theorem}

Here for any $f\in \kX{R}$, $s_f:\kX{R}\to \kX{R}$ is the $R$-linear map defined by $s_f(w)=d_f(w)+fw$. Note that this allows to write
\begin{align*}
\{f,g\}=s_f(g)-s_g(f), \qquad f,g \in \kX{R}.
\end{align*}

A more detailed proof of Theorem \ref{Racinet theorem} 1) is given in \cite[Appendix A]{fu} and also in the first author's master thesis \cite{bu_master}.

From Theorem \ref{Racinet theorem} 2), one obtains natural bijections
\begin{align} \label{bijection dm DM}
\widehat{\dm}(R)&\to \DM_0(R), \\
\psi&\mapsto \exp(s_\psi)(\one). \nonumber
\end{align}
We deduce the following.

\begin{Corollary} \label{DM0 ags} (\cite{ra}) The functor $\DM_0$ is a pro-unipotent affine group scheme with Lie algebra functor 
\[\widehat{\dm}:\QQ\operatorname{-Alg}\to\operatorname{Lie-Alg},\ R \mapsto\widehat{\dm}(R). \] 
\vspace{-1cm} \\ \qed \end{Corollary} 

The group multiplication on $\DM_0(R)$ is given in Theorem \ref{thm:GL_for_Ihara} or Remark \ref{rem:GLP_ihara_grouplike}, so $\DM_0(R)$ is a subgroup of $(\operatorname{Grp}(\widehat{\mathcal{U}}(\fX\otimes R)),\glp)$.

Moreover, it is shown by Racinet that the affine schemes $\DM_\lambda$ are $\DM_0$-torsors.

\begin{Theorem} \label{thm:dm_torsor} (\cite[Section IV, Corollary 3.13]{ra})
The affine group scheme $\DM_0$ acts freely and transitively on the affine schemes $\DM_\lambda$ by left multiplication.  \hfill \qed
\end{Theorem}

In other words, for each $R\in \QAlg$ we obtain a map
\begin{align*}
\glp: \DM_0(R) \times \DM_\lambda(R) &\to \DM_\lambda(R), \\
(\Phi,G)&\mapsto \Phi\glp G, \nonumber
\end{align*}
such that $-\glp G:\DM_0(R)\to \DM_\lambda(R)$ is bijective for any $G\in \DM_\lambda(R)$. Combining those maps gives rise to natural maps 
\begin{align} \label{eq:map_DM_DM0}
\glp:\DM_0(R)\times \DM(R)\to \DM(R).
\end{align}

An application of Yoneda's Lemma to the isomorphism $\exp: \widehat{\dm}\to \DM_0$ of affine schemes given in \eqref{bijection dm DM}, yields the following main result.
\begin{Corollary} \label{Z^f polynomial algebra} (\cite{ra}, Chapter IV, Corollary 3.14) There is an algebra isomorphism
\[ \Zf\simeq \QQ[\zf(2)]\otimes_\QQ\mathcal{U}(\dm)^\vee. \qedhere\]
\end{Corollary}
By \eqref{eq:iso_US}, we have proved Ecalle's theorem \cite{ecalle:tale}.

\begin{Theorem} The algebra $\Zf$ of formal multiple zeta values is a free polynomial algebra.
\end{Theorem}

\subsection{Goncharov coproduct and Goncharov-Brown coaction} \label{Gon coproduct Gon-Brown coaction}

Let 
\[\Af=\Zf/(\zf(2))\]
be the quotient algebra of $\Zf$ by the principal ideal generated by $\zf(2)$ and write
\[
    \amap: \Zf \to \Af
\]
for the natural projection. We define a coproduct on $\Af$ by
\begin{align}
 \label{eq:Gon_on_Zf}
\DeltaGon(a) = \Big( \amap \circ \zf  \otimes\,   \amap \circ \zf  \Big) \big( \DeltaGon( w) \big) , 
\end{align}
where $w \in \QX$ is a lift for $a \in \Af$, i.e. $w$ is any element such that $ \amap(\zf(w)) = a$.

\begin{Theorem} The triple $\big(\Af,\cdot,\DeltaGon \big)$ is a weight-graded Hopf algebra.
\end{Theorem}    
\begin{proof}
  In Corollary \ref{DM0 ags}, we have seen that $\DM_0$ is a subgroup scheme of the affine group scheme given in \eqref{eq: ags from grouplike} induced by the grouplike elements of $(\widehat{\mathcal{U}}(\fX),\co)$. Since $\DM_0$ is represented by the algebra $\Af$ (Proposition \ref{DM represented by}), and the affine group scheme in \eqref{eq: ags from grouplike} is represented by the Hopf algebra $(\QX,\shuffle,\DeltaGon)$ (Lemma \ref{lem:glp_gon_dual}, Proposition \ref{prop:QX_gon_hopf-alg}), also $\Af$ must be a Hopf algebra. The coproduct on $\Af$ is induced by the Goncharov coproduct under the projection
    \[\amap\circ\zf:\QX\to \Af. \qedhere\] 
\end{proof}
 
\begin{Remark} In particular we deduce from the proof of the above theorem that the coproduct  \eqref{eq:Gon_on_Zf} is well-defined. 
We like to point to the fact that our approach
to the coproduct relies on the general theory of post-Lie algebras together with the work of Racinet, whereas
the original definition of Goncharov in \cite{gon} 
was based on topological considerations for the path algebra. The latter is directly related to the representation of multiple zeta values by iterated integrals.
\end{Remark} 

Following Proposition \ref{indecomposables Lie cobracket}, the Goncharov coproduct $\DeltaGon$ induces a Lie cobracket $\delta$ on the space of indecomposables $\Lf$ of the quotient algebra $\Af$. 
By Proposition \ref{duality Lie algebra to ags and indecomposables}, there is a canonical isomorphism of Lie algebras
\[\dm\simeq (\Lf)^\vee.\]

Summarizing the previous results leads to the following diagram, which should be seen as a special case of the diagram \eqref{diagram big picture general}
\begin{equation}  \label{diagram big picture MZV}
\begin{tikzcd}[baseline=(current  bounding  box.center)]
\big(\Af,\cdot,\DeltaGon\big) \arrow[dddd] &&&&&& \big(\mathcal{U}(\dm),\glp,\co\big) \arrow[llllll,"\sim","\text{graded dual}" '] 
\\ \\ 
&&& \big(\DM_0,\glp\big) \arrow[uulll,leftrightarrow,"1:1", "\text{Prop \ref{DM represented by}}"'] &&&
\\ \\ 
\big(\Lf,\delta\big) &&&&&& \big(\dm,\{-,-\}\big)
\arrow[llllll, "\sim", "\text{graded dual}"'] \arrow[uuuu, hookrightarrow] \arrow[uulll,leftrightarrow,"1:1","\text{Cor \ref{DM0 ags}}"']
\end{tikzcd} \end{equation} 

\begin{Definition}\label{def:gon_coaction} 
Mimicking the construction of Brown \cite{br} we define a coaction 
\begin{equation*}
\DeltaGon: \Zf \rightarrow \Af \otimes \Zf
\end{equation*}
on the whole algebra $\Zf$ of formal multiple zeta values by setting
\[
\DeltaGon(\xi) = \Big( \amap \circ \zf  \otimes\,     \zf  \Big) \big( \DeltaGon( w) \big) , 
\]
where $w\in\QX$ is a lift for $\xi$, i.e. $w$ is any element such that $ \zf(w) = \xi$.
\end{Definition}

\begin{Theorem} The coaction $\DeltaGon:\Zf\to \Af\otimes\Zf$ is well-defined.
\end{Theorem} 

\begin{proof}  By \eqref{eq:map_DM_DM0}, we have a morphism of affine schemes
\begin{align*}
\DM_0 \times \DM &\to \DM\\
(\Phi , G) & \mapsto \Phi\glp G.
\end{align*}
Applying Yoneda's Lemma (Theorem \ref{Yoneda}) yields an algebra morphism 
\begin{align*}
    \Zf&\to\Af\otimes\Zf, \\
    \xi&\mapsto\Big(\amap\circ\zf\otimes\zf\Big)\big(\DeltaGon(w)\big), 
\end{align*}
where $w \in \QX$ satisfies $\zf(w)=\xi$. The group structure of $\DM_0$ gives rise to the coassociativity and the counitarity of this morphism.
\end{proof}

By Example \ref{exm:gon_coprod_visu}, we have 
\begin{align*}
\DeltaGon(x_0x_1) &= x_0x_1 \otimes \one + \one \otimes x_0x_1,
\end{align*}
and hence we obtain 
\begin{equation} \label{eq:def_gon_coprod_zf2}
\DeltaGon(\zf(2)) = 1\otimes \zf(2) \in \Af \otimes \Zf.     
\end{equation}

This coaction is an extension of $\DeltaGon$ on $\Af$ given in \eqref{eq:Gon_on_Zf} to $\Zf$. 
In order to give an explicit formula for it we set  for $\ve_0,\ve_{n+1} \in \X$ and  a word $\ve_1\cdots \ve_n \in \X^*$ 
\begin{equation*}
\zf(\ve_0;\ve_1\cdots\ve_n;\ve_{n+1}) = \zf(\Iformal(\ve_0;\ve_1\cdots\ve_n;\ve_{n+1})),
\end{equation*}
and then by Definition \ref{def:gon} for the Goncharov coproduct we get the explicit formula 
\begin{align}\label{eq:gon_coaction}
&\DeltaGon(\zf(\ve_1\cdots\ve_n))=
\SumGonP 
\amap \big(\zf(\ve_{i_p}; \ve_{i_p+1}\cdots\ve_{i_{p+1}-1};\ve_{i_{p+1}})\big) \otimes \zf(\ve_{i_1}\cdots\ve_{i_k})
\end{align}
where $i_0=0,\, i_{k+1}=n+1$ and $\ve_{0}=x_1,\, \ve_{n+1}=x_0$.
  
\begin{Remark}
In Definition \ref{def:gon_coaction}, it is necessary to consider $\Af$ instead of $\Zf$ in the first tensor product factor.
Recall from Example \ref{exm:gon_coprod_visu} that
\begin{align*}
\DeltaGon(x_0x_1) &= x_0x_1 \otimes \one + \one \otimes x_0x_1
\intertext{and therefore on the one hand side we get}
\DeltaGon(x_0x_1 \shuffle x_0x_1) &= \DeltaGon(x_0x_1)^{\shuffle 2} = (x_0x_1\otimes \one + \one \otimes x_0x_1)^{ \shuffle 2} \\
&= (x_0x_1\shuffle x_0x_1) \otimes \one + 2\, x_0x_1 \otimes x_0x_1 + \one \otimes (x_0x_1\shuffle x_0x_1).
\intertext{On the other side we have}
\DeltaGon(x_0x_0x_0x_1) &= x_0x_0x_0x_1\otimes \one - 2 x_0x_0 \otimes x_0x_1 + \one\otimes x_0x_0x_0x_1.
\end{align*}
Since $\frac{2}{5} \zf(2)^2 = \zf(4)$ as shown in Proposition \ref{prop:zeta_222}, both expressions should coincide up to the scalar $\frac{2}{5}$ after passing to formal multiple zeta values. However, the terms in the middle differ as $\zf(x_0x_1) = \zf(2) \neq 0$  but $\zf(x_0x_0) = 0$.
\end{Remark}

\begin{Example}\label{exm:coaction_zf5}
    We want to compute $\DeltaGon(\zf(5))$ by using the graphical interpretation of ''eating worms'' (cf. Figure \ref{fig:eating_worms}). Since $\zf(5) = \zf(x_0^4 x_1)$ we consider subwords of $x_0^4x_1$. The terms in \eqref{eq:gon_coaction} corresponding to $\one$ and $x_0^4 x_1$ are visualized, respectively, as
    \begin{center}
        \includegraphics[width=.7\textwidth]{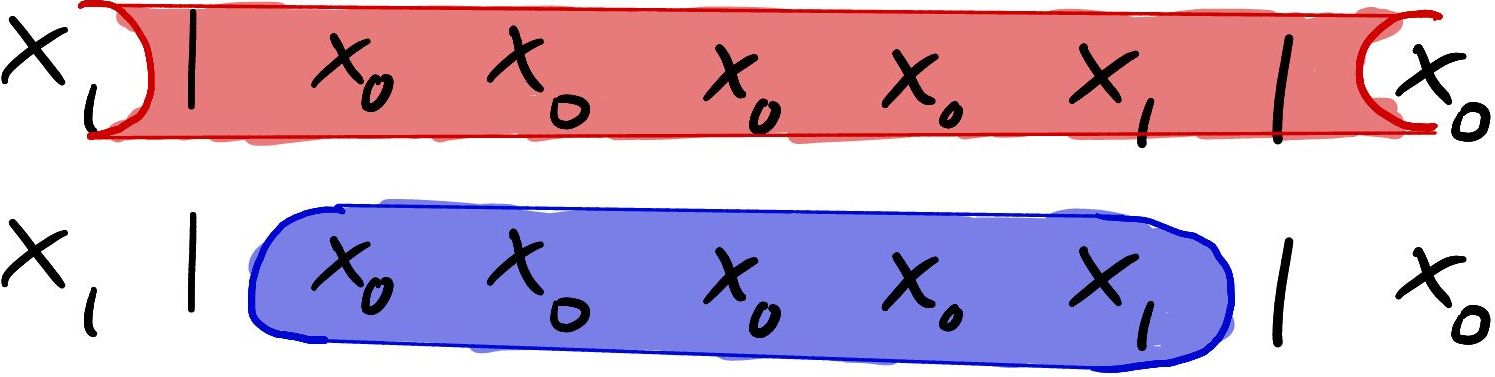}
    \end{center}
    and we obtain the summands $\amap(\zf(5))\otimes \one$ and $\one \otimes \zf(5)$, respectively. Now consider the terms that are depicted with a ''red worm'' ending on the right bound $x_0$, e.\,g.
    \begin{center}
        \includegraphics[width = .8\textwidth]{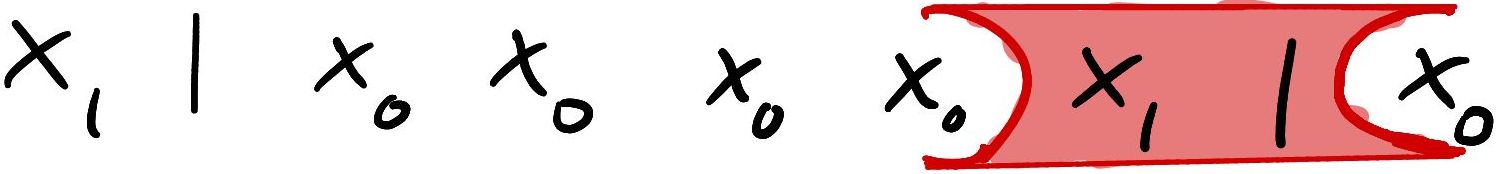}
    \end{center}
    The factor that arises from this ''red worm'' vanishes as the left bound is also $x_0$ and thus the summand does not contribute to $\DeltaGon(\zf(5))$. Similarly, the factors that correspond to a ''red worm'' starting in any $x_0$ vanish, e.\,g. 
    \begin{center}
        \includegraphics[width = .8\textwidth]{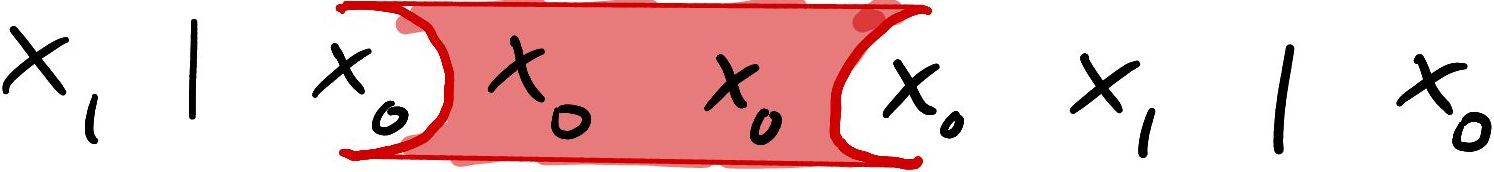}
    \end{center}
    Note that if the worm ends in $x_1$ the factor is $(-1)^i\zf(x_0^i)$ for some $i\leq 3$ which vanishes by \eqref{eq:xxxvanishes}. Hence only the first two terms contribute non-trivially and we obtain
    \begin{equation*}
        \DeltaGon(\zf(5)) = \amap(\zf(5)) \otimes \one + \one \otimes \zf(5)
    \end{equation*}
\end{Example}

This argument can be generalized to show the following lemma.

\begin{Lemma}\label{lem:coaction_primitive}
For $n\in\N$ we have 
\begin{equation*}
\DeltaGon(\zf(2n+1)) =  \amap(\zf(2n+1)) \otimes \one + \one\otimes \zf(2n+1).
\end{equation*}
\end{Lemma}

\begin{proof}
Since $\zf(2n+1) = \zf(x_0^{2n}x_1)$ each summand in $\DeltaGon \zf(2n+1)$ corresponds to a subword of $x_0^{2n}x_1$ (cf. Definition \ref{def:gon_coaction}).
Similar to Example \ref{exm:coaction_zf5}, the subwords $\one$ and $x_0^{2n}x_1$ contribute with $\amap(\zf(2n+1))\otimes\one$ and $\one\otimes\zf(2n+1)$, respectively. It thus suffices to show that all other subwords of $x_0^{2n}x_1$ contribute trivially.
So let $k\in\{1,\dots,2n\}$ and let $\ve_{i_1}\cdots\ve_{i_k}$ be a subword of $x_0^{2n}x_1$. 
If $i_k<2n+1$, then $\ve_{i_1}\cdots\ve_{i_k} = x_0^k$ and the term vanishes since $\zf(x_0^k)=0$ by \eqref{eq:xxxvanishes}.
Thus, we restrict to the case $i_k=2n+1$, i.\,e. $\ve_{i_k} = x_1$. If $k=1$ then the summand vanishes since $\zf(x_1)=0$. 
If $k>1$ then $\ve_{i_1}\cdots\ve_{i_k} = x_0^{k-1}x_1$ and the product in the formula for $\DeltaGon$ contains the factor $\zf(x_0^l)=0$ for some $l>0$. Hence all other subwords contribute trivially.
\end{proof}

\subsection{On the \texorpdfstring{\ogc}{free odd generation conjecture}} 
 
Central for this notes is the following well-known conjecture for $\dm$, which is motivated by conjectures of Deligne (\cite{de}) and Y. Ihara (\cite[p. 300]{ih})
in the context of certain Galois actions and of Drinfeld \cite{dr} on his Grothendieck-Teichm\"uller Lie algebra. By work of Furusho \cite{fu}, we know that the Grothendieck-Teichm\"uller Lie algebra embedds into $\dm$.
 
 \begin{Conjecture} \label{conj:odd_generators} The double shuffle Lie algebra $\dm$ is a free Lie algebra with exactly one generator in each odd weight $w\geq 3$, i.e.
 	\begin{equation*}
 		\dm\simeq\Lie(S)
 	\end{equation*}
 where $S=\{s_3,s_5,\ldots,s_{2n+1},\ldots\}$.
 \end{Conjecture}
 We call this conjecture the \emph{\ogcdot}

 \begin{Remark} 
 	Of course the Lie algebra
 	$\Lie(S)$
 	is also in the heart of the theory of motivic multiple zeta values (\cite{br},\cite{de},\cite{gon}, \cite{degon}). 
 	It occurs as the Lie algebra of the motivic fundamental group of $\mathbb{P}^1 \setminus \{0,1,\infty\}$ and in fact Brown proved that it also equals the Lie algebra of the motivic Galois group attached to the category of mixed Tate motives. 
 \end{Remark}

 An immediate consequence of the \ogc would be the truth of Zagier's dimension conjecture \ref{Zagier dimension conjecture} for formal multiple zeta values. 
 
 \begin{Proposition}\label{prop:ogc_implications}
 	Under the assumption of the \ogccomma one obtains
 	\[H_{\Zf}(x)=\sum_{w\geq0} \dim_\QQ\Zfw{w}x^w=\frac{1}{1-x^2-x^3}.\]
 \end{Proposition}
 
 \begin{proof} 
 	Under the assumption of the \ogc \ref{conj:odd_generators}, one obtains the following Hilbert-Poincare series for the universal enveloping algebra of $\dm$ 
 	\begin{align*} 
 		H_{\mathcal{U}(\dm)}(x)=\frac{1}{1-x^3-x^5-x^7-\dots}.
 	\end{align*}
 	By Corollary \ref{Z^f polynomial algebra}, there is an isomorphism $\Zf\simeq\QQ[\zf(2)]\otimes \mathcal{U}(\dm)^\vee$. Hence, we deduce
 	\begin{align*}
 		H_{\Zf}(x)=H_{\QQ[\zf(2)]}(x)H_{\mathcal{U}(\dm)}(x)=\frac{1}{1-x^2}\frac{1}{1-x^3-x^5-x^7-\cdots}=\frac{1}{1-x^2-x^3}.
 	\end{align*}
 \end{proof} 

Recall from Example \ref{ex:odd_free_lie}, that we have for the Lie algebra in the \ogc \ref{conj:odd_generators}
\[\mathcal{U}(\Lie (S))=(\QQ\langle s_3,s_5,\ldots\rangle,\conc,\co),\] 
and the graded dual is given by
\begin{equation*}
\CUdual = (\Q\langle s_3, s_5, s_7,\dots\rangle, \shuffle, \dec).
\end{equation*}
In Definition \ref{def:dec_Uf_general} we defined $\CUf = \CUdual\otimes \Q[s_2]$
and extended the deconcatenation coproduct on $\CUdual$ to a coaction 
\[\dec: \CUf\to \CUdual\otimes \CUf\]
via $\dec(s_2) = \one\otimes s_2$. 
We now also set $s_{2n} = b_n s_2^n$, where the rational numbers $b_n$ are given in Proposition \ref{prop:zeta_222}.

\begin{Theorem} 
\label{prop:Uf_cong_Zf}
Assume the \ogc \ref{conj:odd_generators} for $\dm$. There is an isomorphism of algebras with coaction 
\begin{equation}\label{eq:zf_iso_extended} 
\Phi\colon (\Zf, \cdot, \DeltaGon) \overset{\sim}{\longrightarrow} (\CUf, \shuffle, \dec)
\end{equation}
satisfying for each $N\geq3$ odd
\begin{equation*}
\Phi(\zf(N)) = s_N.
\end{equation*}
\end{Theorem}

\begin{proof} Assume the \ogc \ref{conj:odd_generators} for $\dm$, then
from Theorem \ref{thm:new-universal} applied to the Ihara bracket, we get a Hopf algebra isomorphism 
\[ 
(\mathcal{U} ( \dm), \glp, \co)) \cong
\mathcal{U}(\Lie ( S))=(\QQ\langle s_3,s_5,\ldots\rangle,\conc,\co).\] 
By Corollary \ref{Z^f polynomial algebra} and by dualization we have 
a Hopf algebra isomorphism
\begin{align} \label{eq:zf_iso}
\Big(\faktor{\Zf}{(\zf(2))},\ \cdot\ ,\DeltaGon\Big)\simeq (\CUdual,\shuffle,\dec).
\end{align}
It is a direct consequence of Corollary \ref{Z^f polynomial algebra} and the compatibily of the construction of the coaction on $\Zf$ via \eqref{eq:def_gon_coprod_zf2} and that on $\CUf$ by Definition \ref{def:dec_Uf_general} that we can extend \eqref{eq:zf_iso} to an isomorphism of algebras with coaction $\Phi:(\Zf,\cdot,\DeltaGon)\xrightarrow{\sim}(\CUf,\shuffle,\dec)$. 

Since we have for any odd $N\geq3$
\begin{align*}
\DeltaGon(\zf(N)) &= \zf(N) \otimes \one+\one\otimes \zf(N) ,\\
\Delta_{\text{dec}}(s_N) &=s_N  \otimes \one+\one\otimes s_N,
\end{align*}
the isomorphism $\Phi$ can be chosen such that $\Phi(\zf(N)) = s_N$.
\end{proof}

\begin{Corollary}\label{cor:z2n+1_non-zero}
Assume the \ogc \ref{conj:odd_generators} for $\dm$. The formal zeta values $\zf(2)$ and $\zf(2r+1)$, $r \in \N$, are non-zero modulo products and  algebraically independent. 
\end{Corollary}
\begin{proof}
Evidently, the letters $s_{2r+1}$ are nonzero modulo products and algebraically independent in $\CUdual$. So the same must hold for their preimages $\zf(2r+1)= \Phi^{-1}(s_{2r+1})$ under the isomorphism \eqref{eq:zf_iso_extended}. By construction of $\CUf=\CUdual\otimes \Q[s_2]$ and the convention $\zf(2)= \Phi^{-1}(s_{2})$, the claim follows.
\end{proof}

\begin{Remark} 
By work of Drinfeld \cite{dr}, Brown \cite{br}, and Furusho \cite{fu}, we have inclusions
\[\Lie(S)\subset \mathfrak{grt}_1\subset \dm.\]
This implies $\mathcal{U}(\Lie(S))\subset \mathcal{U}(\dm)$, and hence $\Phi:\Zf\to \CUf$ is surjective. In particular, we obtain Corollary \ref{cor:z2n+1_non-zero} without assuming the \ogcdot
\end{Remark}

\subsection{The Kernel conjecture}
 
Recall that in Subsection \ref{subsec:hopf_deriva} we introduced the derivations $D_{w}$ and their extension to algebras with particular coaction.  
With the notation from Subsection \ref{subsec:gon_deriva} and \ref{Gon coproduct Gon-Brown coaction} we get the commutative diagram 
\begin{equation}
\begin{aligned}
    \xymatrix{
      D_{2r+1}\colon \hspace{-1cm} &\QX \ar[rr]^-{\DeltaGon'}\ar[dd]^-{\zf} && \bigoplus_{w\geq1} \QX_w \otimes \QX \ar[rr]^-{\pi_{2r+1}\otimes \operatorname{id}}\ar[dd]^-{\amap\circ \zf\otimes \zf} && \LX_{2r+1} \otimes \QX \ar[dd]^-{\overline{\amap\circ \zf} \otimes \zf}\\
    \\
  &\Zf \ar[rr]^-{\DeltaGon'} && \bigoplus_{w\geq1} \Af_w \otimes \Zf \ar[rr]^-{\pi_{2r+1}\,\otimes \operatorname{id}} && \Lf_{2r+1} \otimes \Zf
    }
\end{aligned}
\end{equation}
where similar to the previous
\begin{align*}
\pi_{2r+1}:\bigoplus_{w\geq1} \Af_w&\to \Lf_{2r+1},\\
\overline{\amap\circ \zf}\colon \CL &\rightarrow \Lf,
\end{align*}
are the canonical projections.

Observe as in Lemma \ref{lem:D_on_H_deriv}, the maps $D_{2r+1}: \Zf\to \Lf_{2r+1}\otimes \Zf$ in the lower line  are again derivations.

\begin{Proposition}\label{prop:zf2n+1_non-zero}
     Assume the \ogc \ref{conj:odd_generators} for $\mathfrak{dm}_0$. For each $r\in\N$ the element
	\begin{equation*}
		\zf_{2r+1}= \pi_{2r+1}(\amap( \zf(2r+1))) \in \Lf_{2r+1}
	\end{equation*}
	is non-zero.
\end{Proposition} 

\begin{proof} 
This claim is just a reformulation of Corollary \ref{cor:z2n+1_non-zero}.
\end{proof}

\begin{Lemma} Let $r,n,k \ge 1$ be natural numbers, then we have
 \begin{equation}\label{eq:derivations_on_products}
		D_{2n+1}(\zf(2r+1) \cdot \zf(2)^k) = \delta_{r,n}\cdot \zf_{2r+1} \otimes \zf(2)^k.
	\end{equation}
where $\delta_{r,n}$ is the Kronecker delta.
 \end{Lemma}

\begin{proof}  Recall $\DeltaGon(\zf(2 )^k) =   \one\otimes \zf(2 )^k$ and
\begin{align*}
\DeltaGon(\zf(2r+1)) = \amap( \zf(2r+1)) \otimes \one + \one\otimes \zf(2r+1).
\end{align*}		
We then compute
\begin{align*}	  
&\DeltaGon'(\zf(2r+1) \zf(2)^k)\\
&=\DeltaGon(\zf(2r+1) \zf(2)^k) - \one\otimes \zf(2r+1) \zf(2)^k \\
&=\DeltaGon(\zf(2r+1)) \DeltaGon(\zf(2)^k) - \one\otimes \zf(2r+1) \zf(2)^k   \\
&=   \Big(\amap(\zf(2r+1))\otimes \one + \one\otimes \zf(2r+1)\Big ) \Big(\one\otimes \zf(2)^k\Big)- \one\otimes \zf(2r+1)\zf(2)^k  \\
&=\amap(\zf(2r+1))\otimes \zf(2)^k,
\end{align*}
which in turn yields
\begin{align*}
D_{2n+1}& ( \zf(2r+1)\zeta(2)^k ) = \pi_{2n+1}(\amap(\zf(2r+1))) \otimes \zf(2)^k. \qedhere
\end{align*}
\end{proof}

\begin{Conjecture}[Kernel conjecture] \label{conj:kernel}
Define
\begin{equation}\label{eq:map_D<N}
	D_{<N} = \bigoplus_{3\leq 2r+1<N} D_{2r+1},
\end{equation}
then we have $ \ker(D_{<N})  \cap \Zf_N  =\Q\, \zf(N)$  for all $N\geq 2$.
\end{Conjecture}

\begin{Remark}
    Observe that $\zf(N)\in\ker(D_{<N})$ for all $N\geq 2$. This follows from Proposition \ref{prop:zeta_222} for even $N$ and from Lemma \ref{lem:coaction_primitive} together with Corollary \ref{cor:z2n+1_non-zero} for odd $N$.
\end{Remark}

	\begin{Theorem}\label{thm:ker_CD} The \ogc \ref{conj:odd_generators} for $\mathfrak{dm}_0$  implies the Kernel conjecture \ref{conj:kernel}, i.e., we have for all $N\geq 2$
		\begin{equation*}
			\ker(D_{<N}) \cap \Zf_N = \Q \zf(N).
		\end{equation*}
	\end{Theorem}
	\begin{proof}
		Recall that $\zf(N)\in \ker(D_{<N})$. 
For the other inclusion we employ the isomorphism of   algebras with coaction
\begin{equation*}
			\Phi\colon (\Zf, \cdot, \DeltaGon) \overset{\sim}{\longrightarrow} (\CUf, \shuffle, \Delta_{\text{dec}})
		\end{equation*}		
from  Proposition \ref{prop:Uf_cong_Zf}. In particular this yields for all $N\in\N$
a commutative diagram
\begin{equation}\label{eq:isomorphism_D}
\xymatrix{ 
\Zf_N \ar[d]^{D_{<N}} \ar[rr]^ \Phi  && \CUf_N \ar[d]^{D_{<N}}\\
\bigoplus\limits_{1 < 2r+1 <N} \!\!\! \Lf_{2r+1}\otimes \Zf   \ar[rr]^{\overline{\Phi} \otimes \Phi} &&\bigoplus\limits_{1 < 2r+1 <N} \!\!\! \LUw{2r+1}\otimes \CUf
}
\end{equation}
		Hence, if  $\xi\in\ker(D_{<N})\cap \Zf_N$, then   \eqref{eq:isomorphism_D} implies $\Phi(\xi)\in\ker(D_{<N})$. So by the Proposition \ref{prop:kernel_U_extended}, we obtain that $\Phi(\xi) = \alpha\, s_N$ for some $\alpha\in\Q$. Hence $\xi = \alpha\, \Phi^{-1}(s_N) = \alpha\, \zf(N)$.
	\end{proof}

\begin{Example}
    The direct sum in \eqref{eq:map_D<N} is empty for $N=3$, hence $D_{<3}\equiv 0$ is the zero map. 
    Conjecture \ref{conj:kernel} would thus imply that $\Zf_3 = \Q\, \zf(3)$, so in particular $\zf(2,1)$ would be a rational multiple of $\zf(3)$ 
    (possibly zero). In fact, in Example \ref{exm:stuffle_shuffle_prod} we saw that Euler's identity $\zeta(3) = \zeta(2,1)$ follows from the 
    extended double shuffle relations, hence $\zf(3) = \zf(2,1)$.
\end{Example}

\begin{Example} We give  here some potential applications of the Kernel conjecture. Let $n\in\NN$ be a positive integer. 
\begin{enumerate}
    \item We compute $D_{<2n}(\zf(\{2\}^n))$. Let $r\in\NN$ with $2r+1<2n$. Each term in $D_{2r+1}(\zf(\{2\}^n))$ corresponds to a consecutive 
    subword of $(x_0x_1)^n\in\X^*$ of weight $2r+1$. Note that such a subword always starts and ends in the same letter for parity reasons. 
    This can be depicted on a semicircle as:
    \def\scale{5}
    \def\n{15}
    \begin{center}
    \begin{tikzpicture}[baseline = (current bounding box.north)]
    	\coordinate (start) at (-\scale,0);
    	\coordinate (end) at (\scale,0);
    	
    	\begin{scope}[nodes={fill=black, circle, scale=.5}]
    		\clip (-\scale-.1,0) rectangle (\scale+.1,\scale+.1);
    		\draw[thick] (start)+(end) circle (\scale);
    	\end{scope}
    		
    	\draw[thick] (start) -- (end);
    	
    	\foreach \num in {0,1,...,\n}{
    		\coordinate (L\num) at (180-\num*180/\n : 1.1*\scale);
    		
    		\coordinate (N\num) at (180-\num*180/\n : \scale);
    	}

        \draw[thick, bend right = 15] (N0) to (N1)
            to (N2)
            to (N3)
            to (N4);
     
    	\draw[thick] (N4) to (N10);

        \draw[thick, bend right = 15] (N10) to (N11)
            to (N12)
            to (N13)
            to (N14)
            to (N15);

    	\foreach \num in {0,1,...,5,7,8,...,\n}{
    		\node[fill=black, circle, scale=.5] at (N\num) {};
    	}

        \node at (L0) {$x_1$};
        \node at (L\n) {$x_0$};
        
        \foreach \num in {2,4,10,12,14}{
    		\node[blue] at (L\num) {$x_1$};
    	}

        \node[red] at (L8) {$x_1$};

        \foreach \num in {1,3,11,13}{
    		\node[blue] at (L\num) {$x_0$};
    	}
        
    	\foreach \num in {5,7,9}{
    		\node[red] at (L\num) {$x_0$};
    	}
    	
    	\node[red] at (L6) {$\dots$};
    \end{tikzpicture}
    \end{center}
    So we have $\zf(\{2\}^n)\in\ker(D_{2r+1})$ for all $2r+1<2n$, hence $\zf(\{2\}^n)\in\ker(D_{<2n})$.
    We therefore expect by Conjecture \ref{conj:kernel} that $\zf(\{2\}^n)\in \Q\,\zf(2n)$. In fact, this holds by Proposition \ref{prop:zeta_222}. 
    A generalisation of this method can be found in \cite{C15}.
    \item We compute $D_{<n+2}(\zf(2,\{1\}^n))$. Let $r\in\NN$ with $2r+1<n+2$. Since $\zf(2,\{1\}^n) = \zf(x_0x_1^{n+1})$ there are only two kind of strict subwords whose boundaries are not both $x_1$. They can be depicted as
    \begin{center}
        \includegraphics[width = .8\textwidth]{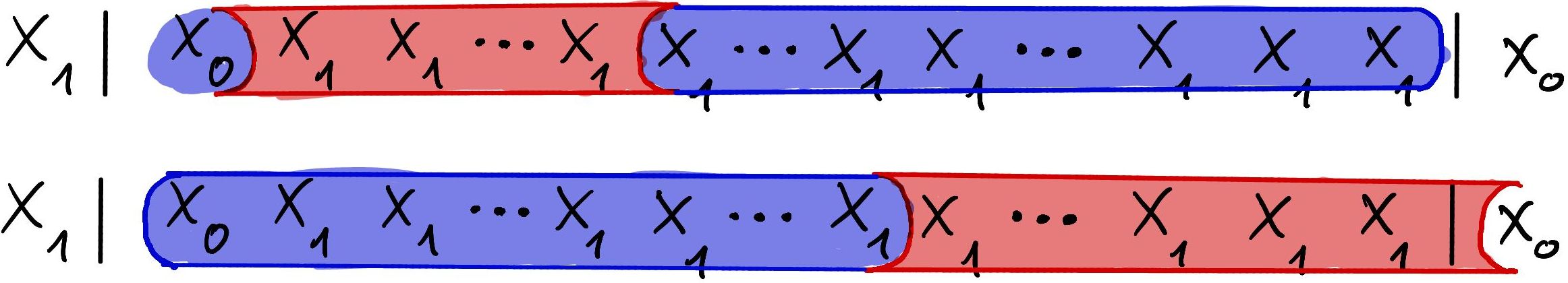}
    \end{center}
    Both corresponding terms vanish since $\zf(x_1^{2r+1}) = 0$ by \eqref{eq:xxxvanishes}. So we have $\zf(2,\{1\}^n)\in\ker(D_{2r+1})$ for all $2r+1 < n+2$, hence $\zf(2,\{1\}^n)\in\ker(D_{<n+2})$. By the Kernel conjecture \ref{conj:kernel} we expect that $\zf(2,\{1\}^n)\in \Q\, \zf(n+2)$. 
    Using the evaluation map would imply
    \begin{equation*}
        \zf(2,\{1\}^n) = \zf(n+2).
    \end{equation*}
  as $\zeta(2,\{1\}^n) = \zeta(n+2)$.  
  Remarkably, there are two independent ways to show that multiple zeta values satisfy this identity. At first  it is an example of \emph{duality relations}, 
    i.\,e. $\zeta(x_0^{k_1-1}x_1\cdots x_0^{k_d-1}x_1) = \zeta(x_0x_1^{k_d-1}\cdots x_0x_1^{k_1-1})$ for integers $k_1>1$ and $k_2,\dots,k_d\geq 1$. 
   Alternatively, as shown in \cite{k}, this identity also follows from the extended double shuffle relations and therefore we also get a unconditional proof without refering to the Kernel conjecture in this case.
    \end{enumerate}
\end{Example}

\section{Brown's theorem for formal multiple zeta values}

In this section we are interested in formulas for the multiple zeta values with level zero and one (cf Definition \ref{def:level}).
We already understand
$H(n)= \zf(\{2\}^n)$, since by Proposition \ref{prop:zeta_222} we have explicit formulae for
\begin{align*}
\zf(2 k) \cdot \,\zf(\{2\}^{n-k}) \in \mathbb{Q} \cdot  \zf(2)^n.
\end{align*} 

\begin{Theorem}[Zagier, 2012]
	Let $a,b\geq 0$. Then we have an identity of multiple zeta values
	\begin{equation*}
		\zeta(\{2\}^a,3,\{2\}^b) = \sum_{r=1}^{a+b+1} c_{a,b}^r\, \zeta(2r+1) \zeta(\{2\}^{a+b+1-r})
	\end{equation*}
	with $c_{a,b}^r$ from Definition \ref{def:coefficients_zagiers_formula}, i.e.
    \begin{equation*}
        c_{a,b}^r = 2\cdot (-1)^r \Bigg(\binom{2r}{2b+2} - (1 - 2^{-2r}) \binom{2r}{2a+1}\Bigg).
    \end{equation*}
\end{Theorem}

Brown showed that this formula also holds for motivic MZVs. We will prove the analogue statement for formal MZVs, assuming the Kernel conjecture \ref{conj:kernel}.

\subsection{Proof of Zagier's formula} 

Zagier's theorem is a vast refinement of Proposition \ref{prop:sum_level_one}, i.e. of the identity
\begin{align*}
\sum_{a+b=n}  \zf(\{2\}^a,3,\{2\}^b) = \sum_{i=0}^{n}  (-1)^i \zf(2i+3) \zf(\{2\}^{n-i}).
\end{align*}
Unfortunately, by now only analytical proofs of his result are known. We follow in our presentation the paper 
\cite{LaiLupuOrr}. 

\begin{Theorem} \label{thm:zagier}
Let $a,b \ge 0$ be integers, then the following three quantities are equal
\begin{align}
    H(a,b) &=  \zeta(\{2\}^{a},3,\{2\}^{b}) \\
     I(a,b) &= \frac{\pi^{2a} 2^{2b+3}}{(2a+1)! (2b+2)!} 
    \int\limits_0^{\pi/2} x^{2b+2}  \Big( 1- \frac{2x}{\pi} \Big)^{2a+1} \cot(x) dx \\
    \hat H(a,b) &= \sum_{r=1}^{a+b+1} c_{a,b}^r\, \zeta(2r+1)  \zeta(\{2\}^{a+b+1-r}).
    \end{align}
\end{Theorem}

Essential ingredients for the proof are the following three lemmas.

\begin{Lemma} \label{lem:lem1_level_one} For any polynomial $P(x) \in \mathbb{C}[x]$ with $P(0)=P(1)= 0$, we have
\begin{align}
\int\limits_0^1 P(x) \cot\big( \frac{\pi x}{2}\big) dx &= 
2 \sum_{ k=1}^{\lfloor \deg(P)/2 \rfloor} (-1)^k 
\Big( \big( 1- 2^{-2k} \big)  P^{(2k)}(1) + P^{(2k)}(0) \Big) 
\frac{\zeta(2k+1)}{\pi^{2k+1}} , 
\end{align}
where $ P^{(2k)}(x)$ denotes the $(2k)$-th derivative of $P(x)$.
\end{Lemma}
\begin{proof}
    This is a special case of \cite[Lemma~2.3]{LaiLupuOrr}.
\end{proof}

\begin{Lemma} \label{lem:lem2_level_one}Let $r\in \mathbb{N}$, then we have
\begin{align}
\frac{\arcsin(x)^{2r}}{(2r)!} = \frac{1}{4^r} 
\sum_{n=1}^\infty \frac{4^n}{n^2 \binom{2n}{n}} \,x^{2n}
\sum\limits_{n_1<n_2< ... <n_{r-1}<n} \frac{1}{ n_1^2 n_2^2 ... n_{r-1}^2}.
\end{align}
\end{Lemma}

\begin{proof}
    For a proof we refer to \cite[Eq.~(2)]{LaiLupuOrr}.
\end{proof}

\begin{Lemma} \label{lem:lem3_level_one} Let $n\geq1$, $a\geq0$, then we have
\begin{align}
\int\limits_0^1 x^{2n-1} \frac{  ( 2 \arccos(x))^{2a+1}}{(2a+1)!} dx = 
 \frac{ \binom{2n}{n} \pi }{2^{2n +1} n}
 \sum_{n=1}^\infty   
\sum\limits_{n< m_1<m_2< ... <m_{a}} \frac{1}{ m_1^2 m_2^2 ... m_{a}^2},
\end{align}
where the sum is set to be $1$ if $a=0$.
\end{Lemma}

\begin{proof}
    For a proof we refer to \cite[Lemma~3.1]{LaiLupuOrr}.
\end{proof}
We are now prepared for the proof of Zagier's theorem. 

\begin{proof}[Proof of Theorem \ref{thm:zagier}] 
At first we observe that via the change of variables $x \mapsto \dfrac{\pi x}{2} $ we obtain
\begin{align*}
I(a,b) &= \frac{\pi^{2a} 2^{2b+3}}{ (2a+1)! (2b+2)!} 
\int\limits_0^{\pi/2} x^{2b+2}  \Big( 1- \frac{2x}{\pi} \Big)^{2a+1} \cot( x) dx \\
& = \frac{\pi^{2a+2b+2} 2}{ (2a+1)! (2b+2)!}
\int\limits_0^{1} x^{2b+2}  \big( 1-  x \big)^{2a+1} \cot\big( \frac{\pi x}{2}\big) dx
\end{align*}
Applying Lemma \ref{lem:lem1_level_one} with $P(x)= x^{2b+2}  (1-x)^{2a+1}$ yields 
\begin{align*}
&\int\limits_0^{1} x^{2b+2}  \big(1-x \big)^{2a+1} \cot\big( \frac{\pi x}{2}\big) dx \\
&=  (2a+1)! (2b+2)! \\
&\cdot 2 \sum_{k=0}^{a+b+1} (-1)^k \Big[ 
\binom{2k}{2b+2} - \big( 1 - 2^{-2k}\big) \binom{2k}{2a+1} \Big] \frac{\zeta(2k+1)}{(2a+2b+3-2k)! \pi^{2k+1}}
\end{align*}
Now using Proposition \ref{prop:zeta_222} the equality $I(a,b)=\hat H(a,b)$ follows by rearranging the above equalities.

To obtain the second equality we start with the change of variables $x \mapsto \arcsin(x) $. Using in addition that $\arccos(x) = \pi/2 - \arcsin(x)$ we get
\begin{align*}
I(a,b) &= \frac{2^{2a+2b+4}}{ (2a+1)! (2b+2)! \pi} 
\int\limits_0^{\pi/2} x^{2b+2}  \Big( \frac{\pi}{2}- x  \Big)^{2a+1} \cot(x) dx \\
& =  \frac{2}{\pi} 
\int\limits_0^{1} \frac{ (2\arcsin(x))^{2b+2}}{(2b +2)!} \frac{(2\arccos(x))^{2a+1}}{(2a+1)!} \frac{1}{x}   dx.
\end{align*}

If we replace the first factor under the integral by its Taylor expansion given in Lemma 
\ref{lem:lem2_level_one}, then the resulting sum of integrals are exactly those of \ref{lem:lem3_level_one}, therefore  we get

\begin{align*}
&\int\limits_0^{1} \frac{ (2 \arcsin(x))^{2b+2}}{(2b +2)!}  \frac{(2 \arccos(x))^{2a+1}}{(2a+1)!} \frac{1}{x}   dx \\
&= \sum_{n=1}^\infty \frac{4^n}{n^2 \binom{2n}{n}}
\sum\limits_{n_1<n_2< ... <n_b<n} \frac{1}{ n_1^2 n_2^2 \cdots n_b^2} 
 \int\limits_0^1 x^{2n-1} \frac{  ( 2 \arccos(x))^{2a+1}}{(2a+1)!} dx\\
&= \frac{\pi}{2} \sum_{n=1}^\infty \frac{1}{n^3} \cdot
\left( \sum\limits_{n_1<n_2< ... <n_b<n} \frac{1}{ n_1^2 n_2^2 ... n_b^2} \right)
\cdot \left( \sum\limits_{n< m_1<m_2< \cdots <m_{a}} \frac{1}{ m_1^2 m_2^2 \cdots m_{a}^2}  \right)\\ 
&= \frac{\pi}{2} \sum\limits_{n_1<n_2< ... <n_{b}<n< m_1<m_2< ... <m_{a} } 
\frac{1}{  n_1^2 n_2^2 ... n_{b}^2 n^3 m_1^2 m_2^2 ... m_{a}^2} \\
&= \frac{\pi}{2} H(a,b). \qedhere
\end{align*}
\end{proof}

\subsection{Lifting Zagier's formula}
\label{subsec:lifting_zagier_to_fMZV}

We will show in a sequence of lemmas that both sides of Zagier's formula $H(a,b) = \hat{H}(a,b)$ in the formal setting have the same image under the map $D_{<2a+2b+3}$. Thus by  Theorem \ref{thm:ker_CD} their difference is a multiple of $\zf(2a+2b+3)$. We determine this multiple
via the canonical map $\Zf\to \Z$.

 \begin{Lemma}  \label{lem:D_on_23}
    For integers $a,b\geq 0$ and $1\leq r\leq a+b$ we have
    \begin{equation*}
        D_{2r+1}\big(\zf(\{2\}^a, 3, \{2\}^b)\big) = \pi_{2r+1}(\amap(\xi_{a,b}^r)) \otimes \zf(\{2\}^{a+b+1-r})
    \end{equation*}
    where
    \begin{align*}
        \xi_{a,b}^r &= \sum_{\substack{0\leq \alpha\leq a\\ 0\leq \beta\leq b\\ \alpha+\beta+1=r}} \zf(\{2\}^\alpha, 3, \{2\}^\beta) - \sum_{\substack{0\leq \alpha < a\\ 0\leq \beta\leq b\\ \alpha+\beta+1=r}} \zf(\{2\}^\beta, 3, \{2\}^\alpha) \\
        &+ 2\, (\I(a\geq r) - \I(b\geq r)) \sum_{i=1}^{r} (-1)^i \zf(2i+1) \zf(\{2\}^{r-i})
    \end{align*}
    where $\I$ denotes the indicator function (i.\,e. $\I(A) = 1$ if $A$ is a true statement and $\I(A) = 0$ else).
\end{Lemma}

\begin{proof}
 	We have
	\begin{equation*}
		\zf(\{2\}^a, 3, \{2\}^b) = \zf((x_0 x_1)^a x_0x_0x_1 (x_0x_1)^b).
	\end{equation*}
	By previous considerations, each summand in $D_{2r+1}(\zf(\{2\}^a, 3,\{2\}^b))$ corresponds to a consecutive subsequence of length $2r+3$ (including boundaries). Let $\ve_0\dots\ve_{2r+2}$ be such a subsequence. \\
	 Cases 1 + 2:  $\ve_{2r+2}$ is to the left or $\ve_0$ is to the right of $x_0x_0$. Then, by parity, $\ve_0 = \ve_{2r+2}$ and the sequence does not contribute.
	The general observation here is that the right tensor factor must be $\zf(\{2\}^{a+b+1-r})$ since $x_0x_0x_1$ must be contained $\ve_0\dots\ve_{2r+2}$. \\
	Cases 3 + 4: The sequence $x_0x_0x_1$ is contained in $\ve_1\dots\ve_{2r+2}$.
	In this case, we obtain factors of the form
	\begin{align*}
		&\zf(x_1; (x_0x_1)^\alpha x_0x_0x_1 (x_0x_1)^\beta; x_0)\vspace{-.3cm} 
		&&\alpha\leq a, \beta\leq b\\
		&\zf(x_0; (x_1x_0)^\alpha x_1x_0x_0 (x_1x_0)^\beta; x_1)
		&&\alpha<a, \beta\leq b
	\end{align*}
	such that $\alpha + \beta + 1 = r$ and $\alpha,\beta\geq 0$. Observe that these factors are $\zf(\{2\}^\alpha, 3, \{2\}^\beta)$ and $-\zf(\{2\}^\beta, 3, \{2\}^\alpha)$, respectively.\\
	Cases 5 + 6:
	Assume that $\ve_0\dots\ve_{2r+2}$ starts or ends in $x_0x_0$. This is only possible, if $b\geq r$ or $a\geq r$. The corresponding factors are $\zf(x_1;w;x_0)$ and $\zf(x_0;w;x_1)$ where $w = (x_0x_1)^r x_0$ is palindromic. Hence $\zf(x_0;w;x_1) = - \zf(x_1;w;x_0)$. The claim now follows from Proposition \ref{prop:sum_level_one}, i.\,e. 
	\begin{equation*}
		\zf((x_0x_1)^r x_0) = 2\sum_{i=1}^r (-1)^{i} \zf(2i+1)\, \zf(\{2\}^{r-i}). \qedhere
	\end{equation*}
\end{proof}

\begin{Lemma}  \label{lem:unique_coeffs}
	Assume the Kernel conjecture \ref{conj:kernel}, and let $a,b\geq 0$ be given integers. Then there exist unique coefficients $\alpha_1,\dots,\alpha_n\in\Q$ such that
    \begin{equation*}
        \zf(\{2\}^a,3,\{2\}^b) = \sum_{i=1}^n \alpha_i\, \zf(2i+1) \zf(\{2\}^{n-i})
    \end{equation*}
    with $n = a+b+1$.
\end{Lemma} 
  
\begin{proof}
    Let $w = (x_0x_1)^a x_0x_0x_1 (x_0x_1)^b$. 
    We proof the statement by induction on $\wt(w)$. The claim is trivial for $\wt(w) = 3$. 
	Let $w = \{2\}^a 3 \{2\}^b$ with $\wt(w) = 2N+1$ and assume the claim holds for all $r < N$. Now Lemma \ref{lem:D_on_23}  implies for all $1\leq r < N$
	\begin{equation*}
		D_{2r+1}(\zf(w)) = \pi_{2r+1}(\amap(\xi_{a,b}^r)) \otimes \zf(\{2\}^{N-r}).
	\end{equation*}
	Recall also from Lemma  \ref{lem:D_on_23} the explicit formula for $\xi_{a,b}^r$ and observe that $\wt(\xi_{a,b}^r) \leq \wt(w)-2$ for all $r\leq a+b$.
 So we can apply the induction hypothesis to all summands in $\xi_{a,b}^r$. Therefore, $\pi_{2r+1}(\amap(\xi_{a,b}^r)) = \alpha_r\, \zf_{2r+1}$.
	Hence
	\begin{equation*}
		D_{2r+1}(\zf(w)) = \alpha_r\, \zf_{2r+1} \otimes \zf(\{2\}^{N-r}).
	\end{equation*}
	In particular, we obtain numbers $\alpha_r\in\Q$ for each $1\leq r \leq N-1$ such that
	\begin{equation*}
		D_{<2N+1}(\zf(w)) = \sum_{r=1}^{N-1} \alpha_r\, \zf_{2r+1} \otimes \zf(\{2\}^{N-r}).
	\end{equation*}
	
	On the other hand, we deduce from \eqref{eq:derivations_on_products} that
	\begin{equation*}
		D_{2r+1} \Big(\sum_{i=1}^{N-1} \alpha_i \zf(2i+1) \zf(\{2\}^{N-i})\Big) = \alpha_r\, \zf_{2r+1} \otimes \zf(\{2\}^{N-r}),
	\end{equation*}
	because $D_{2r+1}$ is a derivation. Therefore,
	\begin{equation*}
		D_{<2N+1}\Big(\zf(w) - \sum_{i=1}^{N-1} \alpha_i\, \zf(2i+1) \zf(\{2\}^{N-i}) \Big) = 0.
	\end{equation*}
	By the assumption $\ker(D_{<2N+1})\cap \Zf_{2N+1} = \Q \zf(2N+1)$ it follows that both sides of the claim differ by $\alpha_n\, \zf(2N+1)$ for some $\alpha_n\in\Q$ and the claim follows.
\end{proof}

	\begin{Remark}
		Note that the linear combination
		\begin{equation*}
			\Phi(\zf(w)) = \sum_{i=1}^n \alpha_i\, s_{2i+1} s_2^{n-i}
		\end{equation*}
		is unique where $\Phi\colon \Zf \overset{\sim}{\longrightarrow}\CUf$ is the conjectured isomorphism from \eqref{eq:zf_iso_extended}.
		This may give an alternative proof.
	\end{Remark}

Recall the numbers 
\begin{equation*}
    c_{a,b}^r = 2\cdot (-1)^r \Bigg( \binom{2r}{2b+2} - (1-2^{-2r}) \binom{2r}{2a+1}\Bigg)
\end{equation*}
from Definition \ref{def:coefficients_zagiers_formula} that we studied in Section \ref{subsec:numbers_c}.

\begin{Theorem}\label{thm:zagier_formal_mzv}
	Assume the Kernel conjecture \ref{conj:kernel}, and let $a,b\geq 0$ be given integers. Then we have
	\begin{equation}\label{eq:formal_zagier_formula}
		\zf(\{2\}^a,3,\{2\}^b) = \sum_{r=1}^{a+b+1} c_{a,b}^r\, \zf(2r+1) \zf(\{2\}^{a+b+1-r}).
	\end{equation}
\end{Theorem}

\begin{proof}
    Induction on the weight. Let $w = (x_0x_1)^a x_0x_0x_1(x_0x_1)^b$ and $N = a+b$. By Lemma \ref{lem:unique_coeffs} we have 
    \begin{equation*}
    	\zf(\{2\}^\alpha, 3,\{2\}^\beta) = \sum_{r=1}^{\alpha+\beta+1} \alpha_r \, \zf(2r+1) \zf(\{2\}^{\alpha+\beta+1-r})
    \end{equation*}
    We assume that $\alpha_r=c_{\alpha,\beta}^r$
    holds for all $\alpha,\beta\geq 0$ such that $\alpha+\beta < N$. By Lemma \ref{lem:D_on_23} we have
    \begin{equation*}
    	D_{2r+1}(\zf(w)) = \pi_{2r+1}(\amap(\xi_{a,b}^r)) \otimes \zf(\{2\}^{a+b+1-r})
    \end{equation*}
    for all $1\leq r\leq N$. 
    Applying the induction hypothesis to the explicit formula for $\xi_{a,b}^r$ given in Lemma \ref{lem:D_on_23} yields modulo products
    \begin{equation*}
    	\pi_{2r+1}(\amap(\xi_{a,b}^r)) = \Big(\sum_{\substack{0\leq \alpha\leq a\\ 0\leq \beta\leq b\\ \alpha+\beta+1=r}} c_{\alpha,\beta}^r - \sum_{\substack{0\leq \alpha < a\\ 0\leq \beta\leq b\\ \alpha+\beta+1=r}} c_{\beta,\alpha}^r + 2\cdot (-1)^r \big(\I(a\geq r) - \I(b\geq r)\big)\Big) \zf_{2r+1}.
    \end{equation*}
    
    By Lemma \ref{lem:binomials} we have
    \begin{equation}\label{eq:binom_identity}
         c_{a,b}^r = \sum_{\substack{0\leq \alpha\leq a\\ 0\leq \beta\leq b\\ \alpha+\beta+1=r}} c_{\alpha,\beta}^r - \sum_{\substack{0\leq \alpha < a\\ 0\leq \beta\leq b\\ \alpha+\beta+1=r}}  c_{\beta,\alpha}^r + 2\cdot (-1)^r \big(\I(a\geq r) - \I(b\geq r)\big).
    \end{equation}    
    
    Hence for all $1\leq r\leq N$ we have
    \begin{equation*}
    	D_{2r+1}\Big(\zf(\{2\}^a, 3, \{2\}^b)\Big) = c_{a,b}^r \, \zf_{2r+1} \otimes \zf(\{2\}^{a+b+1-r}).
    \end{equation*}
    This implies with formula \eqref{eq:derivations_on_products}  applied to the products on the right hand side of the claimed identity \eqref{eq:formal_zagier_formula} that
    \begin{equation*}
    	D_{<2N+3}\Big( \zf(\{2\}^a,3,\{2\}^b) - \sum_{r=1}^{a+b+1} c_{a,b}^r\, \zf(2r+1) \zf(\{2\}^{a+b+1-r})\Big) = 0.
    \end{equation*}
    By assumption, there is some $\alpha\in\Q$ such that
    \begin{equation}\label{eq:factor_alpha}
    	\alpha\, \zf(2N+3) = \zf(\{2\}^a,3,\{2\}^b) - \sum_{r=1}^{a+b+1} c_{a,b}^r\, \zf(2r+1) \zf(\{2\}^{a+b+1-r}).
    \end{equation}
    Applying the surjective morphism $\Zf\twoheadrightarrow \CZ$ from \eqref{evaluation map} to \eqref{eq:factor_alpha} yields $\alpha = 0$, because we have Zagier's theorem.
\end{proof}

Applying the projection modulo products to \eqref{eq:formal_zagier_formula}, we deduce the following for the sets $\mathcal{B}^1$ and $\B$ from Definition \ref{def:B_1} and \ref{def:level}.

\begin{Corollary}\label{cor:level_one_mod_prod_fmzv}
    Assuming the Kernel conjecture \ref{conj:kernel}, we have 
    \[  \pi_{2r+1}\big(\amap\big(\zf(\mathcal{B}^1)\big)\big) = \Q\, \zf_{2r+1},
    \]
    and in particular
    \[
    \pi_{2r+1}\big(\amap\big(\zf(\gr_1^F(\B))\big)\big) = \Q\, \zf_{2r+1}.
    \]
\end{Corollary}

\subsection{The final step in Brown's proof}

\begin{Definition}
	 We set
	\begin{equation*}
		\ZfH = \operatorname{span}_\Q\{\zf(k_1,\dots,k_d) \mid d\geq 0, k_i\in\{2,3\}\}.
	\end{equation*}
\end{Definition}

On our way to prove that the \ogc implies $\ZfH = \Zf$ the next 
ingredient is a detailed study of the coaction on $\ZfH$.  
We have $\zf( \B) = \ZfH$ and thus 
by Lemma \ref{lem:coaction_B23} $\DeltaGon$ restricts to
\begin{equation*}
	\DeltaGon\colon \ZfH \rightarrow \Af \otimes \ZfH,
\end{equation*}
 
 According to Definition \ref{def:level}, we are interested in the following filtration.
\begin{Definition}
	We define a \emph{level filtration} on $\ZfH$ by setting 
	\begin{equation*}
		F_\ell \ZfH = \operatorname{span}_\Q\Big\{\zf(k_1,\dots,k_d) \in\ZfH \;\Big|\; \deg_3(x_0^{k_1-1}x_1\cdots x_0^{k_d-1}x_1) \leq \ell\Big\}
	\end{equation*}
	for all $\ell\in\N_0$.
\end{Definition}

 As a corollary of Lemma \ref{lem:coaction_B23} we get 
	\begin{equation*}
		\DeltaGon\colon F_\ell \ZfH \rightarrow \Af \otimes F_\ell\ZfH.
	\end{equation*}
and Proposition \ref{prop:level_red_derivations_X} 1. implies
 \begin{equation}\label{eq:lem:D_reduces_level_H23}
		D_{2r+1}(F_\ell \ZfH) \subseteq \Lf_{2r+1} \otimes F_{\ell-1} \ZfH.
\end{equation}
 
Again we follow Definition \ref{def:ass_graded_B23}, and introduce the following.
 
\begin{Definition}
	For all $\ell\in\N$ we set
	\begin{equation*}
		\gr_\ell^F(\ZfH) = \faktor{F_\ell \ZfH}{F_{\ell-1} \ZfH}
	\end{equation*}
	and  $\gr_0^F= F_0\ZfH$. For all $\ell\in\N$ we denote by $\pi_\ell^F\colon F_\ell \ZfH \twoheadrightarrow \gr_\ell^F\ZfH$
the natural projections.
\end{Definition}	
	
By  \eqref{eq:lem:D_reduces_level_H23} the map $D_{2r+1}^{(\ell)} = (\operatorname{id}\otimes \pi_{\ell-1}^F)\circ D_{2r+1}\rvert_{F_\ell\ZfH}$ induces a map 
	\begin{equation}\label{eq:D_reduces_level_H23} 
		D_{2r+1}^{(\ell)} \colon \gr_\ell^F \ZfH \rightarrow \Lf_{2r+1} \otimes \gr_{\ell-1}^F \ZfH
	\end{equation}
for all $r,\ell\geq 1$.	We have the following refinement of   \eqref{eq:D_reduces_level_H23}, that combines
Theorem \ref{thm:graded_image_D_B23} with the lift of Zagier's formula \eqref{eq:formal_zagier_formula}.

\begin{Theorem}\label{thm:graded_image_D_H23} 
Assume the Kernel conjecture \ref{conj:kernel} holds.  
Let $r,\ell\geq 1$. Then for each weight $N\geq 2r+1$, we have 
	\begin{equation*}
		D_{2r+1}^{(\ell)}\colon \gr_{\ell}^F (\ZfH)_N \rightarrow \Q\, \zf_{2r+1} \otimes \gr_{\ell-1}^F (\ZfH)_{N-2r-1}.
	\end{equation*}
\end{Theorem}

\begin{proof}
    Recall the derivations
    \[
        D_{2r+1}^{(\ell)}\colon \gr_\ell^F(\B)_N \rightarrow \pi_{2r+1}(\mathcal{B}^1) \otimes \gr_{\ell-1}^F (\B)_{N-2r-1}
    \]
    from Proposition \ref{prop:level_red_derivations_X}. 
    As a consequence of lifting Zagier's formula \eqref{eq:formal_zagier_formula}, we observed in Corollary \ref{cor:level_one_mod_prod_fmzv} that $\pi_{2r+1}\big(\amap\big(\zf(\mathcal{B}^1)\big)\big) = \Q\, \zf_{2r+1}$.
    The claim now follows by applying the maps $\zf$ and $\overline{\amap\circ \zf}$ to the dervations above since we have $\overline{\amap\circ\zf}(\pi_{2r+1}(\mathcal{B}^1)) \subseteq \pi_{2r+1}(\amap(\zf(\mathcal{B}^1)))$. 
\end{proof}

Assuming the \ogc \ref{conj:odd_generators}, Proposition \ref{prop:zf2n+1_non-zero} states that $\zf_{2r+1}$ is non-zero.
Therefore we get the identifications 
\[
\rho_r\colon\Q\,\zf_{2r+1} \otimes \gr_{\ell-1}^F (\ZfH)_N \overset{\sim}{\longrightarrow} \gr_{\ell-1}^F (\ZfH)_N
\]
induced by $\zf_{2r+1} \mapsto 1$ for all $r,\ell,N\in\N$. 

\begin{Definition}\label{def:zf_DleN_maps}
    Assuming the \ogccomma we define for all $N,\ell\geq 1$ the maps
	\begin{equation*}
		D_{<N}^{(\ell)} \colon \gr_\ell^F (\ZfH)_N  \rightarrow \bigoplus_{3\leq2r+1\leq N} \gr_{\ell-1}^F (\ZfH)_{N-2r-1} 
	\end{equation*}
	via
	\begin{equation*}
		D_{<N}^{(\ell)} = \bigoplus_{3\leq 2r+1\leq N} \rho_r \circ  D_{2r+1}^{(\ell)}.
	\end{equation*}
\end{Definition}

Observe that 
\begin{equation*}
	\bigoplus_{3\leq2r+1\leq N} \gr_{\ell-1}^F (\ZfH)_{N-2r-1} = \gr_{\ell-1}^F (\ZfH)_{<N-1}
\end{equation*}
since the weight gives a grading on $\Zf$.

\begin{Lemma} \label{lem:diag_level-reduction_B23_Zf}
If \ogc \ref{conj:odd_generators} holds, then
the  diagram
\begin{equation}\label{eq:diag_level-reduction_B23_Zf}
\begin{aligned}
\xymatrix{
\gr_\ell^F (\B)_N \ar[d]^{\zf} \ar[rr]^-{\partialphi}  && \bigoplus\limits_{3\leq2r+1\leq N}  \gr_{\ell-1}^F (\B) _{N-2r-1} \ar[d]^{\bigoplus\limits_{3\leq 2r+1\leq N} \zf} \\
\gr_\ell^F (\ZfH)_N  \ar[rr]^{\hspace{-3em}D_{<N}^{(\ell)}}   && \bigoplus\limits_{3\leq2r+1\leq N} \gr_{\ell-1}^F (\ZfH)_{N-2r-1} 
}
\end{aligned}
\end{equation}
is well-defined and commutative.
\end{Lemma}
\begin{proof}
    Recall $\partialphi$ is given in Definition \ref{def:filtered_partial_leN} and by Theorem \ref{thm:matrix_invertible}, the upper row is an isomorphism. The lower horizontal map is given by Definition \ref{def:zf_DleN_maps} and it is well-defined by the discussion before its definition. Since $\partialphi$ and $D_{<N}^{(\ell)}$ are compatible by Remark \ref{rem:diagr_partial_D}, it suffices to show for $w\in\mathcal{B}^1$ that 
    \begin{equation*}
        \phi(w) \zf_{2r+1} = \pi_{2r+1}(\amap(\zf(w))).
    \end{equation*}
    But now, observing that projecting the identities \eqref{eq:palindromic_fmzv} and \eqref{eq:formal_zagier_formula} modulo products and weight gives exactly the coefficients $\phi(w)$ from Definition \ref{def:B_1_to_c_w}.
\end{proof}

\begin{Remark} We may change the assumptions for
   Lemma \ref{lem:diag_level-reduction_B23_Zf} a little bit, namely by imposing  that the Kernel conjecture \ref{conj:kernel} holds and that $\zf_{2r+1}\neq 0\in \Lf_{2r+1}$ for all $r\in\N$.
\end{Remark}

\begin{Theorem}\label{thm:linearly_independence_fmzv}
	Assume the \ogc \ref{conj:odd_generators}.
	Then the elements
	\begin{equation*}
		\{\zf(k_1,\ldots,k_d) \mid k_i\in\{2,3\}\}
	\end{equation*}
	are linearly independent.
\end{Theorem}

\begin{proof}
    Recall that diagram \eqref{eq:diag_level-reduction_B23_Zf} commutes by Lemma \ref{lem:diag_level-reduction_B23_Zf} because of our assumption of the \ogcdot By Theorem \ref{thm:matrix_invertible}, the upper row is an isomorphism and obviously the vertical maps are epimorphisms. We prove by induction on the level $\ell\in\N$ that they are in fact isomorphisms. \\
    We start with the case $\ell=1$ and $N\in\N$ arbitrary. Since $\Zf$ is a weight-graded algebra, the elements $\zf(\{2\}^n)$ are linearly independent for all $n\geq0$. We deduce that the right-hand side of \eqref{eq:diag_level-reduction_B23_Zf} is an isomorphism. Thus the left-hand side must also be an isomorphism. \\
    Clearly, the same lines of argument would hold for all $\ell>1$, thus \eqref{eq:diag_level-reduction_B23_Zf} is a commutative diagram of isomorphisms. Now, there are no relations in different weights and if there were a relation in different level, then it must induce a relation in the corresponding graded pieces. But we have shown that this is impossible, thus the claim follows.
\end{proof}

\begin{Theorem}
    Assume the \ogc \ref{conj:odd_generators}, then
    \[
        \Zf = \Zf_{2,3}.
    \]
\end{Theorem}

\begin{proof}
Under the assumption of the \ogccomma Theorem \ref{thm:linearly_independence_fmzv} implies that $\{\zeta(k_1,\ldots,k_d)\mid k_i\in \{2,3\}\}$ forms a basis of $\ZfH$. Counting the number of such indices in a given weight $N$ yields
 \[ \sum_{N\geq 0} \dim\big((\Zf_{2,3})_N\big)\, x^N = \frac{1}{1-x^2-x^3}. \]
 This agrees with the computed dimensions for $\Zf$ in Proposition \ref{prop:ogc_implications}. We deduce that the two spaces $\Zf$ and $\ZfH$ must agree.
\end{proof}

\begin{Remark} Finally we observe that Brown's proof for the generators also works in the formal setup:  Assume the \ogc holds, then
any element in $\Zf$ can be written uniquely as a polynomial of 
elements $ \zf(w)$, where $w$  is a Lyndon word in $\two$ and $\three$. 
\end{Remark}

\cleardoublepage

\phantomsection

\addcontentsline{toc}{section}{References}

\bibliographystyle{amsalpha}

\bibliography{agzt_lectures}
 
\end{document}